\documentclass[final,3p,sort&compress]{elsarticle}

\usepackage{amssymb}
\usepackage{amsthm}
\usepackage{amsmath}
\usepackage{mathtools}
\usepackage{subfig}
\usepackage{showkeys}
\usepackage{booktabs}

\usepackage{tikz}
\usepackage{pgfplots}
\usetikzlibrary{calc,patterns}
\usetikzlibrary{external}
\tikzsetfigurename{figure_}

\newsavebox{\tempbox}

\def \u{{\bf u}}
\def \uo{{\overline{\bf u}}}
\def \po{{\overline{p}}}
\def \vo{{\overline{\bf v}}}
\def \qo{{\overline{q}}}
\def \v{{\bf v}}
\def \w{{\bf w}}
\def \x{{\bf x}}
\def \vv{{\bf V}}

\def \ff{{\bf f}}

\def\resm{\mathbf{R}^{\mathrm{M}}_h}
\def\resc{\mathrm{R}^{\mathrm{C}}_h}

\newcommand{\btau}{{\boldsymbol \tau}}

\begin{document}

\begin{frontmatter}

\title{Numerical comparisons of finite element stabilized methods for\\ high Reynolds numbers vortex dynamics simulations}

\author[wa]{Naveed Ahmed} \ead{naveed.ahmed@wias-berlin.de}
\author[us]{Samuele Rubino} \ead{samuele@us.es}

\cortext[cor]{Corresponding author}
\address[wa]{
Weierstrass Institute for Applied Analysis and Stochastics (WIAS),
Mohrenstr. 39, 10117 Berlin, Germany 
}
\address[us]{Department EDAN \& IMUS, University of Seville, Avda. Reina Mercedes s/n, 41012 Seville, Spain
}

\begin{abstract}
In this paper, we consider up-to-date and classical Finite Element (FE) stabilized methods 
for time-dependent incompressible flows. All studied methods belong to the Variational 
MultiScale (VMS) framework. So, different realizations of stabilized FE-VMS methods are 
compared in high Reynolds numbers vortex dynamics simulations. In particular, a fully 
Residual-Based (RB)-VMS method is compared with the classical Streamline-Upwind 
Petrov--Galerkin (SUPG) method together with grad-div stabilization, a standard 
one-level Local Projection Stabilization (LPS) method, and a recently proposed LPS method 
by interpolation. These procedures do not make use of the statistical theory of equilibrium 
turbulence, and no ad-hoc eddy viscosity modeling is required for all methods. Applications 
to the simulations of high Reynolds numbers flows with vortical structures on relatively 
coarse grids are showcased, by focusing on two-dimensional plane mixing-layer flows. Both 
Inf-Sup Stable (ISS) and Equal Order (EO) FE pairs are explored, using a second-order 
semi-implicit Backward Differentiation Formula (BDF2) in time. Based on the numerical 
studies, it is concluded that the SUPG method using both ISS and EO FE pairs performs best 
among all methods. Furthermore, there seems to be no reason to extend SUPG method by the 
higher order terms of the RB-VMS method. 
\end{abstract}
\end{frontmatter}

\section{Introduction} \label{sec;intro}
In this paper, we consider up-to-date and classical Finite Element (FE) stabilized methods for time-dependent incompressible flows fulfilling the 
incompressible Navier--Stokes Equations (NSE). Let \(\Omega \in \mathbb{R}^d\), \(d \in \{2, 3\}\), be a bounded domain with Lipschitz boundary \(\Gamma\)
and \((0,T)\) be a bounded time interval. 
The incompressible NSE read 
as follows:  \medskip

Find a velocity field \(\u: (0,T] \times \Omega \rightarrow \mathbb{R}^d \) and a pressure field 
\(p: (0,T] \times \Omega \rightarrow \mathbb{R} \) such that

\begin{equation}\label{eq:nse}
\begin{array}{rcll}
\partial_t \u -\nu\Delta \u+ (\u\cdot\nabla)\u +\nabla p &=&
\ff & \mbox{in } (0,T] \times \Omega, \\
\nabla \cdot \u &=&0 & \mbox{in } [0,T] \times \Omega,\\
\u&=&{\bf 0} & \mbox{on } [0,T] \times \Gamma,\\
\u(0,{\bf x}) &=& \u_0 & \mbox{in } \Omega,
\end{array}
\end{equation}
where \(\nu\) is the kinematic viscosity that is assumed to be positive and constant,
\(\ff \) is the given body force, and \(\u_0 \)
is the given initial velocity field, assumed to be divergence-free. 
For simplicity of presentation, we consider the case of homogeneous Dirichlet 
boundary conditions on the whole boundary. 

The main contribution of this paper is a comprehensive and thorough numerical study in the FE stabilized framework of two-scales 
fully Residual-Based (RB) and local projection-based Variational MultiScale (VMS) methods for time-dependent 
high Reynolds numbers incompressible flows with a strong dynamic vortical structure. The derivation of efficient and accurate numerical schemes 
for the simulation of turbulent incompressible flows is a very active field of research.
In particular, various realizations of VMS methods for simulating turbulent incompressible flows have been proposed in the past fifteen years 
(see~\cite{ACJR17} for a recent detailed review). All of these realizations obey the basic principles of VMS methods: they are based on the 
variational formulation of the incompressible NSE and the scale separation is defined by projections. However, apart from these common basic 
features, the various VMS methods look quite different. In this paper, our main goal is to focus on two-scales VMS methods, and provide a thorough 
numerical investigation of up-to-date and classical FE stabilized methods belonging to this category when applied to a relevant fixed setup for 
numerical studies such as the 2D Kelvin--Helmholtz instability problem. Indeed, even if VMS methods, despite their relatively recent development, 
are already well-established, and considered state-of-the-art in turbulence modeling that provides a promising and successful alternative to classical 
Large Eddy Simulation (LES) models, in the literature there is no so much about a structured comparison of them in terms of numerical studies. Up to 
our knowledge, the first (and only) attempt to go towards this research direction has been performed in~\cite{JK10_cmame}, where the authors studied
different realizations of VMS methods within the framework of FE in turbulent channel flow simulations. However, they just focus on three-scales VMS 
models, in which the effect of the unresolved scales on the resolved ones is modeled by means of an eddy viscosity term of Smagorinsky type that only
acts directly on the small resolved scales. In the present paper, we aim at complementing and extending this research avenue, by mainly focusing on 
two-scale VMS methods, which use a direct modeling of the subgrid scale flow by numerically approximating the related equations. Thus, they do not 
need any modeling of the subgrid scales by statistical theories of turbulence, and in particular they do not include eddy viscosity. The numerical 
diffusion inherent to those stabilized models basically plays the role of the eddy diffusion. 
In this way, the present paper aims at giving a thorough numerical investigation, similar to the one performed in~\cite{JK10_cmame}, but for two-scale 
VMS methods. A structured presentation is provided in this framework, with special emphasis on experience in numerical studies.
Once reached almost ``definitive'' conclusions within this paper, a comparison of the selected ``best performing'' two-scale VMS method with three-scale 
VMS methods that use eddy viscosity (in a more or less sophisticated manner) to model the effect of subgrid scales shall appear in a forthcoming paper. 
In this way, the numerical performances of different VMS methods would be assessed. Up to our knowledge, this is the first time that such a numerical 
study is conducted in the literature in a unified VMS framework.

The RB-VMS method was introduced in~\cite{BCCHRS07}. A straightforward simplification 
of the RB-VMS method leads to the classical Streamline-Upwind Petrov--Galerkin (SUPG) method~\cite{BH82, HB79}.
Also, another variant of the RB-VMS method, which is not fully consistent, but of optimal order with respect to
the FE interpolation, is given by the so-called Local Projection Stabilization (LPS) methods~\cite{beckerbraack}.
So, different realizations of stabilized FE-VMS methods are compared in high Reynolds numbers vortex dynamics 
simulations in this paper. In particular, the RB-VMS method~\cite{BCCHRS07} is compared with the classical 
SUPG method~\cite{BH82, HB79} together with grad-div stabilization, a standard one-level LPS 
method~\cite{matskytob}, and a recently proposed LPS method by interpolation~\cite{chamavisa, Naveed15}. 
To our best knowledge, a comparison of these methods is so far not available.
To keep the paper self-contained, a brief presentation of the cited numerical methods, which provides the basic 
concepts, will be considered hereafter. For more details on their derivation, see the up-to-date review on VMS 
methods for the simulation of turbulent incompressible flows~\cite{ACJR17}.

To assess the different numerical methods, applications to the simulations 
of high Reynolds numbers flows with vortical structures on relatively coarse grids are showcased, by 
focusing on two-dimensional plane mixing-layer flows as benchmark problem, since it presents a wide range of flow scales and an interesting 
time evolution of the flow field. Starting from a perturbed initial condition, the transition to the development of small vortices takes place, 
which then pair to larger vortices until one single eddy finally remains, rotating at a fixed position.
In particular, we analyze different quantities of interest associated to this problem (i.e, temporal evolution of vorticity field/thickness, 
kinetic energy, enstrophy, palinstrophy) in order to judge the performance of all studied methods, and draw some definitive conclusions. All 
the numerical results are benchmarked against a reference simulation, consisting of a Pressure Stabilized Petrov--Galerkin (PSPG) 
method~\cite{HFB86} with grad-div stabilization, computed with a finer space resolution. However, note that this model problem is 
very sensitive to small perturbations that are almost unavoidable in numerical simulations, thus some targets, such as a conclusive 
prediction of the final pairing into one single eddy, seems to be somehow not achievable, even among the simulations with the highest resolutions.

In this paper, we mainly focus on efficient spatial and temporal discretizations, for which both Inf-Sup Stable (ISS) and Equal Order (EO) low-order 
FE pairs are explored, using a second-order semi-implicit Backward Differentiation Formula (BDF2) in time, where the linearization of the fully 
discrete problem at each time step is done by means of temporal extrapolation. In contrast to a fully implicit scheme, this approach yields a unique 
linear system of equations to be solved at each time step. Altogether, performing simulations with semi-implicit schemes uses less computing time than 
fully implicit schemes. However, while a fully implicit approach is generally yielding a stable time discretization scheme, a semi-implicit approach 
may require a time step restriction due to the stability issue of the time stepping scheme. For this reason, we performed simulations using larger and 
smaller time step lengths, and evaluating the corresponding numerical results, we noticed effectively that physical consistency (e.g., monotone 
decline of kinetic energy) is lost to some extent when considering larger time step for some methods. Note that semi-implicit BDF schemes for the 
numerical simulation of NSE with VMS turbulence modeling have already been investigated in the literature, see for instance~\cite{FortiDede15}, 
and also~\cite{Rubino18} for a stable velocity-pressure segregation version. 

The paper is organized as follows. In Section \ref{sec:VMS}, VMS methods are described, with a special focus on the derivation of two-scale VMS 
methods analyzed in the present work, that are RB-VMS methods and LPS methods. In Section \ref{sec:TimeDisc}, a semi-implicit approach for the 
time discretization, applying the two-step BDF (BDF2) in order to get the corresponding fully discrete schemes, is detailed for each studied method, 
together with some numerical implementation aspects. In Section \ref{sec:numerics}, the studied methods are numerically compared on the simulation 
of two-dimensional Kelvin--Helmholtz instabilities in the high Reynolds number regime. Here, several quantities of interest are presented, evaluated 
and discussed. Finally, Section \ref{sec:conclusions} summarizes the main conclusions of the paper and gives an outlook.

\section{Variational multiscale methods}\label{sec:VMS}
As already mentioned, VMS methods are based on the
variational formulation of the incompressible NSE \eqref{eq:nse}.
To define the variational formulation of \eqref{eq:nse}, the velocity space \(\vv = [H^1_0(\Omega)]^d\) and the 
pressure space \(Q = L^2_0(\Omega)\) are introduced. Let \((\cdot,\cdot)\) denote the \(L^2\) inner product with respect 
to the domain \(\Omega\).
The variational formulation of \eqref{eq:nse} reads as follows: \medskip

Find \((\u,p) : (0,T) \rightarrow \vv \times Q \) such that for all \((\v,q) \in \vv\times Q \)
\begin{equation}\label{eq:weak_form}
\frac{d}{dt}(\u,\v )+\nu (\nabla\u, \nabla\v) + ((\u\cdot\nabla)\u, \v) - (p, \nabla \cdot \v)  + 
(\nabla \cdot \u,q)=
\langle \ff,\v\rangle \qquad \mbox{in } \mathcal{D}^\star(0,T), 
\end{equation}
with \(\u(0,{\bf x}) = \u_0({\bf x}) \) in \(\Omega\), where \(\langle \cdot, \cdot \rangle\) denotes the duality pairing 
between the velocity space \(\vv \) and its dual \(\vv^\star\) and \(\mathcal{D}^\star(0,T) \) is the space 
of distribution on \((0,T)\).

In standard conforming Finite Element (FE) formulations, the infinite-dimensional spaces \((\vv , Q)\) are replaced with finite 
dimensional-subspaces \((\vv_h, Q_h)\) consisting of typically low-order piecewise polynomials with respect to a triangulation 
\(\mathcal{T}_h\) of \(\Omega\). In this paper, both Inf-Sup Stable (ISS, \cite{Bab71, Bre74}) and Equal Order (EO) conforming FE 
pairs are explored, which are not exactly divergence-free, by considering in general the popular Taylor--Hood FE pair 
\(\mathbf{P}_{k}/\mathbb{P}_{k-1}\) \cite{HT74} and the EO FE pair \(\mathbf{P}_{k}/\mathbb{P}_{k}\), respectively, 
with \(k\geq 2\), where \(\mathbb{P}_k\)  denotes the space of continuous functions whose restriction to each mesh 
cell \(K\in \mathcal{T}_h\) is the Lagrange polynomial of degree less than or equal to \(k\), and $\mathbf{P}_{k}=[\mathbb{P}_{k}]^d$.

\subsection{Two-scale VMS methods}
This section discusses basic concepts of two-scale VMS methods.
Starting point of two-scale VMS methods is the separation of the
flow field into resolved scales \((\uo,\po)\) and unresolved scales 
\((\u', p')\) such that \(\u = \uo+\u'\) and  \( p = \po + p^\prime\). 
Analogously, a direct-sum decomposition of velocity space
\(\vv = \overline{\vv} \oplus \vv^\prime \) and pressure space \(Q =\overline{Q} \oplus Q^\prime\)
is considered. 
It should be
emphasized that although this approach is in
principle the same as in Large Eddy Simulations (LES), it is well known that the definition of the scales is
different. A variational projection, either \(L^2\) projection or elliptic projection, for the separation of scales and spaces is performed in VMS methods.

Note that the VMS methodology allows further
decompositions of the resolved scales. The most common
approach of this kind is a decomposition of these scales
into large resolved scales (or large scales) and small
resolved scales, leading finally to a so-called three-scale
VMS method. In this case, the effect of the unresolved scales on
the resolved ones is modeled by means of an eddy
viscosity term that only acts directly on the small
resolved scales (cf. \cite{GWR04, GWR05, JKK08}). However, in the present paper we just focus on the comparison between
VMS methods that use a direct modeling of the subgrid scale flow by approximating the related equations,
for which no eddy viscosity is introduced to model the effect of the subgrid scales. This is the reason why we restrict to two-scale VMS 
methods. Once reached almost ``definitive'" conclusions within this paper, a comparison of the selected ``best performing" two-scale VMS 
method with three-scale VMS methods that use eddy viscosity (in a more or less sophisticated manner) to model the effect of subgrid scales 
shall appear in a forthcoming paper.  

For clearness of presentation, the weak formulation \eqref{eq:weak_form} of the NSE
will be expressed in a short form as follows:\medskip

Given \(\u(0,{\bf x})=\u_0({\bf x})\), find \((\u,p):(0,T)\rightarrow \vv \times Q\) satisfying 
\begin{equation}\label{eq:weak_form_n}
A\left(\u;(\u,p),(\v,q)\right) = \ff(\v) \quad \forall\ (\v,q)\in \vv\times Q.
\end{equation}

Decomposing also the test functions into two scales and using the linearity with respect to the test functions, 
the variational formulation \eqref{eq:weak_form_n} leads to the coupled set of equations: 
\begin{itemize}
	\item an equation for the resolved scales
	\begin{equation}\label{eq:resolved_scales}
	A\left(\u;\left(\uo,\po\right),\left(\vo,\qo\right)\right)  +
	A\left(\u;\left(\u',p'\right),\left(\vo,\qo\right)\right) 
	= \ff\left(\vo\right), 
	\end{equation}
	\item and an equation for the unresolved scales
	\begin{equation}
	\label{eq:unresolved_scales}
	A\left(\u;\left(\uo,\po\right),\left(\v',q'\right)\right)  +
	A\left(\u;\left(\u',p'\right),\left(\v',q'\right)\right) 
	=  \ff\left(\v'\right). 
	\end{equation}
\end{itemize}
The form \(A(\cdot;\cdot,\cdot)\) is decomposed into its linear part and the trilinear convective term as
\begin{equation*}
A\left(\u;{\bf U},{\bf W}\right) = A_{\mathrm{lin}}\left({\bf U},{\bf W}\right) + ((\u\cdot \nabla \u),\v)
\end{equation*}
where the abbreviations \({\bf U} = (\u,p)^T \) and \({\bf W} = (\v, q)^T\) are used for simplicity. Then, the equation
\eqref{eq:unresolved_scales} for the unresolved scales can be written in the form
\begin{equation}\label{eq:basic_small_0}
A_{\bf U}\left({\bf U}',{\bf W}'\right) + \left((\u'\cdot \nabla)\u',\v'\right) = \left\langle \mathbf{R}\left(\overline{\bf U}\right),{\bf W}'
\right\rangle
\end{equation}
with 
\begin{align*}
A_{\bf U}\left({\bf U}',{\bf W}'\right) &=  A_{\mathrm{lin}}\left({\bf U}',{\bf V}'\right) + \left((\u'\cdot \nabla) \uo,\v'\right) 
+ \left((\uo\cdot \nabla )\u',\v'\right),\\
\left\langle \mathbf{R}\left(\overline{\bf U}\right),{\bf V}'
\right\rangle & =  \ff(\v') 
- A_{\mathrm{lin}}\left(\overline{\bf U},{\bf W}'\right) - \left((\uo\cdot \nabla \uo),\v'\right),
\end{align*}
where \(A_{\bf U}\left({\bf U}',{\bf W}'\right)\) is the G\^ateaux derivative
of $A(\cdot;\cdot,\cdot)$ at ${\bf U}$ in the direction of \({\bf U}'\). 
The solution of \eqref{eq:basic_small_0} can be formally represented as
\begin{equation}\label{eq:small_scales}
{\bf U}' = F_{\bf U}\left(\mathbf{R}\left(\overline{\bf U}\right)\right),
\end{equation}
which can be interpreted as the unresolved scales which are driven as a function of the residual of the resolved scales. 
Finally, inserting expression \eqref{eq:small_scales} in the resolved scales equations \eqref{eq:resolved_scales}
leads to a single set of equations for the resolved scales. 

Two-scale VMS methods aim to approximate \(F_{\bf U}\) by models which do not rely on considerations from the physics
of turbulent flows, but are derived just with mathematical arguments. In the next subsections, concrete approaches
will be presented.

\subsection{Residual-based VMS method}
The main idea in the derivation of the two-scale RB-VMS method is based on a perturbation series with respect 
to the norm of the residual associated to the resolved scales. It is proposed in \cite{BCCHRS07} to truncate the series
after the first term and to apply some modeling of this term. The resulting method can be considered as a generalization of classical stabilization methods for the NSE.

A perturbation series for a potentially small quantity 
\(\varepsilon = \|\mathbf{R}(\overline{\bf U})\|_{(\vv^\prime \times Q^\prime)^{*}}\) is considered. It is assumed that the larger the space \((\vv \times Q)\), 
the better \(\overline{\bf U}\) approximates \({\bf U}\), and the smaller is \({\bf R(\overline{\bf U})}\).
The perturbation series is of the form
\begin{equation}\label{eq:pert_series}
{\bf U'} = \varepsilon {\bf U_1'} + \varepsilon^2 {\bf U'_2} + \ldots = \sum_{i=1}^\infty \varepsilon^i  {\bf U'_i}.
\end{equation}
In particular, if \(\varepsilon=0\), i.e. \({\bf R(\overline{\bf U})}=0\), then \({\bf U'}=F_{\bf U}\left(\mathbf{R}\left(\overline{\bf U}\right)\right) = 0\) from \eqref{eq:small_scales}-\eqref{eq:pert_series} .
Inserting the perturbation series \eqref{eq:pert_series} 
in the terms of \eqref{eq:unresolved_scales} for the unresolved scales gives
\begin{equation*}
A_{\bf U}\left( \sum_{i=1}^\infty \varepsilon^i {\bf U}'_i,{\bf W'}\right) =  \sum_{i=1}^\infty \varepsilon^i A_{\bf U}\left( {\bf U}'_i,{\bf W'}\right)
\end{equation*}
and
\begin{align*}
\left( \left(\sum_{i=1}^\infty\varepsilon^i \u'_i\cdot\nabla\right) \sum_{i=1}^\infty\varepsilon^i \u'_i,\v'\right)
&= \varepsilon^2 \left((\u'_1\cdot \nabla)\u'_1,\v'\right) + \varepsilon^3  \left[\left((\u'_1\cdot \nabla)\u'_2,\v'\right) +  \left((\u'_2\cdot \nabla ) \u'_1,\v'\right)\right]
+ \ldots \\
&= \sum_{i=2}^\infty \varepsilon^i \left(\sum_{j=1}^{i-1} \left((\u'_j\cdot \nabla)\u'_{i-j},\v'\right)\right). 
\end{align*}
Substituting these terms into \eqref{eq:unresolved_scales} yields
\begin{equation*}
\sum_{i=1}^\infty \varepsilon^i A_{\bf U}\left( {\bf U}'_i,{\bf W'}\right) + 
\sum_{i=2}^\infty \varepsilon^i \left(\sum_{j=1}^{i-1} \left((\u'_j\cdot \nabla)\u'_{i-j},\v'\right)\right) = 
\varepsilon \left\langle \frac{\mathbf{R}\left(\overline{\bf U}\right)}
{\left\|\mathbf{R}\left(\overline{\bf U}\right)\right\|_{(\vv^\prime \times Q^\prime)^{*}}}, {\bf W}'
\right\rangle.
\end{equation*}
Collecting similar terms with respect to \(\varepsilon\) leads to a system of variational problems which are coupled 
through the right-hand side, that is
\begin{align*}
A_{\bf U}\left( {\bf U}'_1,{\bf W}'\right) &= \left\langle \frac{\mathbf{R}\left(\overline{\bf U}\right)}
{\left\|\mathbf{R}\left(\overline{\bf U}\right)\right\|_{(\vv^\prime \times Q^\prime)^{*}}},{\bf W}'
\right\rangle,\\
A_{\bf U}\left( {\bf U}'_i,{\bf W}'\right) &= -\sum_{j=1}^{i-1} \left((\u'_j\cdot \nabla )\u'_{i-j},\v'\right)  \quad i\ge 2.
\end{align*}
In the modeling of the unresolved scales, it is suggested in \cite{BCCHRS07} to truncate the series \eqref{eq:pert_series} after the first term,
and to use a linear approximation of the so-called fine-scale Green's operator that formally represent \({\bf U}'_1\) 
\begin{align}\label{eq:res_model}
{\bf U}' \approx \varepsilon {\bf U}_1' &= \|\mathbf{R}(\overline{\bf U})\|_{(\vv^\prime \times Q^\prime)^{*}} {\bf U}'_1 
\approx \btau  \mathbf{R}\left(\overline{\bf U}\right) = \btau \mathbf{R}\begin{pmatrix} \u_h\\ p_h\end{pmatrix} \nonumber\\
& =  \begin{pmatrix}
\btau_m\left(\ff_{h} - \partial_t \u_h + \nu \Delta \u_h - (\u_h\cdot \nabla) \u_h - \nabla p_h\right) \\ \\
- \tau_{\rm c} \left(\nabla \cdot \u_h\right) 
\end{pmatrix}
= \begin{pmatrix}
\resm \\ \\
\resc
\end{pmatrix}
\end{align}
where \(\btau\) is a \(4 \times 4 \) diagonal tensor-valued function, and the approximation of the resolved scales is computed in a standard FE space. 

The RB-VMS FE formulation is obtained by inserting the approximation \eqref{eq:res_model}
into the large scales equation \eqref{eq:resolved_scales}, omitting the models of the terms 
\((\partial_t \u' , \v_h) \) and \(\nu (\nabla{\u'},\nabla{\v_h} )\), and integrating by parts the continuity equation with 
respect to the unresolved scale in \eqref{eq:resolved_scales}, assuming that \(\u'=0\) on \(\Gamma\): 

Find \(\u_h \ :\ (0,T)\to\vv_h, \ p_h \ : \ (0,T)\to Q_h\) satisfying
\begin{align}\label{eq:res_vms00}
\left(\partial_t\u_h,\v_h\right) &+ \nu\left(\nabla{\u_h},\nabla{\v_h}\right) + 
\left((\u_h\cdot\nabla)\u_h,\v_h\right) - \left(p_h,\nabla\cdot\v_h\right)
+ \left(\nabla \cdot \u_h,q_h\right)  \nonumber \\
&+b\left(\resm ,\u_h,\v_h\right)
+b\left(\u_h ,\resm,\v_h\right)
+b\left(\resm ,\resm,\v_h\right) \nonumber\\
& - \left(\resc,\nabla \cdot \v_h\right)
- \left(\resm,\nabla q_h\right) 	
= (\ff_h, \v_h)
\end{align}
for all \((\v_h,q_h) \in \vv_h\times Q_h\), where $b$ in \eqref{eq:res_vms00} denotes the trilinear convective form 
given by \(b(\u,\v,\w)=\left((\u\cdot\nabla)\v,\w\right),\quad\u,\v,\w\in \vv\). 

Concerning the actual choice of \(b\), it
is advisable from the practical point of view that one does
not need to compute a derivative of the residual of the
momentum equation. For this reason, it is suggested to use the following form, which is
obtained from the divergence form with integration by parts:
\begin{align}\label{eq:conv}
b(\u, \v, \w) = (\nabla \cdot (\u \v^T), \w) = -(\u \v^T, \nabla \w).
\end{align}
The two terms \(b(\resm, \u_h, \v_h)\) and \(b(\u_h, \resm, \v_h)\) are known 
as cross-stress terms, and \(b\left(\resm ,\resm,\v_h\right)\)
as the subgrid (or Reynolds-stress) term. Using \eqref{eq:conv}, \((\u \v^T, \nabla \w) = (\v, (\nabla \w)^T \u)\) and 
\((\nabla \v)\u=(\u\cdot \nabla )\v \), one gets for the first cross-stress term in \eqref{eq:res_vms00}: 
\begin{align}\label{eq:cross1}
b\left(\resm, \u_h, \v_h\right) = -\left(\resm (\u_h)^T, \nabla \v_h\right) = - \left(\u_h,(\nabla \v_h)^T \resm\right)
= - \left(\resm,(\nabla \v_h)\u_h\right) = -\left( \resm, (\u_h\cdot \nabla )\v_h\right),
\end{align}
which together with the last term in the left-hand side of \eqref{eq:res_vms00} gives:
\begin{align}\label{eq:SUPG}
b\left(\resm, \u_h, \v_h\right) - (\resm , \nabla q_h) = 
-\left( \resm, (\u_h\cdot \nabla )\v_h + \nabla q_h\right).
\end{align}
This term corresponds to the well known stabilization term of the Streamline-Upwind Petrov-Galerkin (SUPG) method for the convection field \(\u_h\). 
One can also observe the contribution of the so-called grad-div stabilization term by inserting the concrete formula of the 
residual of the continuity equation into \eqref{eq:res_vms00}, that is:
\begin{equation}\label{eq:grad-div}
\left(\tau_{\rm c}\nabla \cdot \u_h, \nabla \cdot \v_h\right).
\end{equation}
Similarly, using \eqref{eq:conv} and \((\u \v^T, \nabla \w) = (\v, (\nabla \w)^T \u)\), 
one obtains for the second cross-stress term and the subgrid term in \eqref{eq:res_vms00}: 
\begin{align}\label{eq:cross2}
b\left(\u_h, \resm, \v_h\right) = -\left(\u_h (\resm)^T, \v_h \right) = -\left(\resm, (\nabla \v_h)^T \u_h \right),
\end{align} 
\begin{equation}\label{eq:subgrid}
b(\resm,\resm, \v_h) = -\left(\resm
(\resm)^T, \v_h \right) = -\left(\resm, (\nabla \v_h)^T 
\resm \right).
\end{equation}
Considering formulas \eqref{eq:cross1} and \eqref{eq:cross2} for the cross-stress terms, and formula \eqref{eq:subgrid} for the subgrid term, the RB-VMS method \eqref{eq:res_vms00} can be 
expressed as: 

Find \(\u_h \ :\ (0,T)\to\vv_h, \ p_h \ : \ (0,T)\to Q_h\) satisfying
\begin{align}\label{eq:res_vms}
\left(\partial_t\u_h,\v_h\right) &+ \nu\left(\nabla{\u_h},\nabla{\v_h}\right) + 
\left((\u_h\cdot \nabla )\u_h,\v_h\right) -	\left(p_h,\nabla \v_h\right)+ \left(\nabla \cdot \u_h,q_h\right)  
- \left(\resm ,(\u_h \cdot \nabla )\v_h + C\nabla q_h\right) \nonumber \\
& - \left( \resm, (\nabla \v_h)^T \u_h\right)
- \left(\resm ,(\nabla \v_h)^T \resm\right)
+ \left(\tau_{\rm c}\nabla \cdot \u_h ,\nabla \cdot \v_h\right)	
= (\ff_h, \v_h),
\end{align}
for all \((\v_h,q_h) \in \vv_h\times Q_h\).
The formulation \eqref{eq:res_vms} provides the complete RB-VMS method, which retains numerical consistency in the FE equations, in the
sense that the continuous solution exactly satisfies the discrete equations, whenever it is smooth enough. 
In this paper, both ISS and EO conforming FE pairs would be explored. For this reason, we have added the constant \(C\) in formulation
\eqref{eq:res_vms}, so that \(C=1\) when using EO FE pairs, and we will drop the 
dependency of the pressure stabilization term from \eqref{eq:res_vms} when using ISS FE pairs by fixing \(C=0\). 
We recall that in \eqref{eq:res_vms}
the terms
\[
\left(\resm ,(\u_h \cdot \nabla )\v_h + \nabla q_h\right) \qquad \text{ and } \qquad
\tau_{\rm c} \left(\nabla \cdot \u_h ,\nabla \cdot \v_h\right)
\]
are the classical stabilization terms of the SUPG and grad-div methods, respectively. In this paper, we are interested in performing 
numerical studies also with a simplified model arising from \eqref{eq:res_vms}, which is the classical SUPG method together with 
grad-div stabilization:

Find \(\u_h \ :\ (0,T)\to\vv_h, \ p_h \ : \ (0,T)\to Q_h\) satisfying
\begin{align}\label{eq:SUPGgd}
\left(\partial_t\u_h,\v_h\right) &+ \nu\left(\nabla{\u_h},\nabla{\v_h}\right) + 
\left((\u_h\cdot \nabla )\u_h,\v_h\right) -	\left(p_h,\nabla \v_h\right)+ \left(\nabla \cdot \u_h,q_h\right)  
- \left(\resm ,(\u_h \cdot \nabla )\v_h + C\nabla q_h\right) \nonumber \\
& + \left(\tau_{\rm c}\nabla \cdot \u_h ,\nabla \cdot \v_h\right)	
= (\ff_h, \v_h),
\end{align}
for all \((\v_h,q_h) \in \vv_h\times Q_h\), again for both ISS (\(C=0\)) and EO (\(C=1\)) FE pairs. 

\subsection{Local projection stabilization methods}
Local Projection Stabilization (LPS) methods are 
stabilization methods that provide specific stabilization of any single operator term that could be a
source of instability for the numerical discretization. 
They were introduced in  \cite{beckerbraack} and they could be viewed as simplifications of the two-scale RB-VMS
method described in the previous section.
Indeed, LPS methods are not fully consistent (only specific dissipative interactions are retained), but of optimal order with 
respect to the FE interpolation. 
The fact
that the stabilization enjoys the right asymptotic behavior
without full consistency allows to decouple the stabilization
of the pressure and the velocity, without having all the residual terms coupled, thus relying on a term-by-term structure. 
This feature could be
considered an important advantage with respect to the more
complex RB-VMS method in view of practical
implementations such as to perform the numerical analysis,
since it leads to a simpler and less expensive structure.               
Different variants of LPS methods have been investigated during the recent years
for incompressible flow problems. The main common feature is that, thanks to local projection,
the symmetric stabilization terms only act on the small scales of the flow, thus
ensuring a higher accuracy with respect to more classical stabilization procedures, such
as penalty-stabilized methods, cf. \cite{Chacon98}. Thus, the effect of LPS is on the one hand to improve the convergence to
smooth solutions. On the other hand, for rough solutions, LPS limits the propagation of
perturbations generated in the vicinity of sharp gradients, potentially maintaining these
schemes as suitable and useful tools for the simulation of turbulent flows.

As a single rule, the structure of LPS method is achieved by considering in the RB-VMS method \eqref{eq:res_vms} just the specific dissipative 
interactions that stabilize convection and pressure gradient, and by introducing local \(L^2\) projections in the approximation of the 
unresolved scales, in such a way the symmetric stabilization terms only act on the small scales of the flow. 
This leads to a family of methods, associated to the choice of the actual local \(L^2\) projection. 

The main derivation of LPS methods will be introduced here for the NSE \eqref{eq:nse}. The stabilization effect is achieved by adding least-squares 
terms that give a weighted control on the fluctuations of the quantity of interest. This control is based upon a projection operation 
\(\pi_h :\,L^2(\Omega) \mapsto  D_h\) onto a discontinuous FE space \(D_h\) (the \lq\lq projection\rq\rq space). This space is built 
on a grid \({\cal M}_h\) formed by macro-elements built from the triangulation \({\cal T}_h\) of \(\Omega\). The component-wise extension 
of \(\pi_h\) to vector functions is 
denoted by \(\boldsymbol{\pi}_h\).  The LPS approximation of the NSE reads: 

Find \(\u_h \ :\ (0,T)\to\vv_h, \ p_h \ : \ (0,T)\to Q_h\) satisfying
\begin{align}\label{eq:lps}
\left(\partial_t\u_h,\v_h\right) &+ \nu\left(\nabla{\u_h},\nabla{\v_h}\right) + 
\left((\u_h\cdot \nabla )\u_h,\v_h\right) -	\left(p_h,\nabla \v_h\right)+ \left(\nabla \cdot \u_h,q_h\right) \nonumber \\
& + \left(\btau_{m}\boldsymbol{k}_{h}((\u_h \cdot \nabla) \u_h) ,\boldsymbol{k}_{h}((\u_h \cdot \nabla) \v_h)\right) 
+ \left(\btau_{m}\boldsymbol{k}_{h}(\nabla p_h),\boldsymbol{k}_{h}(C\nabla q_h)\right)
+ \left(\tau_{\rm c}\nabla \cdot \u_h ,\nabla \cdot \v_h\right)	
= (\ff_h, \v_h),
\end{align}
for all \((\v_h,q_h) \in \vv_h\times Q_h\).
In \eqref{eq:lps}, \(\boldsymbol{k}_{h} = \boldsymbol{I} - \boldsymbol{\pi}_h\) is the \lq\lq fluctuation\rq\rq operator, being \(\boldsymbol{I}\) 
the identity operator. Also, the additional grad-div term stabilizing term has been added, since not exactly divergence-free FE pairs would be explored. 
As before, we have added the constant \(C\) in formulation \eqref{eq:lps}, so that \(C=1\) when using EO FE pairs, and we will drop the 
dependency of the pressure stabilization term from \eqref{eq:lps} when using ISS FE pairs by fixing \(C=0\). 

The stability of LPS methods is based upon local inf-sup conditions (see \cite{ACJR17}, Section 6.2): The local restriction \(\vv_h(M)\) of the velocity 
space \(\vv_h\) (the \lq\lq approximation\rq\rq space) to any macro-element \(M \in {\cal M}_h\) must be rich enough in degrees of freedom with respect 
to \(D_h(M)\), much as in mixed methods the global velocity space \(\vv_h\) must be rich enough with respect to the pressure space \(Q_h\) to achieve 
the standard discrete inf-sup condition \cite{Bab71, Bre74}. 
With this purpose, two main approaches of LPS methods have been proposed
(see \cite{hetobiska}): In the one-level approach, the approximation space is enriched such that the local inf-sup condition holds and both \(\vv_h\) 
and \(D_h\) are built on the same mesh. In the two-level approach, the projection space is built on a coarser mesh level to satisfy the local inf-sup 
condition. It is possible to consider overlapping sets of macro elements (see \cite{beckerbraack2}). In this work, we will restrict numerical studies 
to the one-level LPS method (defined on a single mesh), considering \(\mathbf{P}_{2}^{\rm bubble}/\mathbb{P}_{1}^{\rm dc}\) ISS FE pair on the one hand, 
and \(\mathbf{P}_{2}^{\rm bubble}/\mathbb{P}_{2}^{\rm bubble}\) EO FE pair on the other hand, with projection space \(D_{h}=\mathbb{P}_{1}^{\rm dc}\), 
i.e. the discontinuous version of \(\mathbb{P}_{1}\).

\subsubsection{Local projection stabilization by interpolation}
A further simplification of LPS schemes is achieved when the local \(L^2\) projection operator \(\boldsymbol{\pi}_h\) is replaced by an interpolation 
operator from \([L^2(\Omega)]^d\) onto a projection space \({\bf D}_h\) formed by continuous FE (see \cite{chamavisa}). To describe this approach, 
assume that the discrete velocity and
pressure spaces \(\vv_h\) and \(Q_h\) are formed by piecewise polynomial functions of degree \(k\) at most, e.g.
\begin{equation}
\label{ded}
\vv_h = \mathbf{P}_{k} \cap \vv, \quad Q_h=\mathbb{P}_{k}\cap Q.  
\end{equation}
It is assumed that \(\boldsymbol{\pi}_h\) is some locally stable approximation operator from \([L^2(\Omega)]^d\) onto \({\bf D}_h=\mathbf{P}_{k-1}\), 
satisfying optimal error estimates. 
In practical implementations, we choose \(\boldsymbol{\pi}_h\) as a Scott--Zhang-like \cite{SZ90} linear interpolation operator in the space 
\(\mathbf{P}_{1}\) (since we consider \(\mathbf{P}_{2}\) as FE velocity space), implemented in the software FreeFem++ \cite{Hecht12}. This 
interpolant may be defined as
\begin{equation*}
\forall \x\in\overline{\Omega},\quad \boldsymbol{\pi}_{h}(\v)(\x)=\sum_{a\in{\cal N}}\Pi_{h}(\v)(a){\boldsymbol{\varphi}}_{a}(\x),
\end{equation*}
where \({\cal N}\) is the set of Lagrange interpolation nodes of \(\mathbf{P}_{1}\), \({\boldsymbol{\varphi}}_{a}\) are the Lagrange 
basis functions associated to \({\cal N}\), and \(\Pi_{h}\) is the interpolation operator by local averaging of Scott--Zhang kind, 
which coincides with the standard nodal Lagrange interpolant when acting on continuous functions (cf.~\cite{chamavisa}, section 4).
This is an interpolant that just uses nodal values, and so is simpler to work out and more computationally efficient than the variant 
of the Scott--Zhang operator introduced in \cite{Badia12} for the Stokes problem, which is instead an operator defined from a 
node-to-element map and requires integration on mesh elements. The LPS method by interpolation is still stated by \eqref{eq:lps}, 
but assuming that the grids \({\cal T}_h\) and \({\cal M}_h\) coincide. The stability of this LPS method by interpolation follows 
from a specific discrete inf-sup condition (see \cite{Naveed15}, Lemma 4.2).

Therefore, this method presents the same structure of the Streamline Derivative-based (SD-based) LPS model 
\cite{BraackBurman06, KnoblochLube09}, but it differs from it because at the same time it uses continuous 
buffer functions, it does not need enriched FE spaces, it does 
not need element-wise projections satisfying suitable orthogonality properties, and it does not need 
different nested meshes. An interpolant-stabilized structure of Scott--Zhang type replaces the projection-stabilized structure of 
standard LPS methods. The interpolation operator takes its values in a continuous buffer space, different from 
the discrete velocity space, but defined on the same mesh, constituted by standard polynomials with one degree 
less than the FE space for the velocity. This approach gives rise to a method with reduced computational 
cost for some choices of the interpolation operator. This method has been recently supported by
a thorough numerical analysis (existence and uniqueness, stability, convergence, error estimates, asymptotic
energy balance) for the nonlinear problem related to the evolution NSE, cf.~\cite{Naveed15}, using a 
semi-implicit Euler scheme for the monolithic discretization in time. In particular, the error analysis 
reveals a self-adapting high spatial accuracy in laminar regions of a turbulent flow that turns to be of
overall optimal high accuracy if the flow is fully laminar. Numerical simulations
of 3D Beltrami flow in laminar regimes \cite{Naveed15} confirm this fact. This also allows to
obtain an asymptotic energy balance for smooth flows.

\section{Time discretization and numerical implementation aspects}\label{sec:TimeDisc}
In this section, we propose a semi-implicit approach for the time discretization, applying the two-step backward 
difference formula (BDF2) in order to get the fully discrete schemes. 
We compute the appro\-ximations \(\u_{h}^{n}\) and \(p_{h}^{n}\) to \(\u^{n}=\u(\cdot,t_{n})\) and \(p^{n}=p(\cdot,t_{n})\), 
respectively, by using temporal schemes based on semi-implicit BDF2, for which the nonlinear terms are extrapolated by means 
of Newton--Gregory backward polynomials \cite{Cellier91}. In order to abbreviate the discrete time derivative, we define the 
operator \(D_{t}^{2}\) by
\begin{equation}\label{eq:discTimeDer}
D_{t}^{2}\u_{h}^{n+1}=\frac{3\u_{h}^{n+1}-4\u_{h}^{n}+\u_{h}^{n-1}}{2\Delta t},\quad n\geq 1.
\end{equation}
We consider the following extrapolation for the convection velocity: \(\widehat{\u}_{h}^{n}=2\u_{h}^{n}-\u_{h}^{n-1}, n\geq 1\), in 
order to achieve a second-order accuracy in time for all methods. For the initialization (\(n=0\)), we consider \(\u_{h}^{-1}=\u_{h}^{0}\), 
being \(\u_{h}^{0}\) the initial condition, so that time schemes reduce to semi-implicit Euler method for the first time step 
\((\Delta t)^{0}=(2/3)\Delta t\). 

For all methods, the following expressions of the stabilization coefficients are used in the fully discrete schemes
\begin{equation}\label{eq:taut}
\btau_{m}^{n}={\rm{diag}}([\tau_{m}^{n}]^{d}),\text{ with }
\tau_{m}^{n}(K)=\left(\frac{\gamma^{2}}{\Delta t^{2}}+d\,c_{1}^{2}\frac{\nu^{2}}{(h_{K}/k)^{4}}+c_{2}^{2}\frac{U_{K}^{n}}{(h_{K}/k)^{2}}\right)^{-1/2},
\end{equation}
and
\begin{equation}\label{eq:taud}
\tau_{\rm c}^{n}(K)=\frac{(h_{K}/k)^{2}}{d\,c_{1}\tau_{m}^{n}(K)},
\end{equation}
by adapting the form proposed in \cite{Codina02, CodinaBlasco02}, designed by a specific Fourier analysis applied in the framework of stabilized methods. 
In \eqref{eq:taut}-\eqref{eq:taud}, \(\gamma\) denotes the order of accuracy in time, \(d\) is the dimension of the problem, \(c_{1}\) and \(c_{2}\) 
are user-chosen positive constants, \(h_{K}\) is the diameter of element \(K\), \(k\) is the polynomial degree of the velocity FE approximation, and 
\(U_{K}^{n}\) is some local speed on the mesh cell \(K\) at time step \(n\), \(n=0,1,\ldots,N-1\). 
In this work, we have \(\gamma=2\), \(d=2\), and \(k=2\). Also, the values of the constants \(c_{1}\) and \(c_{2}\) are chosen to be \(c_{1}=4\), 
\(c_{2}=\sqrt{c_{1}}=2\) (cf. \cite{Codina01}), and we set \(U_{K}^{n}=||\widehat{\u}_{h}^{n}||^{2}_{{\bf L}^{2}(K)}/|K|\), with \(|K|\) denoting 
the surface (or volume, if \(d=3\)) of element \(K\). Thus, the stabilization coefficients reads
\begin{equation}\label{eq:tauti}
\tau_{m}^{n}(K)=\left(\frac{4}{\Delta t^{2}}+32\frac{\nu^{2}}{(h_{K}/2)^{4}}+4\frac{||\widehat{\u}_{h}^{n}||^{2}_{{\bf L}^{2}(K)}/|K|}{(h_{K}/2)^{2}}\right)^{-1/2},
\end{equation}
and
\begin{equation}\label{eq:taudi}
\tau_{\rm c}^{n}(K)=\frac{(h_{K}/2)^{2}}{8\tau_{m}^{n}(K)}.
\end{equation}
In the following subsections, we specify in detail how it reads the fully discrete scheme for one of each considered method.

\subsection{Semi-implicit BDF2 RB-VMS scheme}
We consider the time discretization of problem \eqref{eq:res_vms} by means of a semi-implicit BDF2 scheme. Similarly to 
\cite{FortiDede15} (section 2), the fully discrete semi-implicit BDF2 RB-VMS scheme consists in solving, for \(n=0,\ldots,N-1\): 

Find \(\u_h^{n+1}\in\vv_h, \ p_h^{n+1} \in Q_h\) satisfying
\begin{align}\label{eq:res_vmsBDF2}
\left(D_{t}^{2}\u_{h}^{n+1},\v_h\right) &+ \nu\left(\nabla{\u_h^{n+1}},\nabla{\v_h}\right) + 
\left((\widehat{\u}_h^{n}\cdot \nabla )\u_h^{n+1},\v_h\right) -	\left(p_h^{n+1},\nabla \v_h\right)+ \left(\nabla \cdot \u_h^{n+1},q_h\right) \nonumber \\
& - \left(\resm(\u_h^{n+1},p_h^{n+1}) ,(\widehat{\u}_h^{n} \cdot \nabla )\v_h + C\nabla q_h\right)
- \left( \resm(\u_h^{n+1},p_h^{n+1}), (\nabla \v_h)^T \widehat{\u}_h^{n}\right) \nonumber \\
& - \left(\resm(\u_h^{n+1},p_h^{n+1}) ,(\nabla \v_h)^T \resm(\widehat{\u}_h^{n},\widehat{p}_h^{n})\right)
+ \left(\tau_{\rm c}^{n}\nabla \cdot \u_h^{n+1} ,\nabla \cdot \v_h\right)	
= (\ff_h^{n+1}, \v_h),
\end{align}
for all \((\v_h,q_h) \in \vv_h\times Q_h\), where
$$
\resm(\u_h^{n+1},p_h^{n+1})=\btau_m^{n}\left(\ff_h^{n+1} - D_t^{2} \u_h^{n+1} + \nu \Delta \u_h^{n+1} -
(\widehat{\u}_h^{n}\cdot \nabla) \u_h^{n+1} - \nabla p_h^{n+1}\right),
$$
and
$$
\resm(\widehat{\u}_h^{n},\widehat{p}_h^{n})=\btau_m^{n}\left(\ff_h^{n+1} - D_t^{2} \widehat{\u}_h^{n} 
+ \nu \Delta \widehat{\u}_h^{n} - (\widehat{\u}_h^{n}\cdot \nabla) \widehat{\u}_h^{n} - \nabla \widehat{p}_h^{n}\right),
$$
with \(\widehat{p}_h^{n}=2 p_{h}^{n}-2 p_{h}^{n-1}\), and \(p_h^{0}=p_h^{-1}\) for \(n=0\), so that 
one has to initialize the pressure (e.g., solve the steady Stokes problem at \(t=0\)).

\subsection{Semi-implicit BDF2 SUPG scheme with grad-div stabilization}
Similarly to \eqref{eq:res_vmsBDF2}, for \(n=0,\ldots,N-1\), the semi-implicit BDF2 SUPG scheme with grad-div stabilization reads:

Find \(\u_h^{n+1}\in\vv_h, \ p_h^{n+1} \in Q_h\) satisfying
\begin{align}\label{eq:SUPGgdBDF2}
\left(D_{t}^{2}\u_{h}^{n+1},\v_h\right) &+ \nu\left(\nabla{\u_h^{n+1}},\nabla{\v_h}\right) + 
\left((\widehat{\u}_h^{n}\cdot \nabla )\u_h^{n+1},\v_h\right) -	\left(p_h^{n+1},\nabla \v_h\right)+ \left(\nabla \cdot \u_h^{n+1},q_h\right) \nonumber \\
& - \left(\resm(\u_h^{n+1},p_h^{n+1}) ,(\widehat{\u}_h^{n} \cdot \nabla )\v_h + C\nabla q_h\right) 
+ \left(\tau_{\rm c}^{n}\nabla \cdot \u_h^{n+1} ,\nabla \cdot \v_h\right)	
= (\ff_h^{n+1}, \v_h),
\end{align}
for all \((\v_h,q_h) \in \vv_h\times Q_h\).

\subsection{Semi-implicit BDF2 LPS schemes}
Apart from the difference in the definition of the projection/interpolation operator \(\boldsymbol{\pi}_h\), the 
semi-implicit BDF2 time discretization of both one-level LPS and LPS by interpolation schemes is given, for \(n=0,\ldots,N-1\), by:

Find \(\u_h \ :\ (0,T)\to\vv_h, \ p_h \ : \ (0,T)\to Q_h\) satisfying
\begin{align}\label{eq:lpsBDF2}
\left(D_{t}^{2}\u_{h}^{n+1},\v_h\right) &+ \nu\left(\nabla{\u_h^{n+1}},\nabla{\v_h}\right) + 
\left((\widehat{\u}_h^{n}\cdot \nabla )\u_h^{n+1},\v_h\right) -	\left(p_h^{n+1},\nabla \v_h\right)+ \left(\nabla \cdot \u_h^{n+1},q_h\right) \nonumber \\
& + \left(\btau_{m}^{n}\boldsymbol{k}_{h}((\widehat{\u}_h^{n} \cdot \nabla) \u_h^{n+1}) ,\boldsymbol{k}_{h}((\widehat{\u}_h^{n} \cdot \nabla) \v_h)\right) 
+ \left(\btau_{m}^{n}\boldsymbol{k}_{h}(\nabla p_h^{n+1}),\boldsymbol{k}_{h}(C\nabla q_h)\right) \nonumber \\
& + \left(\tau_{\rm c}^{n}\nabla \cdot \u_h^{n+1} ,\nabla \cdot \v_h\right)	
= (\ff_h^{n+1}, \v_h),
\end{align}
for all \((\v_h,q_h) \in \vv_h\times Q_h\), where we recall that \(\boldsymbol{k}_{h} = \boldsymbol{I} - \boldsymbol{\pi}_h\) is the fluctuation operator.

\section{Numerical studies: 2D Kelvin--Helmholtz instability}\label{sec:numerics}
In this section, we present the numerical study of a two-dimensional mixing layer problem
evolving in time at Reynolds number \(Re=10^4\). All computations have been performed with
the FE package ParMooN~\cite{ParMooN}, except for the LPS method by interpolation, 
for which we used the FE software \textit{FreeFem++}~\cite{Hecht12}. 

\subsection{Model problem and monitored quantities of interest}
Following a similar setup as described in~\cite{GWR05,Naveed15,SL17}, we briefly summarize the setting 
of the model problem. The problem is defined in \(\Omega = (0,1)^2\). Free slip boundary conditions
are imposed at \(y=0\) and \(y=1\). At \(x=0\) and \(x=1\), periodic boundary conditions are 
prescribed. There is no external forcing, that is \(\ff=0\). The initial velocity field is given by
\[ \u_0 = 
\left( \begin{array}{cc}
U_\infty \tanh ((2y-1)/\delta_0) \\ 
\\
0
\end{array}
\right) + c_{n} U_\infty \left(
\begin{array}{cc}
\partial_{y} \psi\\
\\
-\partial_{x}\psi
\end{array}\right),
\]
where \(U_\infty\) is a reference velocity, \(\delta_0\) is the initial vorticity thickness that will be defined later, \(c_{n}\) is 
a parameter giving the strength of perturbation, and the stream function is given by
\[
\psi = \exp \left(- ((y-0.5)/\delta_0)^2 \right) \left(\cos (8\pi x) + cos(20\pi y)\right).
\]
Let the initial vorticity thickness \(\delta_0 = 1/28\), \(U_\infty = 1\), the scaling/noise factor 
\(c_n=10^{-3}\), and the inverse of viscosity \(\nu^{-1} = 28 \times 10^4\). Thus, the Reynolds number 
associated with the flow is \(Re = U_\infty \delta_0/\nu = 10^4 \). The mixing layer problem is known to be inviscid 
unstable, thus the chosen small viscosity makes the solution very sensitive. Slight perturbations of the initial condition 
are amplified by the so-called Kelvin--Helmholtz instabilities. Because of the unstable nature of the problem, this is a 
challenging test case for the study of 2D turbulence and vortex dynamics in free shear layers of incompressible flows 
(cf.~\cite{Lesieur88}).

Several attempts have been made in the literature to numerically investigate the Kelvin--Helmholtz instabilities caused 
by slight perturbations in the initial condition of the described model problem (both in 2D and 3D).
In particular, it has been deeply discussed in~\cite{Lesieur88}, where a direct numerical simulation of a two-dimensional
temporal mixing layer problem was performed, applying a second-order finite difference method at the
high resolution of \(256^2\) grid points with a uniform spacing in each direction. Further numerical studies for this 
problem, including LES, VMS and stabilized models, may be found, e.g., in~\cite{Boersma97,Burman07,GWR05,John05,Naveed15,SL17}. 
The corresponding three-dimensional case has been numerically analyzed, e.g., in~\cite{Balaras01,John05}.

For the evaluation of computational results, we consider the vorticity of the flow 
\[\omega = \nabla \times \u = \partial_x u_2 - \partial_y u_1. \] 
The vorticity thickness is defined by 
\[
\delta(t_n) = \frac{2U_\infty}{\sup_{y\in[0,1]} | \langle \omega \rangle(y,t_n) |},
\]
where \( \langle \omega\rangle (y,t_n) \) is the integral mean in the periodic direction and is defined as 
\[
\displaystyle	\langle \omega\rangle (y,t_n)  = \frac{\displaystyle\int_0^1 \omega ({\bf x}, t_n)dx}{\displaystyle\int_0^1 dx }
 = \int_0^1 \omega ({\bf x}, t_n) dx.
\]
In the computations, this integral was evaluated discretely for all grid lines parallel to the \(x\)-axis (cf.~\cite{John05}), 
and the maximum of the computed values was taken to obtain \(\delta(t_n)\). In the evaluation of computations, we considered 
the vorticity thickness relative to \(\delta_0:\delta(t_n)/\delta_0\). 

The understanding of the physical evolution of the flow is either done qualitatively, by visualizing the evolution of the vorticity 
field through meaningful instants, or quantitatively, by determining the evolution of the relative vorticity thickness. 
Some conclusions can be drawn depending on the pairing, position of the eddies, time at which the pairing happens, and values of 
the peaks of the relative vorticity thickness, corresponding to the pairing of eddies. The general behavior of the vorticity field 
is as follows. Starting from the noisy initial condition \(\u_0\), four primary eddies are developed, which then pair to two larger
secondary eddies that are standing for a long time. Finally, the pairing of secondary eddies leads to one larger eddy, rotating at 
a fixed position. It can be found in the literature that, depending on the numerical method used for the simulations, the position 
of the final eddy is located either at the center of the domain~\cite{Burman07,SL17} or at the periodic 
boundaries~\cite{GWR05,John05,Naveed15}. 
In Section \ref{subsec:VortField}, prior to a quantitative analysis, plots of the vorticity are shown, obtained by a reference 
simulation, which consists of a PSPG method with an additional grad-div stabilization term, computed on a high resolution level. 
A comparison with results from the literature is performed. Complementing the visualization of the vorticity field, the temporal 
evolution of the relative vorticity thickness obtained with the studied methods on different refinement levels and on different 
time step lengths is discussed in Section \ref{subsec:VortThick}.

In addition to the relative vorticity thickness, we are interested in studying also the temporal evolution of the following quantities of interest. 
The kinetic energy of the flow is the most frequently monitored quantity, given by
\[
\text{kinetic energy} \qquad E_{\rm Kin} = \frac12\|\u(t)\|_{L^2(\Omega)}^2=\frac12 \int_{\Omega} | \u(t,\x)|^2 d\x.
\]
For the studied problem, the physically correct behavior of \(E_{\rm Kin}\) is that it strongly monotonically decreases. 
In Section \ref{subsec:EKin}, we will illustrate the temporal evolution of \(E_{\rm Kin}\) in our conducted numerical simulations 
for all the studied methods on different refinement levels and on different time step lengths. 

The next studied quantity of interest is the enstrophy, defined as
\[
 \text{enstrophy} \qquad \mathcal{E} = \frac12 \|\nabla \times \u(t)\|_{L^2(\Omega)}^2 = \frac12\|\omega(t)\|_{L^2(\Omega)}^2
 =\frac12 \int_\Omega |\nabla \times \omega(t,\x)|^2 d\x .
\]
Similar to the kinetic energy, the enstrophy cannot increase. Numerical studies presented 
in~\cite{SL17} shows that the physically correct behavior is a monotone decline from its initial value. 
Furthermore, a more accurate method with a higher resolution leads to a later decrease in 
enstrophy~\cite{SL17}. 
This quantity of interest has been investigated also by other several authors, for details see~\cite{LM96,SF00,EAJ07}.
In Section \ref{subsec:Enst}, we will illustrate the temporal evolution of \(\mathcal{E}\) in our conducted numerical 
simulations for all the studied methods on different refinement levels and on different time step lengths. 

Finally, we will investigate another important and challenging quantity of interest to be monitored, known as palinstrophy, 
which in the context of 2D turbulence drives the dissipation process. Palinstrophy is defined by 
\[
 \text{palinstrophy} \qquad \mathcal{P} = \frac12 \|\nabla \omega(t) \|_{L^2(\Omega)}^2 
 =\frac12 \int_\Omega |\nabla \omega(t,\x)|^2 d\x.
\]
Note that, in contrast to \(E_{\rm Kin}\) and \(\mathcal{E}\), \(\mathcal{P}\) can increase in time (cf.~\cite{DoeringGibbon95}, Section 3.3).
In Section \ref{subsec:Palinst}, we will illustrate the temporal evolution of \(\mathcal{P}\) in our conducted numerical simulations 
for all the studied methods on different refinement levels and on different time step lengths. 

Note that all quantities of interest will be compared with the reference solution. 

\subsection{Preliminaries to numerical simulations}
Our calculations were carried on uniform triangular grids where the coarsest grid (Level 0) is obtained by 
dividing the unit square into two triangles. This grid is refined uniformly and the number of degrees of 
freedom on finer grids is given in Table~\ref{dofs} for different FE spaces used in the simulation. 
We show how sensitive the solution is towards mesh refinement, by comparing three different refined levels 
of resolution (Level 5, 6 and 7) that represent under-resolved to well-resolved
situations (see Table~\ref{dofs}).

The time discretization is performed for all methods with the semi-implicit BDF2 schemes 
described in the previous section, using equidistant time steps of 
length \(\Delta t=1.25\times 10^{-2}\) and \(\Delta t = 3.125\times 10^{-3}\).
For the simplicity of presentation, we will use \(\Delta t_1\) and \(\Delta t_2\) as abbreviation for 
the large and small time step lengths. The final time is set to be \(T=7.15\).

For the one-level variant of LPS method that needs enriched FE spaces for velocities, 
we used mapped FE spaces \cite{ParMooN}, where the enriched space on the reference cell 
\(\widehat{K} = (-1,1)^2\) is defined by 
\[
\mathbb{P}_2^{\rm bubble} (\widehat{K}) =\mathbb{P}_2(\widehat{K}) + \widehat{b}_\triangle \mathbb{P}_1(\widehat{K}),
\]
with \(\widehat{b}_\triangle\) the cubic bubble on the reference triangle. 
Together with the choice \(D_h(M) = \mathbb{P}_1^{dc}(M)\) for the projection space, this space is suited for classical one-level LPS methods. 
Also, for the one-level variant of LPS method, numerical studies concerning the choice of stabilization parameters 
suggests that a good choice is \(\btau_{m}=C_0 h_K\) and \(\tau_c=C_0 h_K\), where \(\delta_0\in 
(0,1)\), see~\cite{AM16}. Based on these studies and on our own experience, the parameter \(C_0\) is 
set to be \(0.1\) in all simulations for the one-level variant of LPS method. 
For all other methods, standard \(\mathbb{P}_2\) FE spaces were used for velocities. 

All monitored quantities of interests for the Kelvin--Helmholtz instability problem in our computational results are 
compared with the reference solution obtained by a PSPG method together with grad-div stabilization using $\mathbf{P}_2/\mathbb{P}_2$ 
FE on a very fine mesh (Level 8), and the small time step length \(\Delta t= 3.125\times 10^{-3}\). 
In addition to that, we will also compare our results with those ones presented in~\cite{SL17}. In~\cite{SL17},
for the same setup of the problem, numerical studies were performed with higher-order divergence-free 
FE on finer meshes. More precisely, exactly divergence-free \(\mathbf{H}\)(div) based 
on Raviart--Thomas FE of order 3 (RT3) were used on four different refinement levels in space for velocities. 
For the time discretization, a multi-step IMEX time stepping scheme based on BDF2 that combines 
BDF2 with a second-order Adams-Bashforth scheme were applied. In comparison to that, 
our computational results are obtained using almost two-times coarser meshes (compared to Level a, b, c in~\cite{SL17}, Table 2), and much cheaper FE and time discretizations.

In the following, each monitored quantity of interest will be discussed and compared separately for all the methods 
presented in the previous sections. Numerical simulations were done both with EO $\mathbf{P}_2/\mathbb{P}_2$ 
and ISS $\mathbf{P}_2/\mathbb{P}_1$ FE for the pair velocity/pressure on different refinement levels. In the case of the one-level LPS method, 
EO $\mathbf{P}_2^{\rm bubble}/\mathbb{P}_2^{\rm bubble}$ and ISS $\mathbf{P}_2^{\rm bubble}/\mathbb{P}_1^{\rm 
dc}$ FE pair are used. We will also analyze
in detail the effect of time step lengths on the computational results.

\begin{table}[t]
\begin{center}
\caption{Overview of meshes and degrees of freedom (d.o.f.).}\label{dofs}
\begin{tabular}{c|c|l|l|l|ll}
Level & \(h\) & \(\mathbf{P}_2\) d.o.f. & \(\mathbb{P}_1\) d.o.f. & \(\mathbf{P}_2^{\rm bubble}\)  d.o.f.
& \(\mathbb{P}_1^{\rm dc}\) d.o.f. &
\\
\hline
5& \(4.419\times 10^{-2}\) & 8 320   & 1 056 & 12 416 & 6 144 &  \\
6&  \(2.210\times 10^{-2}\) & 33 024 & 4 160 & 49 408 & 24 576 & \\
7&  \(1.105\times 10^{-2}\)  & 131 584 & 16 512 & 197 120 & 98 304 & \\
\hline
\end{tabular}
\end{center}
\end{table}

\subsection{Evolution of the flow}\label{subsec:VortField}
\begin{figure}[t!]
	\centering
	\includegraphics[scale=0.45]{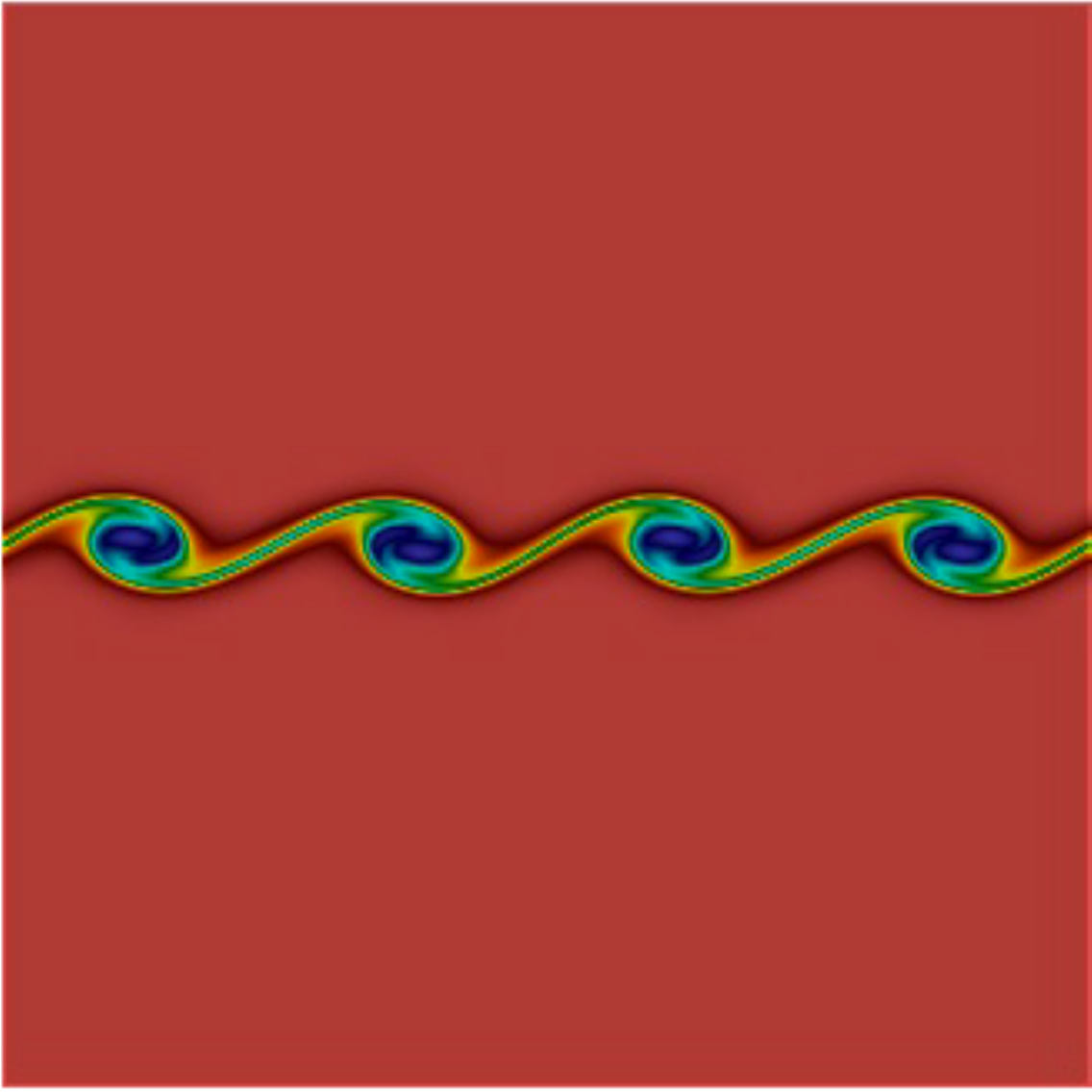}
	\includegraphics[scale=0.45]{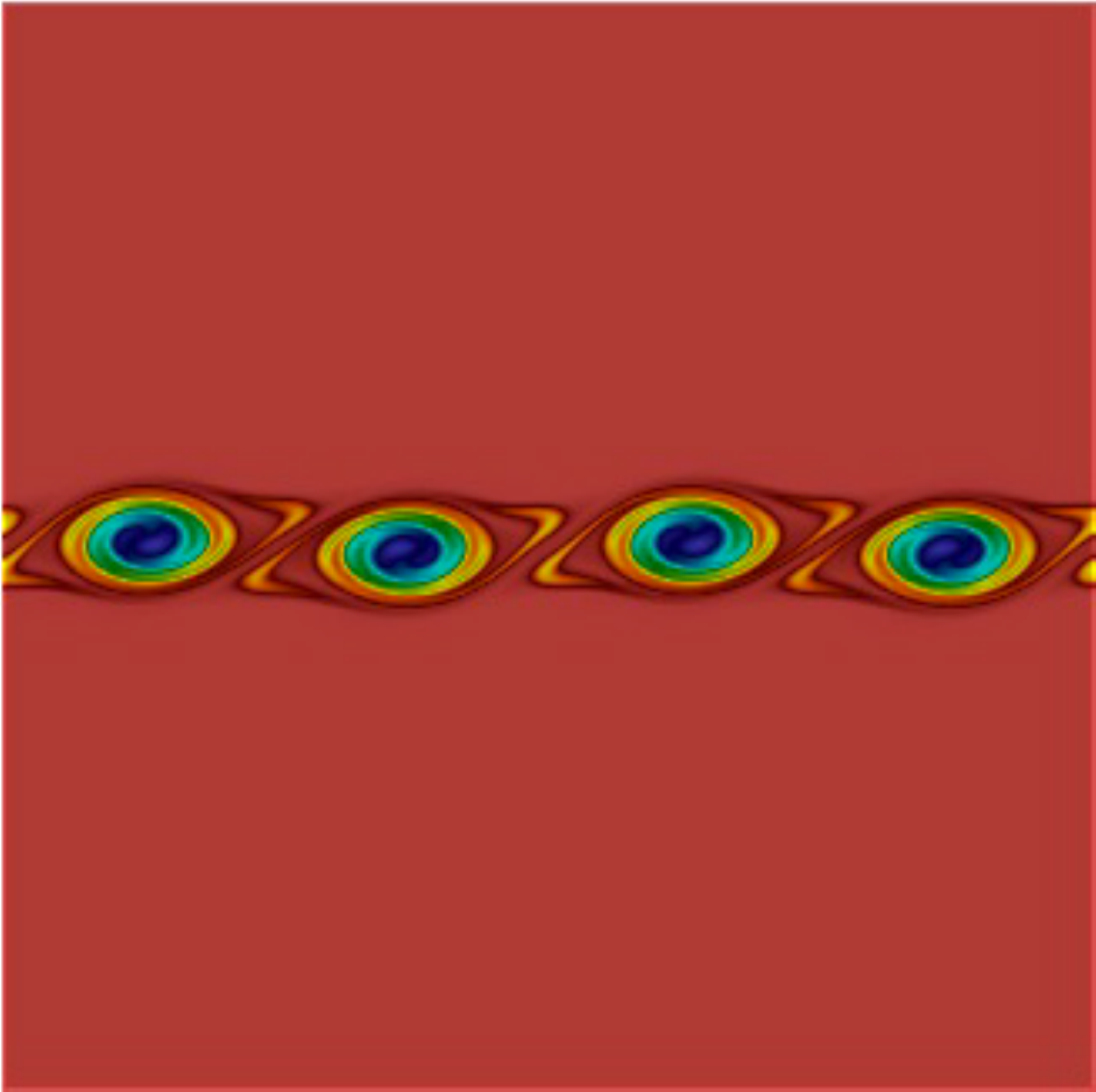}
	\includegraphics[scale=0.45]{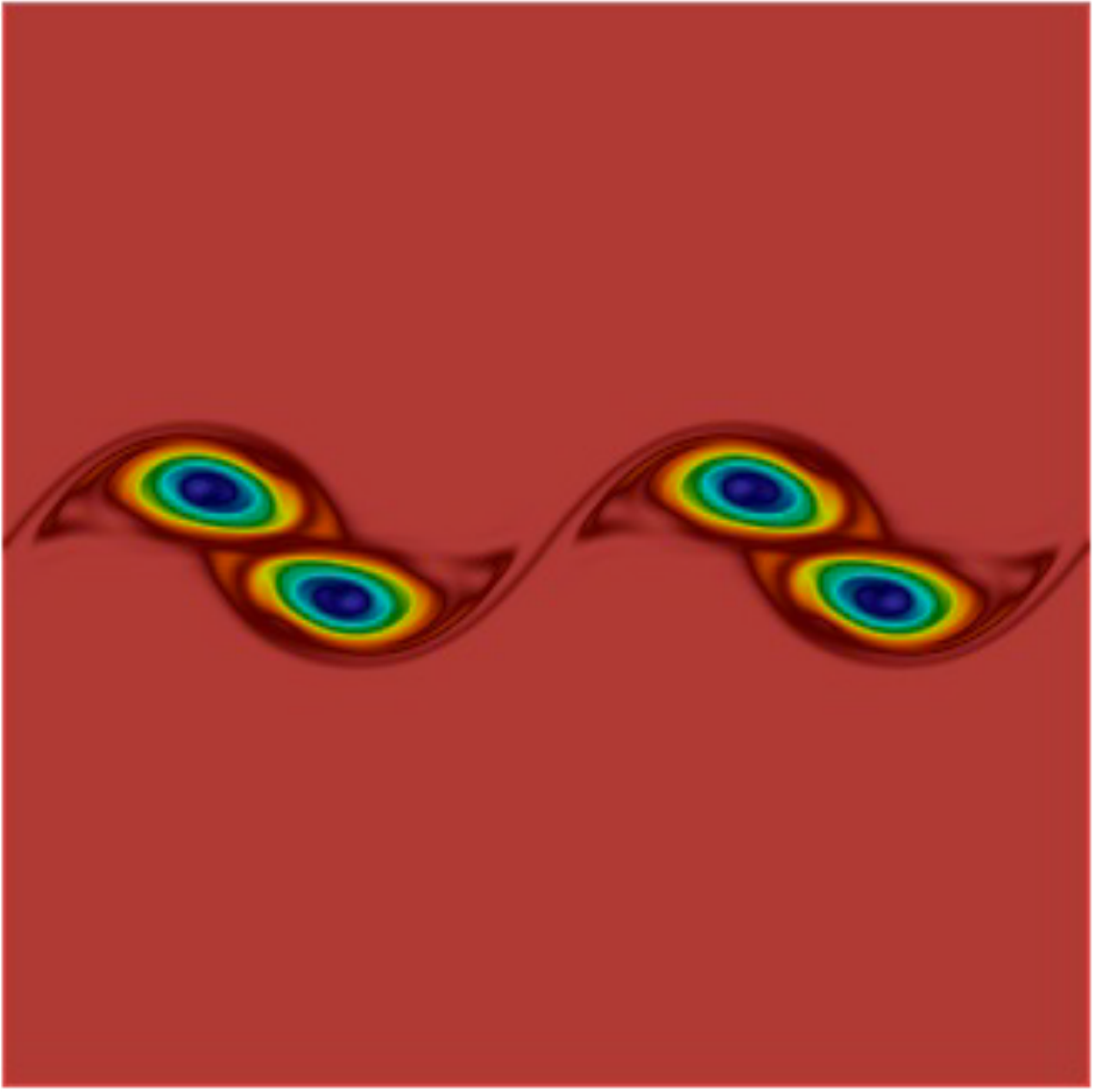}\\
	\includegraphics[scale=0.45]{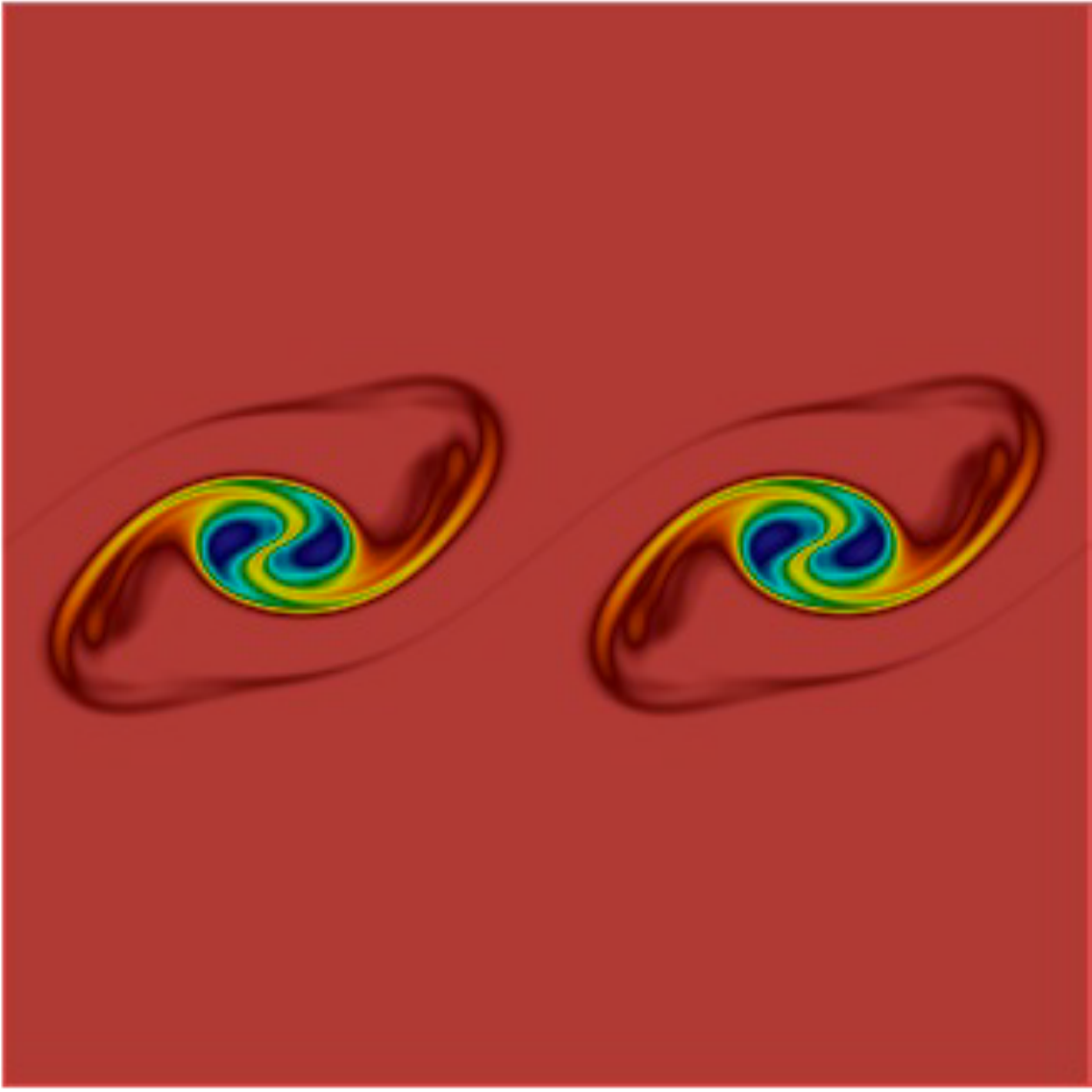}
	\includegraphics[scale=0.45]{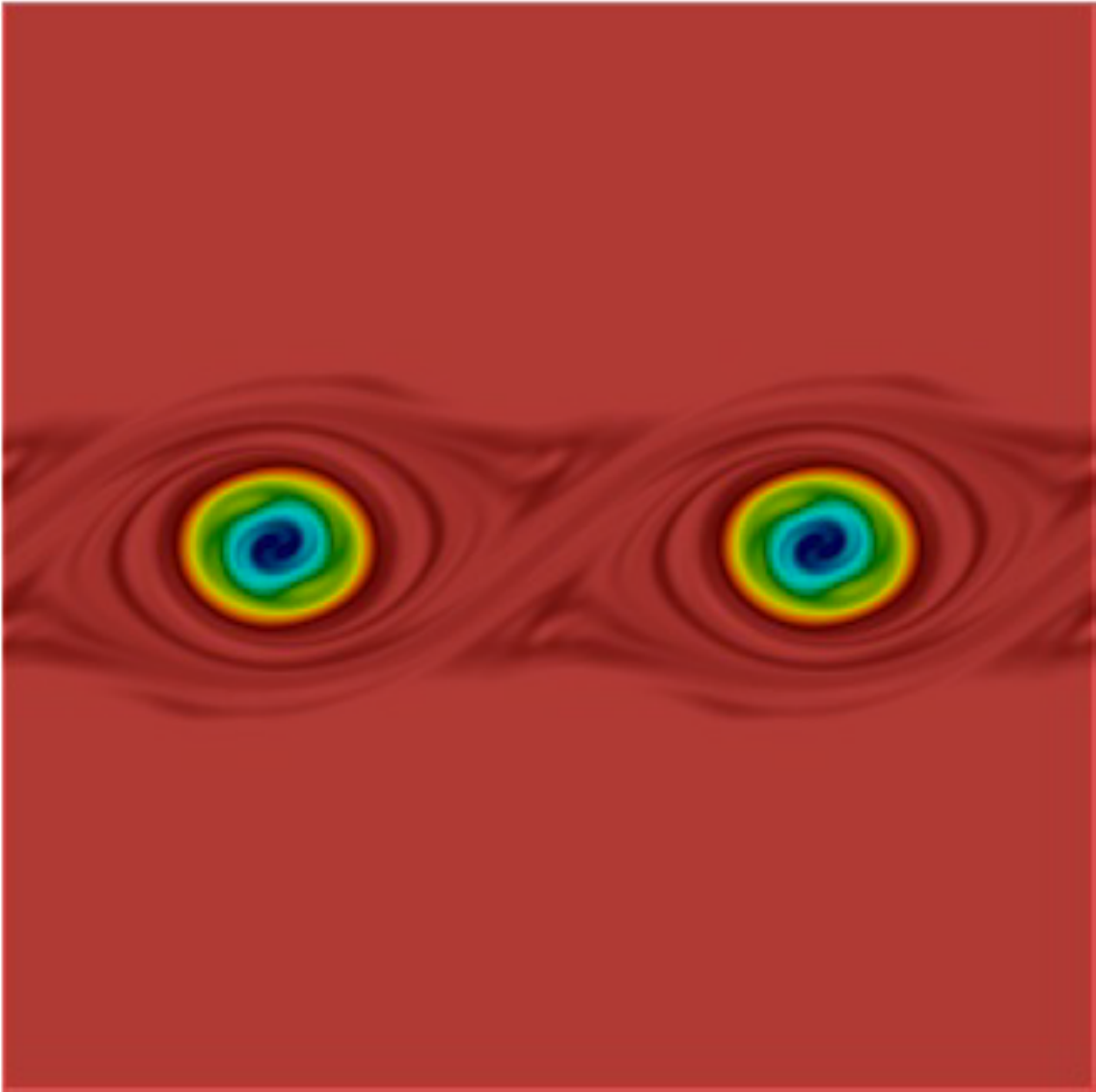}
	\includegraphics[scale=0.45]{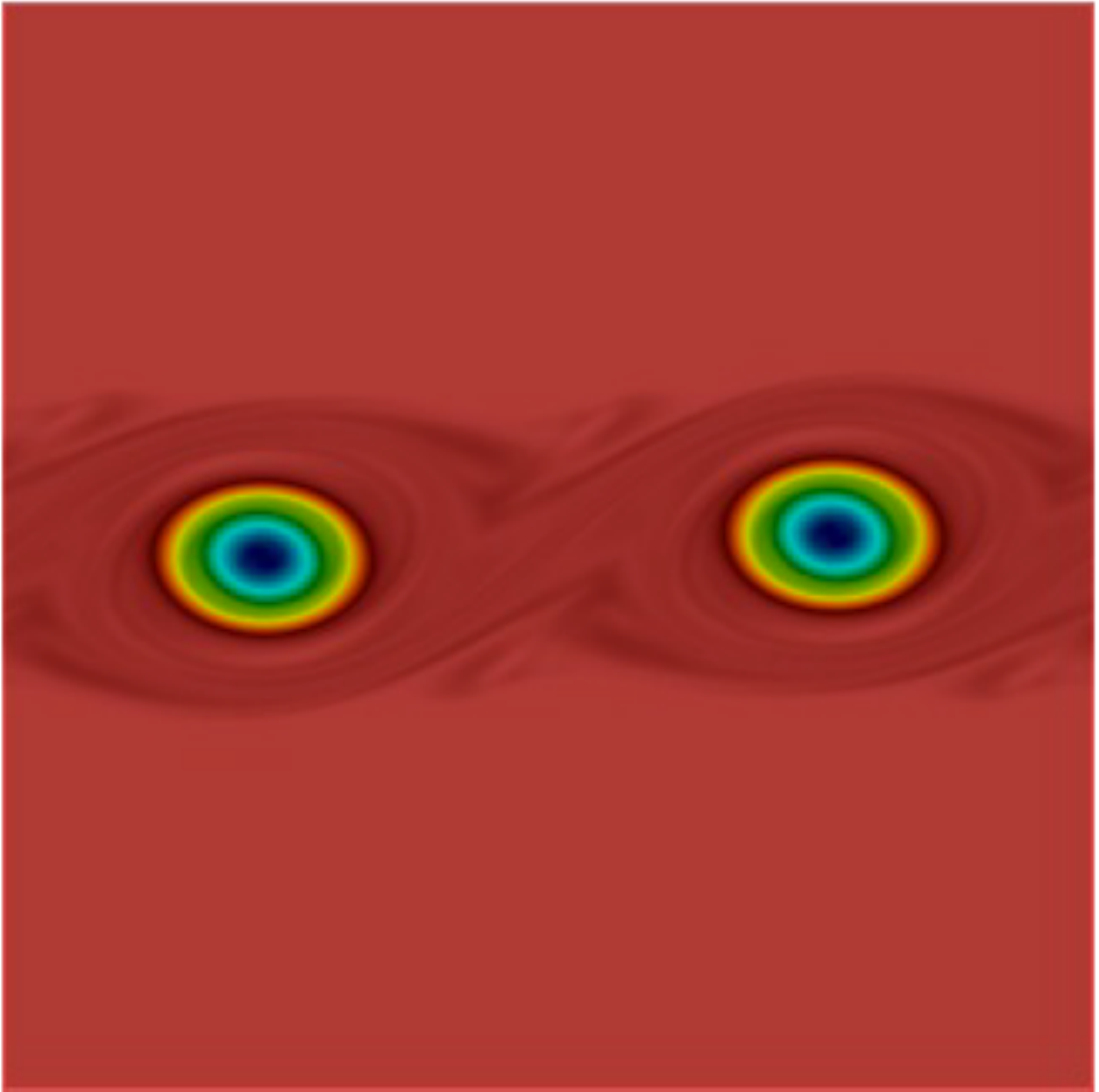}\\
	\includegraphics[scale=0.45]{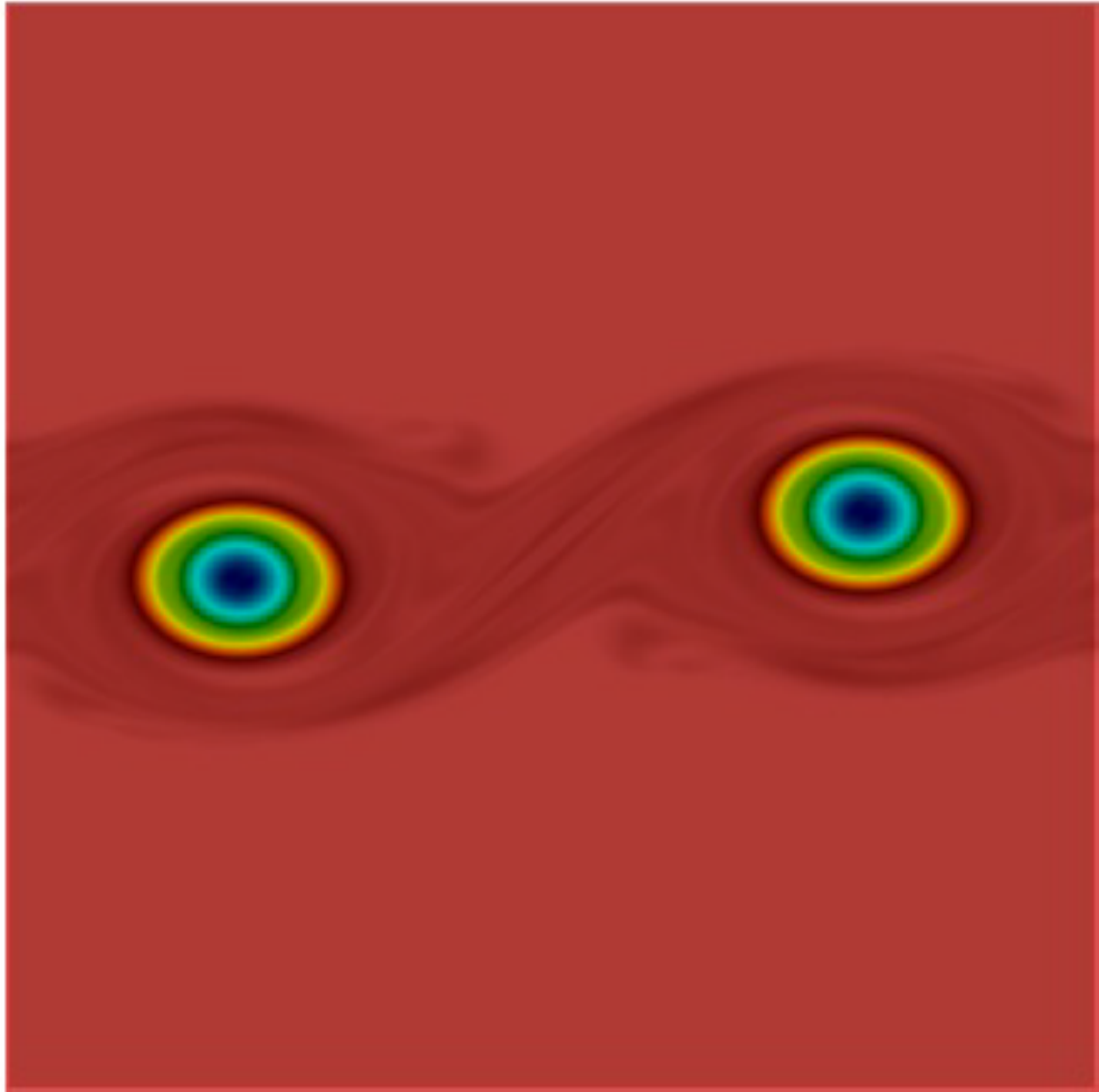}
	\includegraphics[scale=0.45]{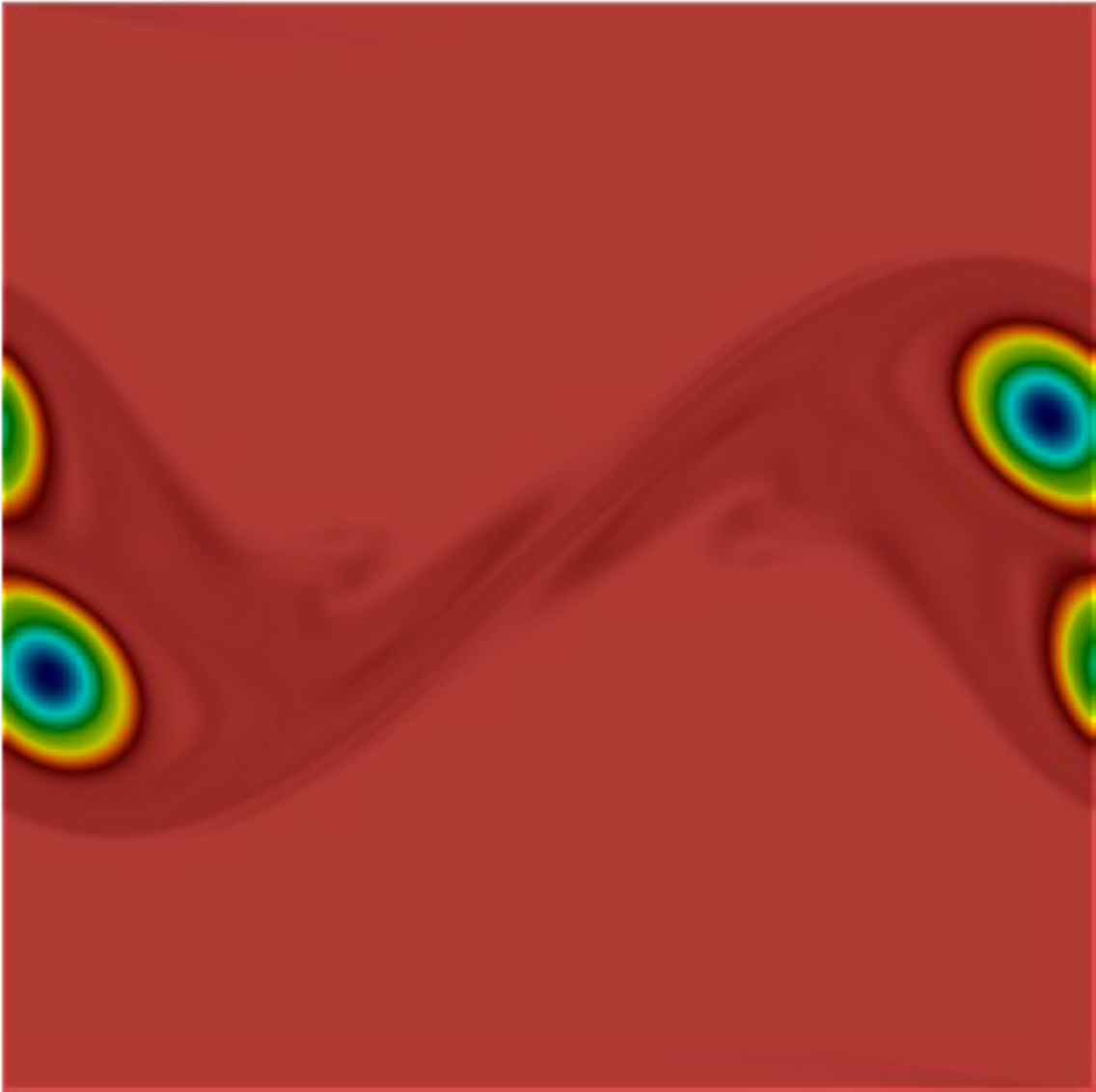}
	\includegraphics[scale=0.45]{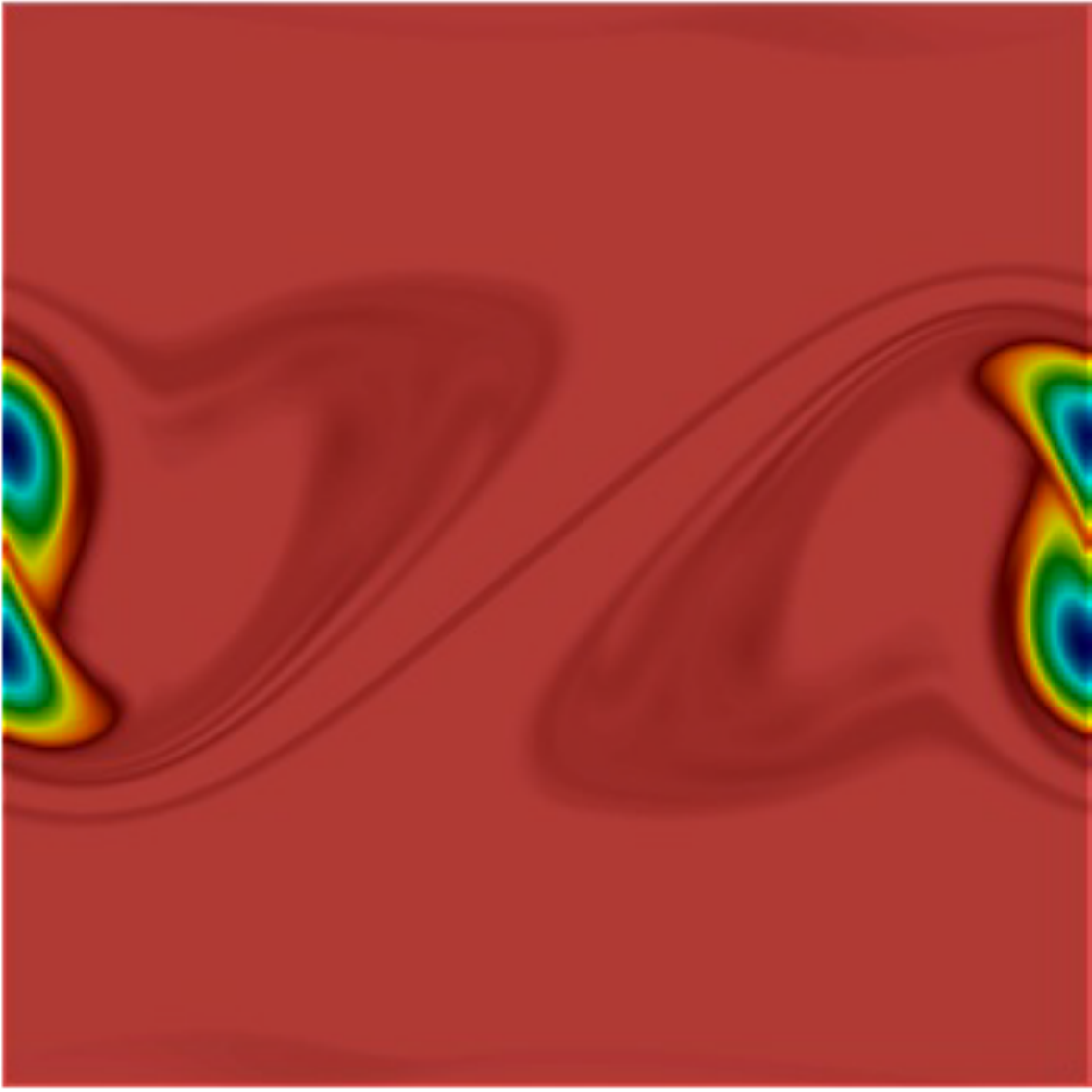}
	\caption{Colored vorticity field (blue: intense vorticity, red: irrotational outer flow) for reference solution (PSPG method with 
	grad-div stabilization on Level $8$ using \(\mathbf{P}_2/\mathbb{P}_2\) FE and semi-implicit BDF2 scheme with 
	\(\Delta t = 3.125\times 10^{-3}\)) at time units \(10, \; 20, \; 30, \; 40, \; 100, \; 155, \;165, \;180, \;200\) (left to right, top to bottom). 
\label{fig:vort_field}}
\end{figure}

The physical evolution of the flow can be described with the help of the 
involved vortices presented in Figure~\ref{fig:vort_field}. These results correspond to the reference solution, obtained using 
the PSPG method with grad-div stabilization on Level 8, with \((\mathbf{P}_2/\mathbb{P}_2)\) FE and 
\(\Delta t = 3.125 \times 10^{-3}\) for the semi-implicit BDF2 scheme. 
In particular, Figure~\ref{fig:vort_field} present the evolution of the vorticity through meaningful time instants. To compare 
the results, the vorticity pictures for our reference solution are shown at the same time instants as in~\cite{SL17}. 
It can be seen that, starting with the initial noise, four primary vortices develop between $10$ and
$20$ time units \(\bar{t}=\delta_{0}/U_{\infty}\), in agreement with~\cite{SL17}. 
Also, the four primary vortices merge at about 35 time units, as observed in~\cite{SL17}.
Always in agreement with~\cite{SL17}, the two secondary vortices are standing for a certain
period of time. However, the instant in time where the second pairing occurs is strongly dependent on which method, solver and resolution are used. 
For instance, considering \(t=155\bar{t}\), the two primary 
vortices in our simulations are still clearly separated and almost aligned parallel to the \(x\)-axis, and they start to approach each other towards
the periodic boundary at \(t= 165 \bar{t}\). On the other side, in~\cite{SL17}, the two secondary vortices are already moving towards each other 
at  \(t=155\bar{t}\), and the last vortex is located in the center of the domain.
Thus, independent of the time instant in which the last pairing occurs,
our results are in agreement with the ones presented in~\cite{John05,GWR05,Naveed15}, where the final vortex rotates 
near the periodic boundary, while in~\cite{Burman07,SL17} the last vortex rotates in the center of the domain, so there is no consensus in the literature 
concerning the location of the last vortex, and one can conclude that different discrete settings generally lead to different final states. Note that 
besides the main vortices, fine-scale flow structures too are captured very well.
Such structures are not so numerically dissipated by the proposed method, which thus
gives a better resolved evolution with respect to~\cite{GWR05,Naveed15}.

\subsection{Vorticity Thickness}\label{subsec:VortThick}
\begin{figure}[t!]
	\centering
	\includegraphics[scale=0.5]{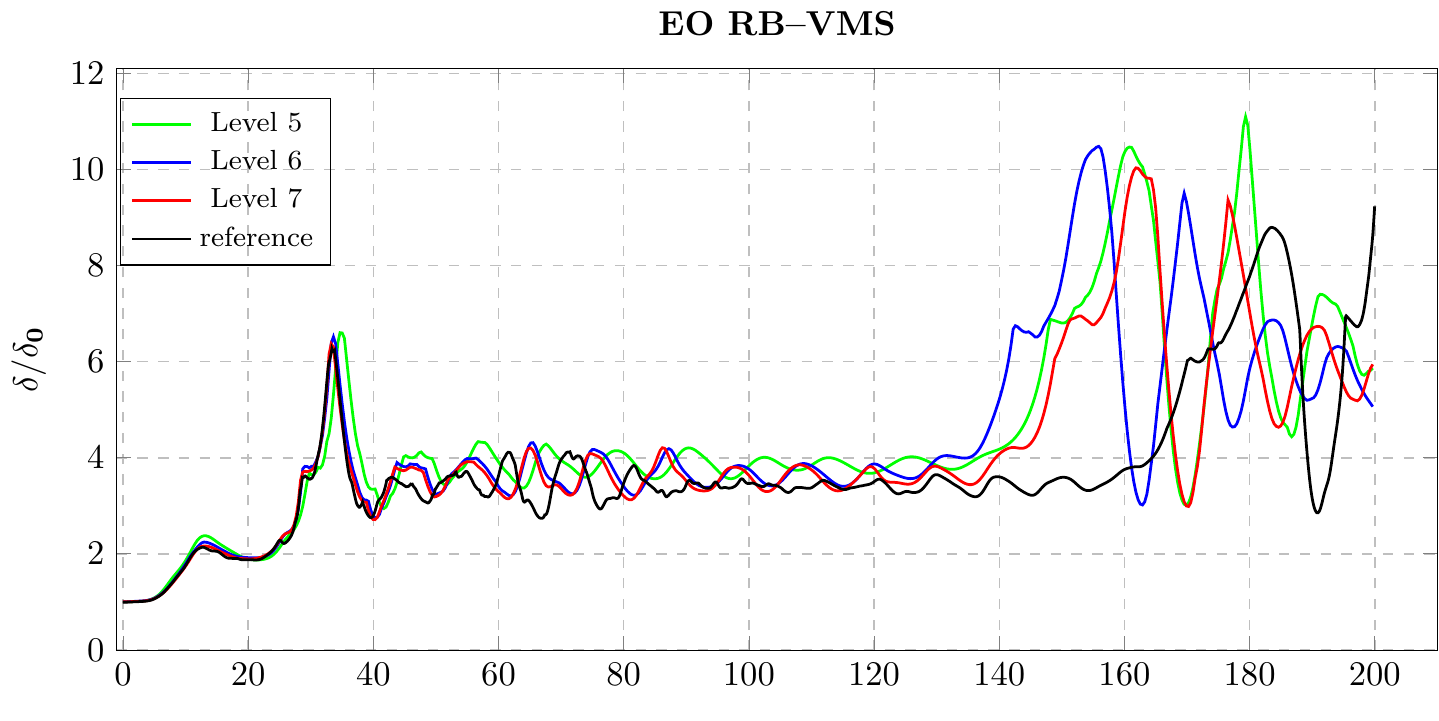}
	\includegraphics[scale=0.5]{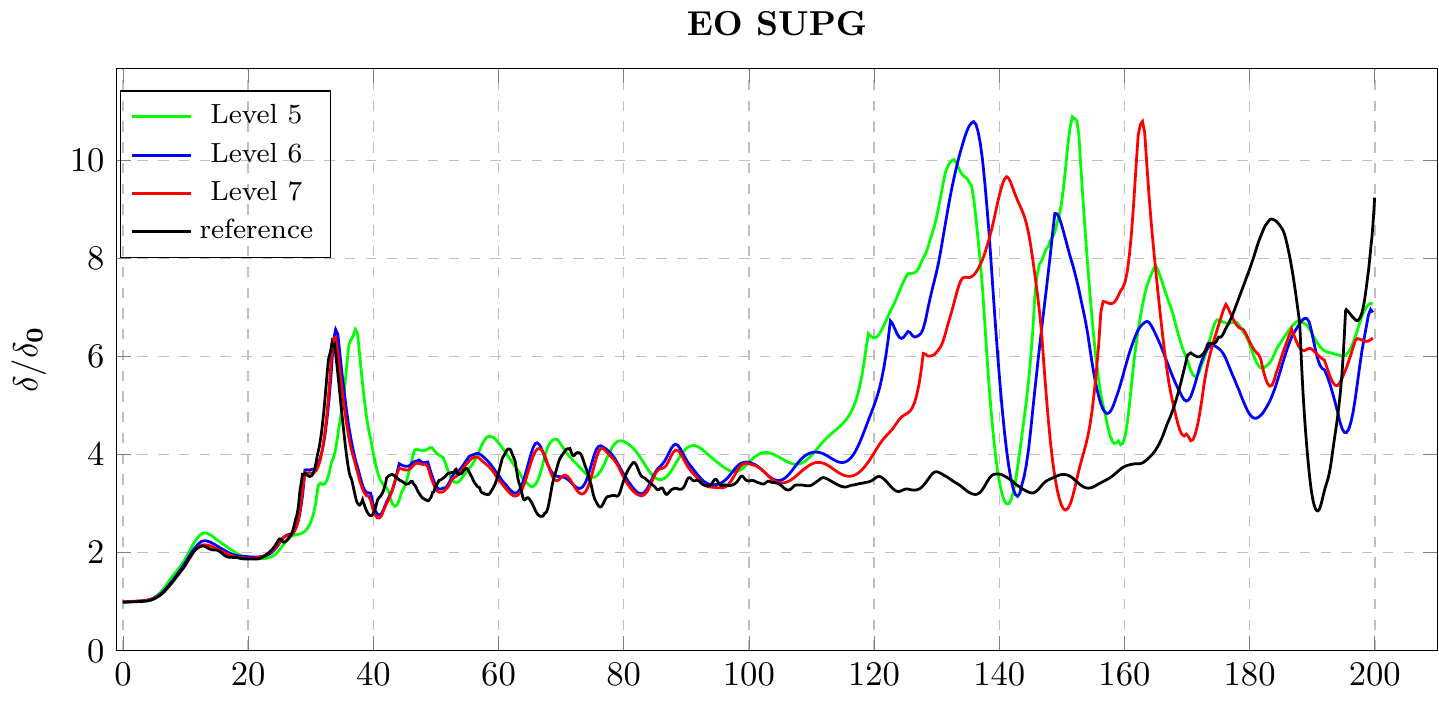}\\
	\includegraphics[scale=0.5]{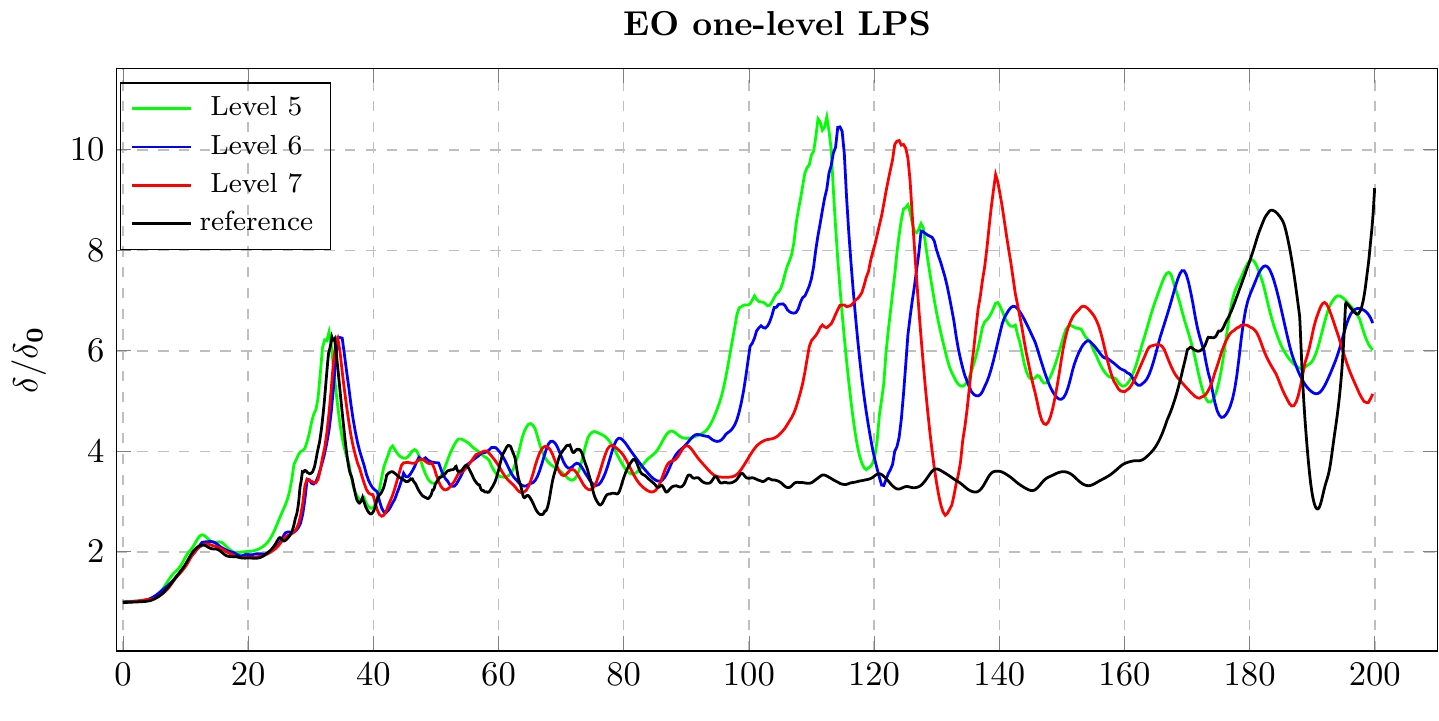}
	\includegraphics[scale=0.5]{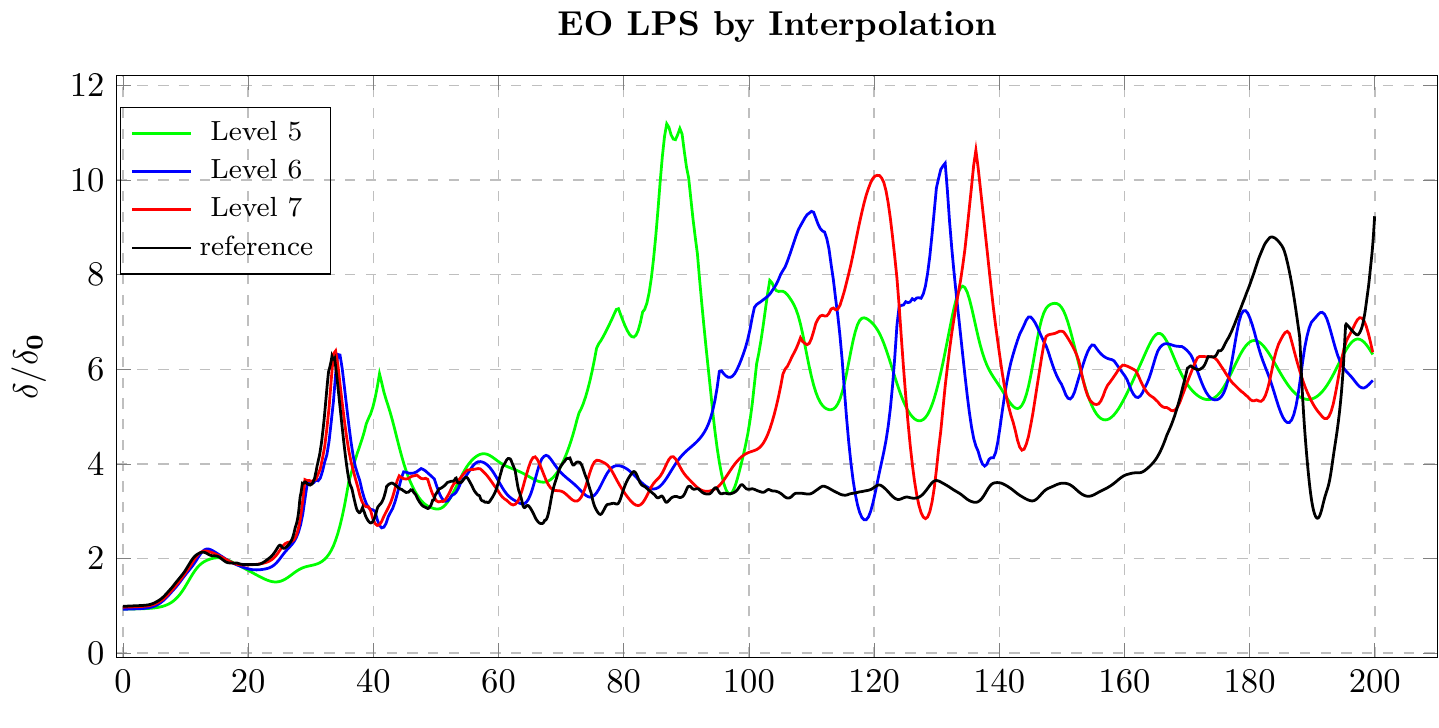}
	\caption{Temporal evolution of vorticity thickness with EO--FE: RB-VMS (top left), SUPG (top right), one-level LPS (bottom left), 
		and LPS by interpolation (bottom right), on different mesh refinement levels, $\Delta t=0.0125$. \label{fig:vort_eo_0.0125}}
\end{figure}

\begin{figure}[t!]
	\centering
	\includegraphics[scale=0.5]{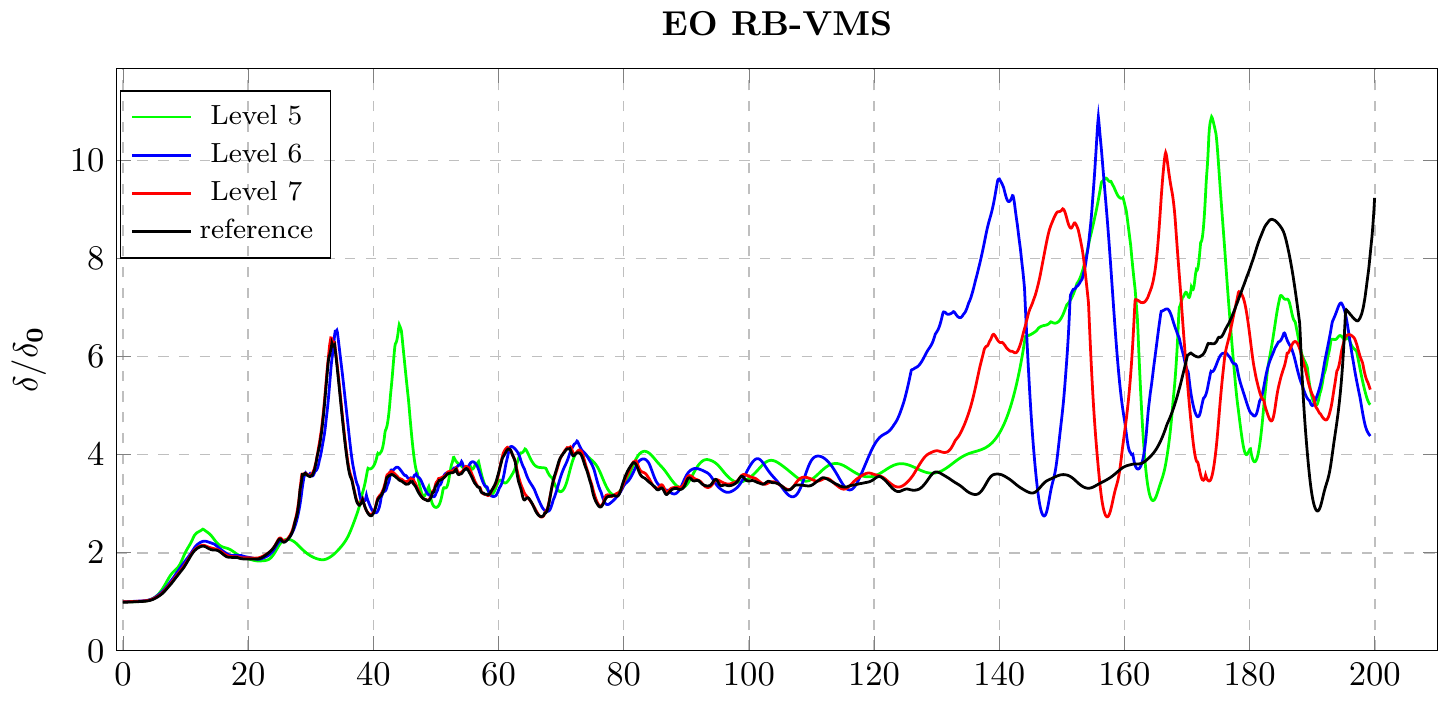}
	\includegraphics[scale=0.5]{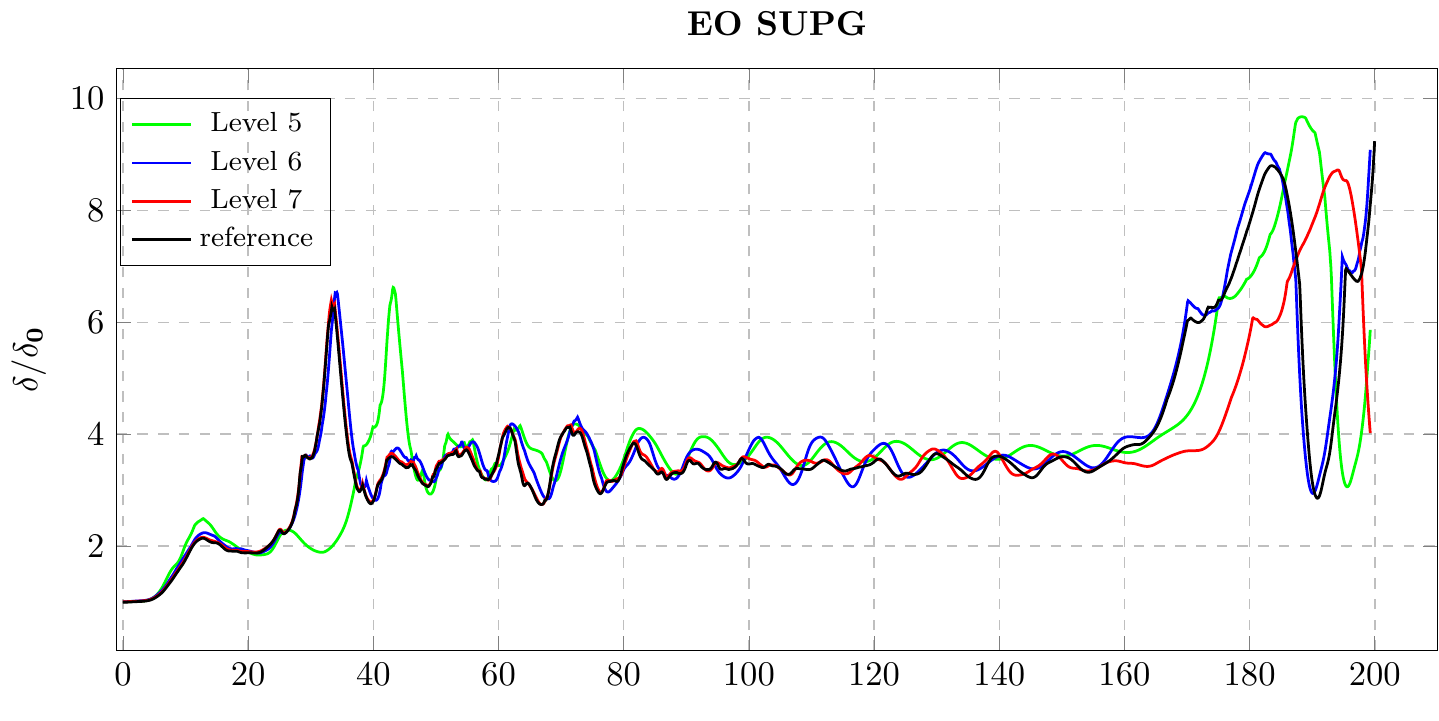}\\
	\includegraphics[scale=0.5]{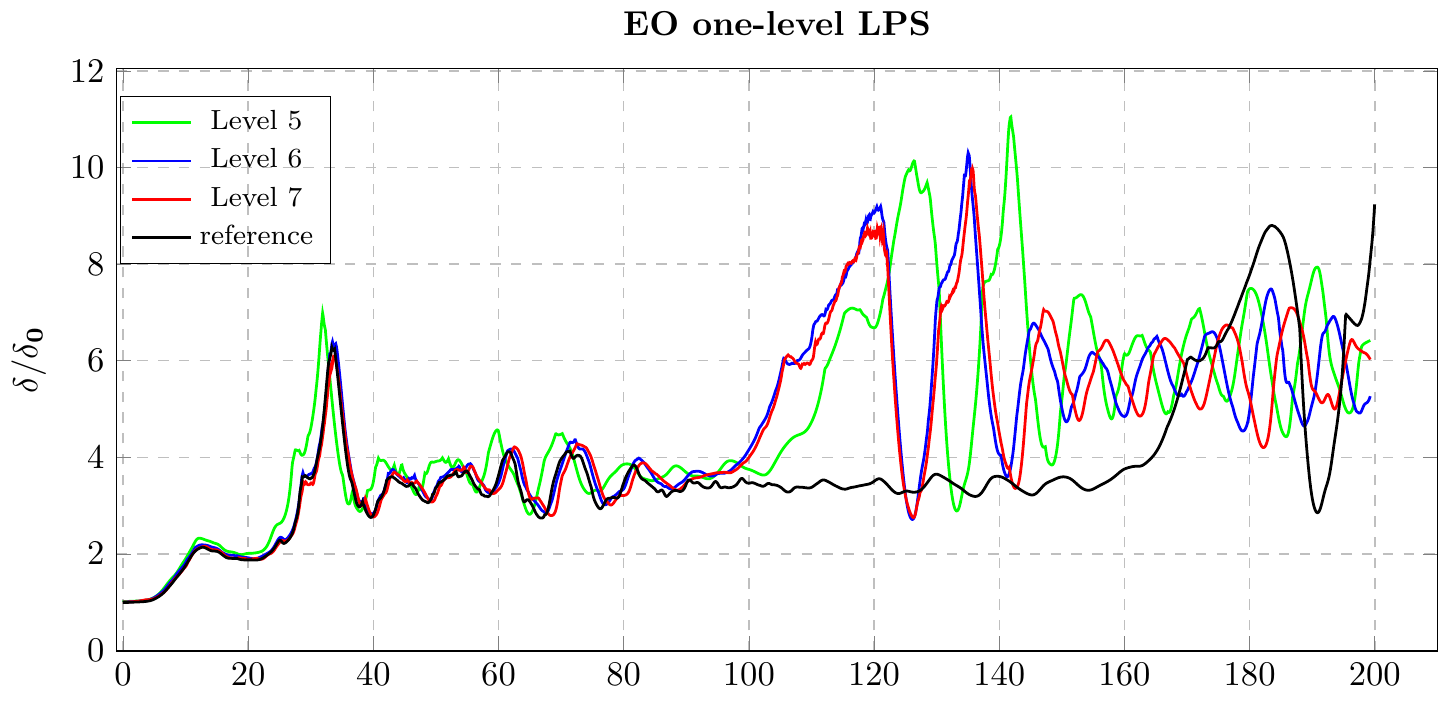}
	\includegraphics[scale=0.5]{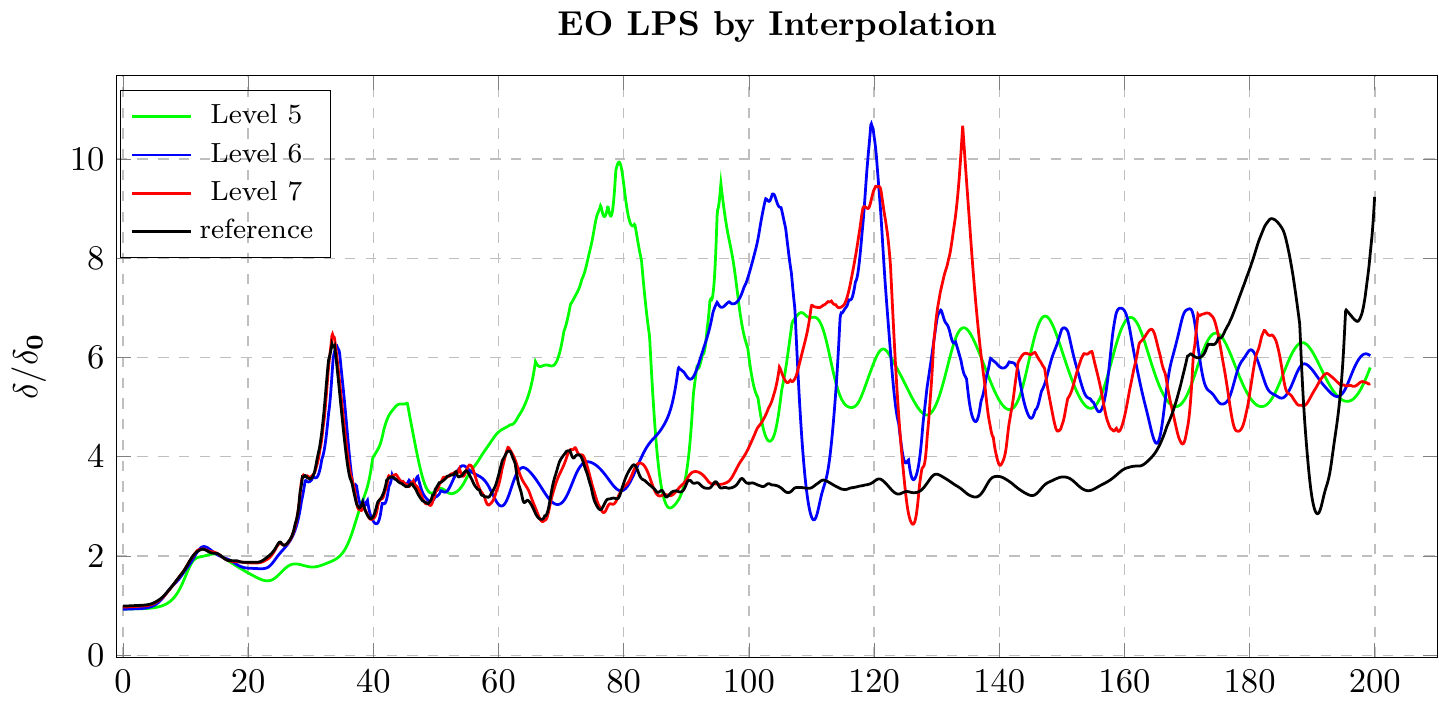}
	\caption{Temporal evolution of vorticity thickness with EO--FE: RB-VMS (top left), SUPG (top right), one-level LPS (bottom left), and LPS by 
		interpolation (bottom right), on different mesh refinement levels, $\Delta t=0.003125$. \label{fig:vort_eo_0.003125}}
\end{figure}

\begin{figure}[t!]
	\centering
	\includegraphics[scale=0.5]{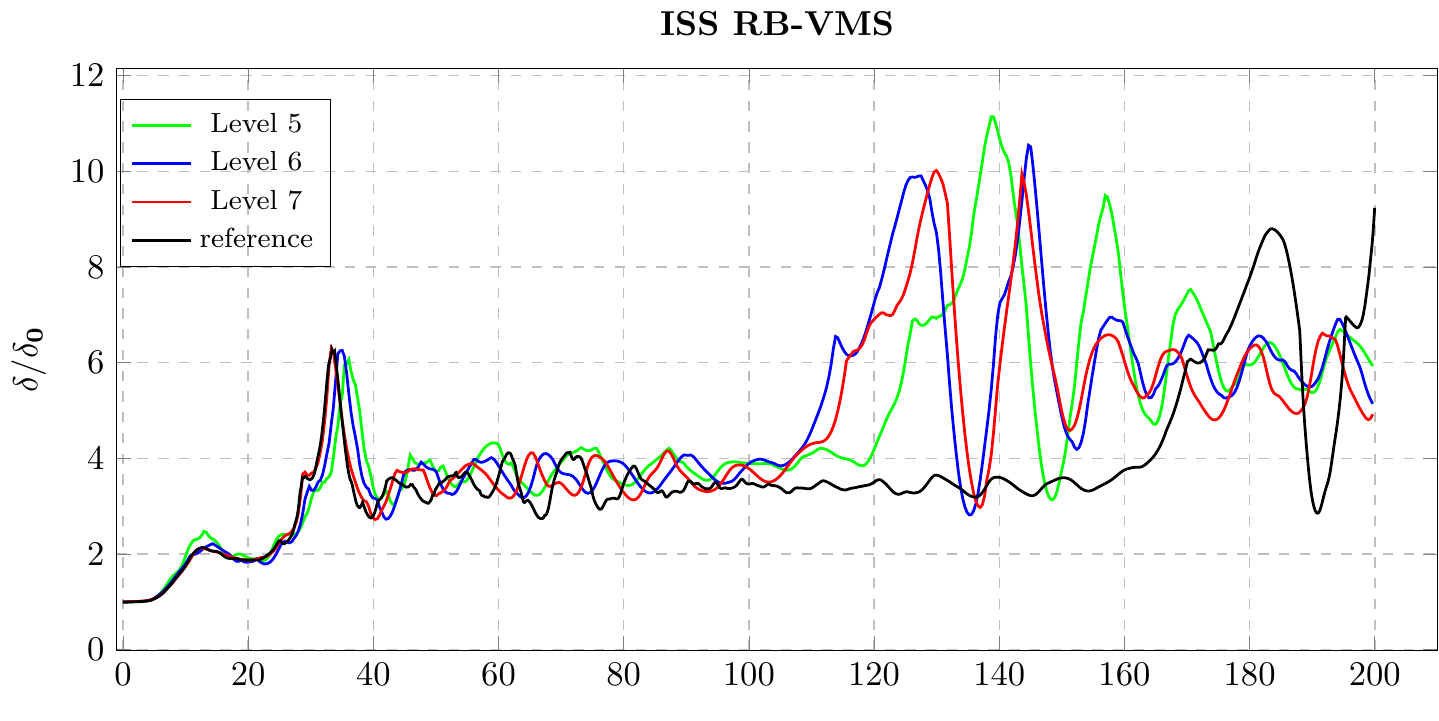}
	\includegraphics[scale=0.5]{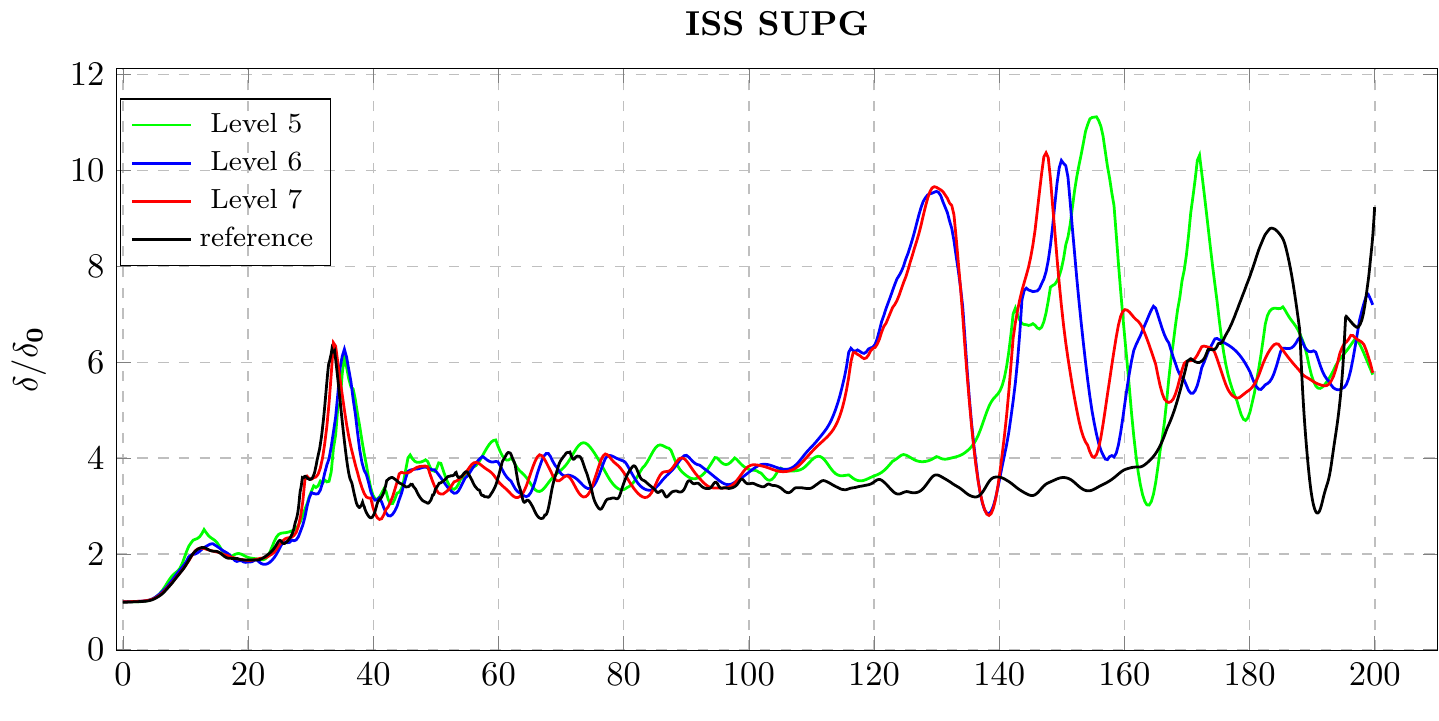}\\
	\includegraphics[scale=0.5]{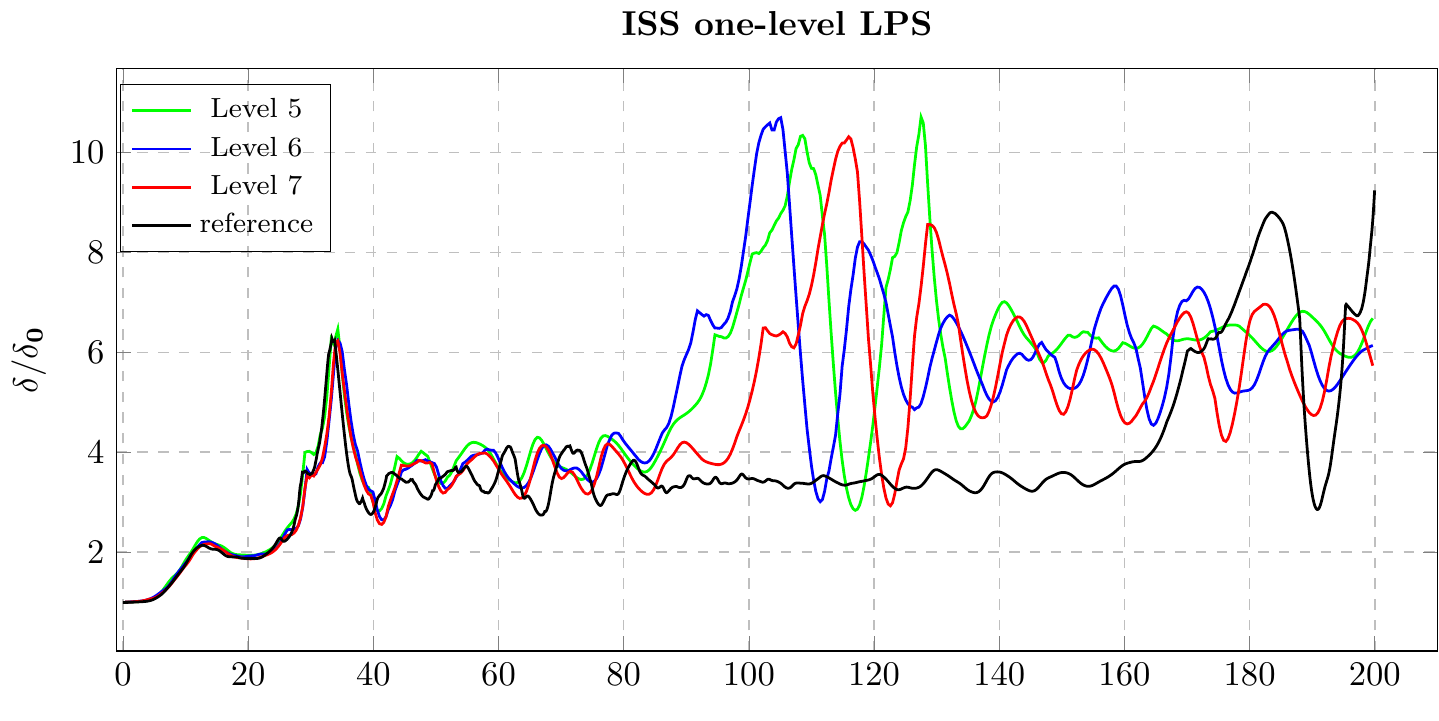}
	\includegraphics[scale=0.5]{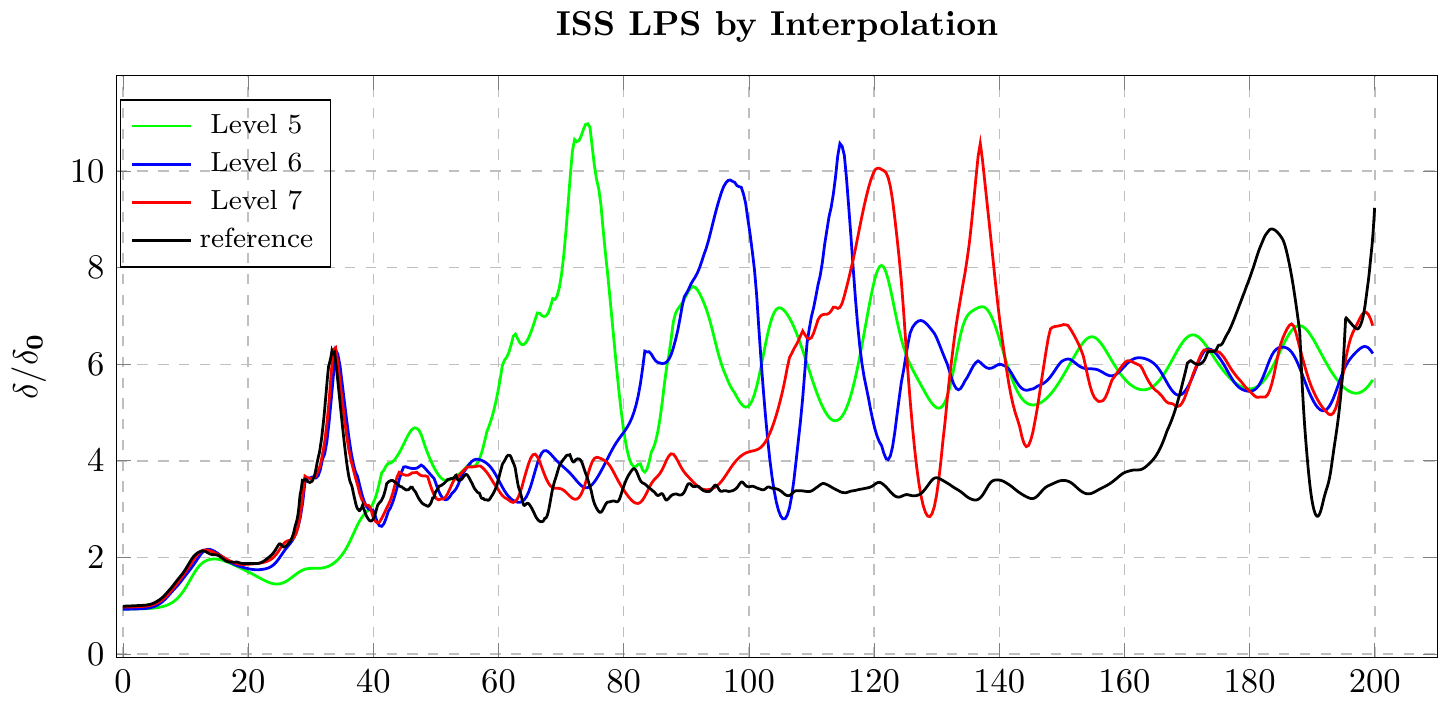}
	\caption{Temporal evolution of vorticity thickness with ISS--FE: RB-VMS (top left), SUPG (top right), 
		one-level LPS (bottom left), and LPS by interpolation (bottom right), on different mesh refinement 
		levels, $\Delta t=0.0125$. \label{fig:vort_is_0.0125}}
\end{figure}

\begin{figure}[t!]
	\centering
	\includegraphics[scale=0.5]{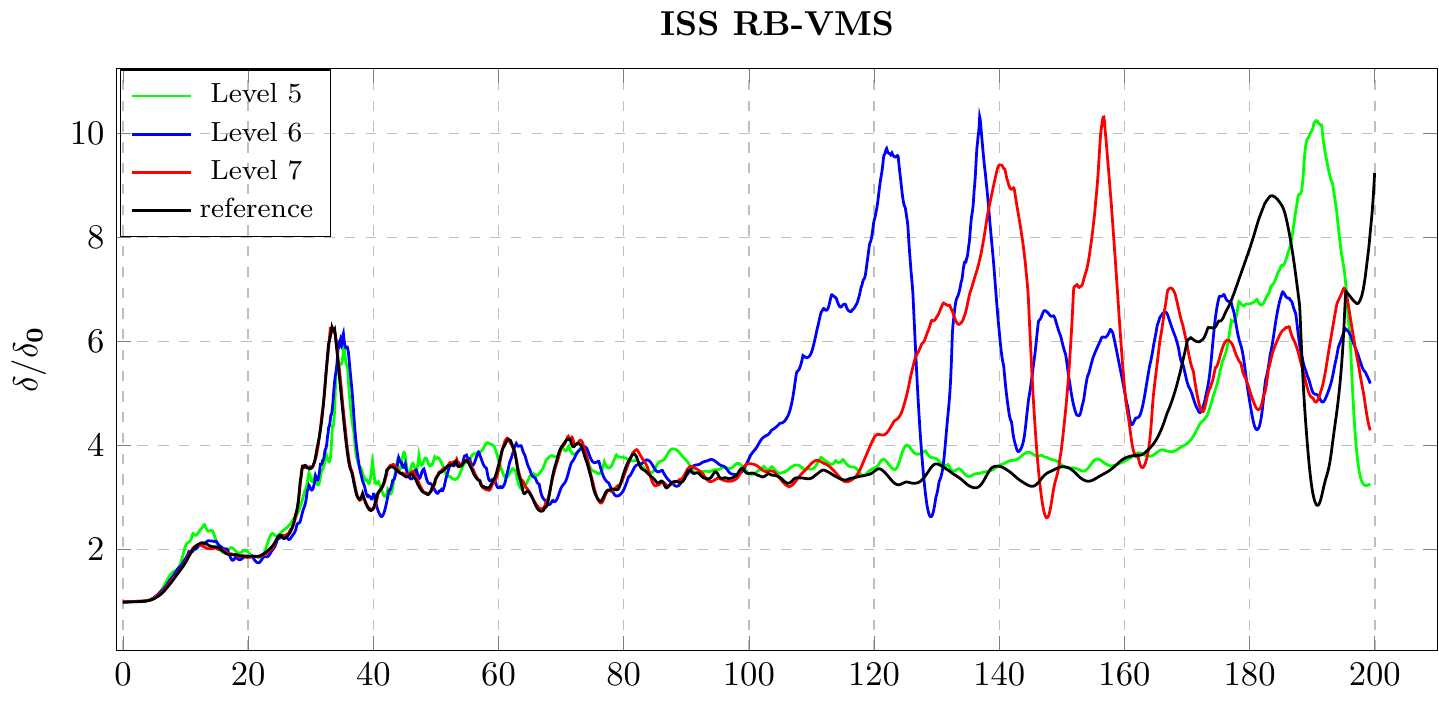}
	\includegraphics[scale=0.5]{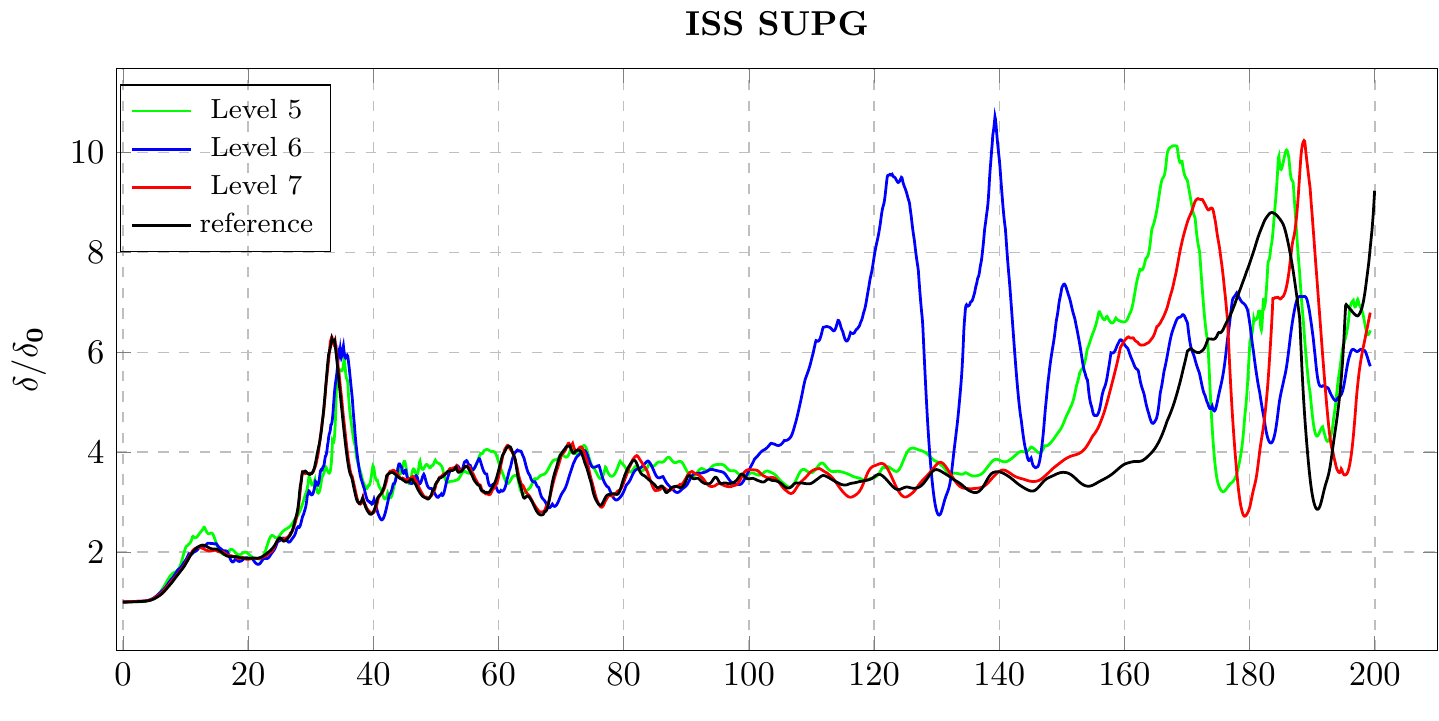}\\
	\includegraphics[scale=0.5]{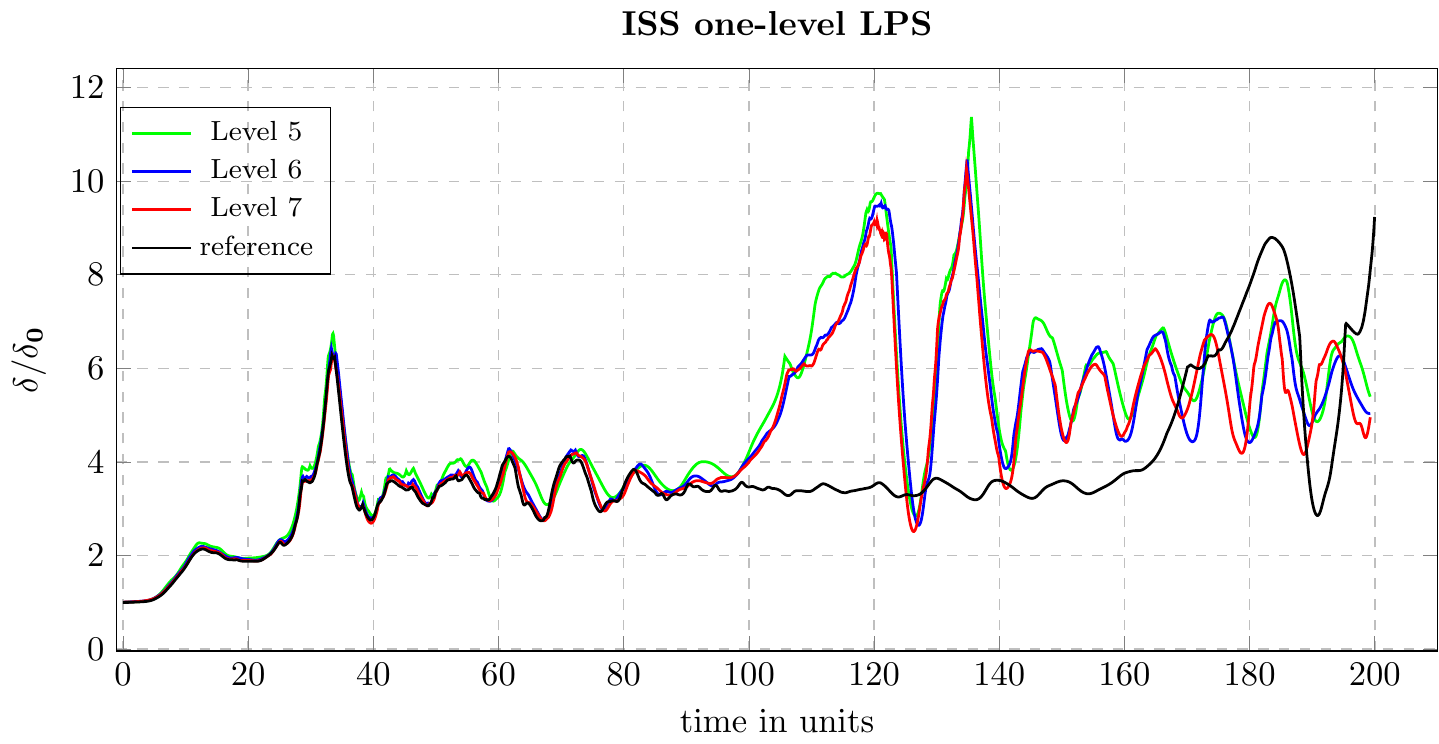}
	\includegraphics[scale=0.5]{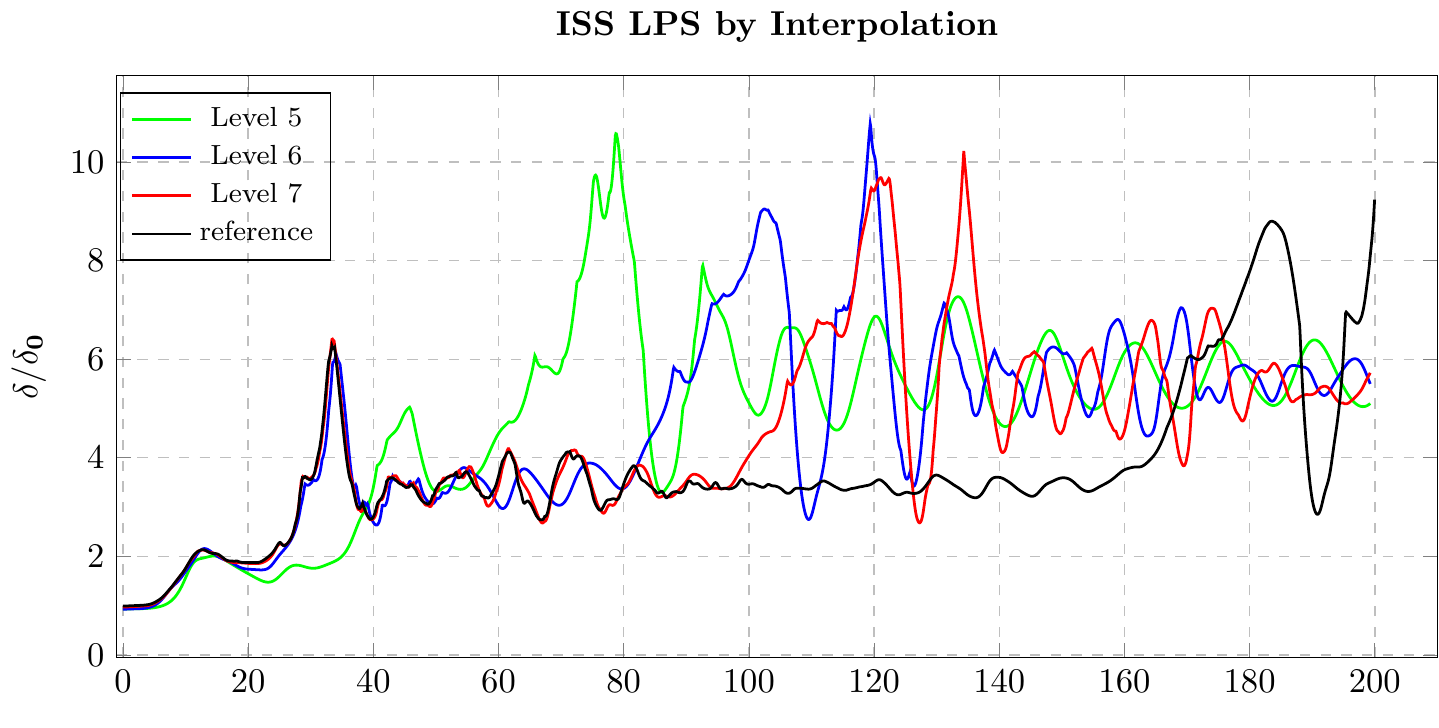}
	\caption{Temporal evolution of vorticity thickness with ISS--FE: RB-VMS (top left), SUPG (top right), 
		one-level LPS (bottom left), and LPS by interpolation (bottom right), on different mesh refinement 
		levels, $\Delta t=0.003125$. \label{fig:vort_is_0.003125}}
\end{figure}

The temporal evaluation of the relative vorticity thickness \(\delta/\delta_0\) for all methods on different 
refinement levels is displayed in Figures~\ref{fig:vort_eo_0.0125} and \ref{fig:vort_eo_0.003125} using EO FE and in 
Figures~\ref{fig:vort_is_0.0125} and \ref{fig:vort_is_0.003125} using ISS FE. The computational results
with time step lengths \(\Delta t_1\) and \(\Delta t_2\) are compared here. The formation of succeeding 
peaks in the evolution of the relative vorticity thickness corresponds to the pairing process of the eddies. 
For the reference solution, the local maximum \(\delta/\delta_0  = 6.2\) at \(t=33.5\bar{t}\) indicates 
the pairing of two eddies from four and this is very much comparable with the first pairing in~\cite{SL17}. 
Comparing the relative vorticity thickness computed with all 
stabilization schemes on different refinement levels clearly indicates that the first pairing occurs
at the same time except for the coarse grid simulation. After that, the relative 
vorticity oscillates until the next pairing of eddies occurs. The pairing of the final eddy happens 
somehow at different time for different stabilization methods. Considering the red-curves (refinement level 7)
in Figure~\ref{fig:vort_eo_0.0125} for \(\Delta t_1\), we clearly see that, for the RB-VMS method, the merging of the 
two secondary eddies into one vortex starts later, at time \(t=160\bar{t}\), in agreement with the finest simulation 
in~\cite{SL17}. On the other side, this is slightly anticipated for SUPG method (\(t=140\bar{t}\)), and one-level and 
interpolation based LPS methods (\(t=120\bar{t}\)). On the other hand, for the small time step length \(\Delta t_2\) 
in Figure~\ref{fig:vort_eo_0.003125}, in the case of 
RB-VMS method, the last pairing occurs a bit earlier than for the large time step length \(\Delta t_2\). However, the 
SUPG method behaves very similar to the reference solution, for which last pairing occurs at \(t=180\bar{t}\), even 
on the coarser meshes (level 6). The two versions of the LPS methods are almost comparable to each other, except for
the coarsest mesh (level 5). Thus, note that the development of the vorticity thickness strongly depends on the mesh refinement. In particular,
the last pairing process, where two secondary eddies merge to become one, is very sensible
with respect to how accurate the simulation is. However, the actual values of the 
amplitudes of the various peaks are almost identical for all refinement levels, and in agreement both with our reference solution and results in~\cite{SL17}.

A similar conclusion can be made for ISS FE by comparing with EO FE case, 
see Figures~\ref{fig:vort_is_0.0125} and \ref{fig:vort_is_0.003125}. 
One can see that the mesh refinement and time step lengths have again a considerable 
influence on the temporal development of the vorticity thickness, but the values
of the amplitude are almost identical. In addition, a mesh convergence can be seen 
for the one-level variant of LPS method in the case of small time step length. 
Altogether, the EO SUPG method with small time step length is 
superior to all other methods, since almost approach the very fine reference solution on relatively coarse grids. 
However, the EO RB-VMS method and ISS SUPG method with small time step also perform quite well, being almost in agreement with the finest simulation in~\cite{SL17}. 

\subsection{Kinetic Energy}\label{subsec:EKin}
\begin{figure}[t!]
	\centering
	\includegraphics[scale=0.5]{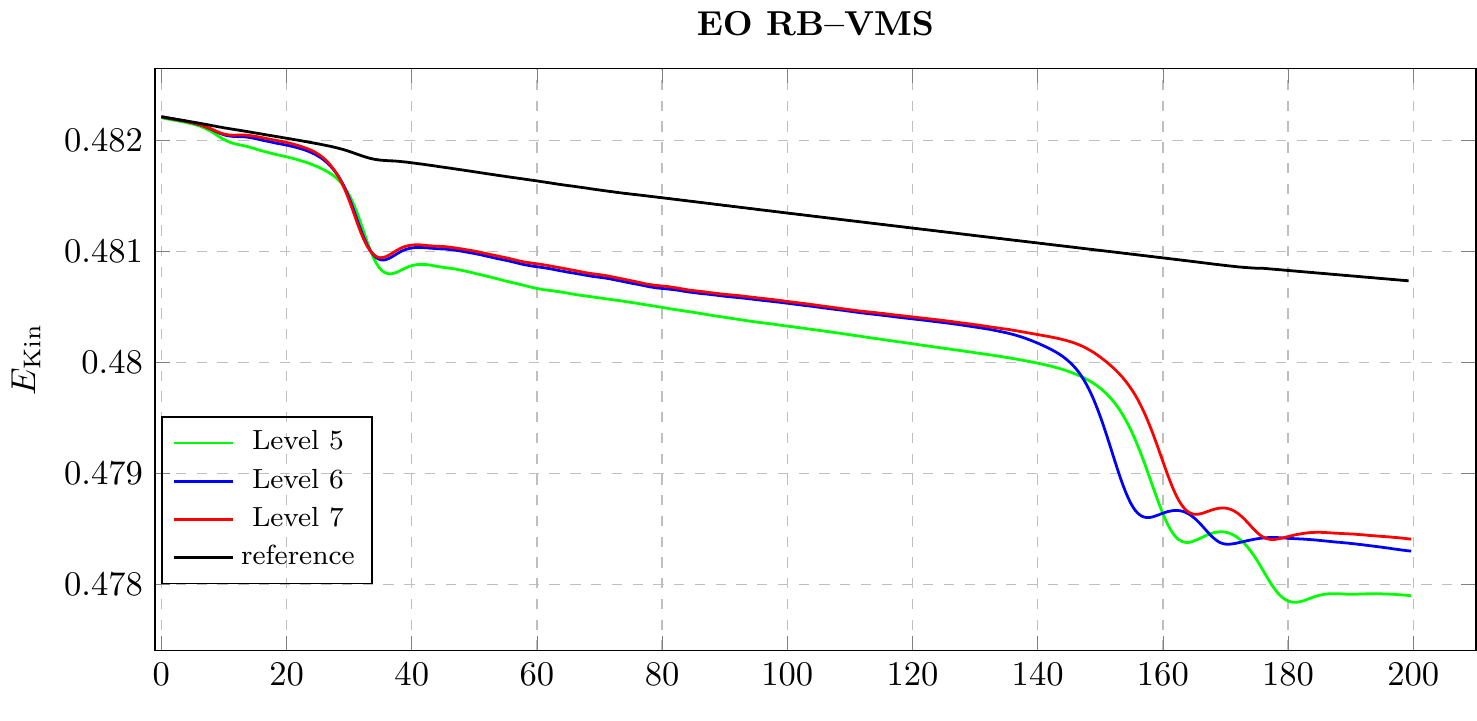}
	\includegraphics[scale=0.5]{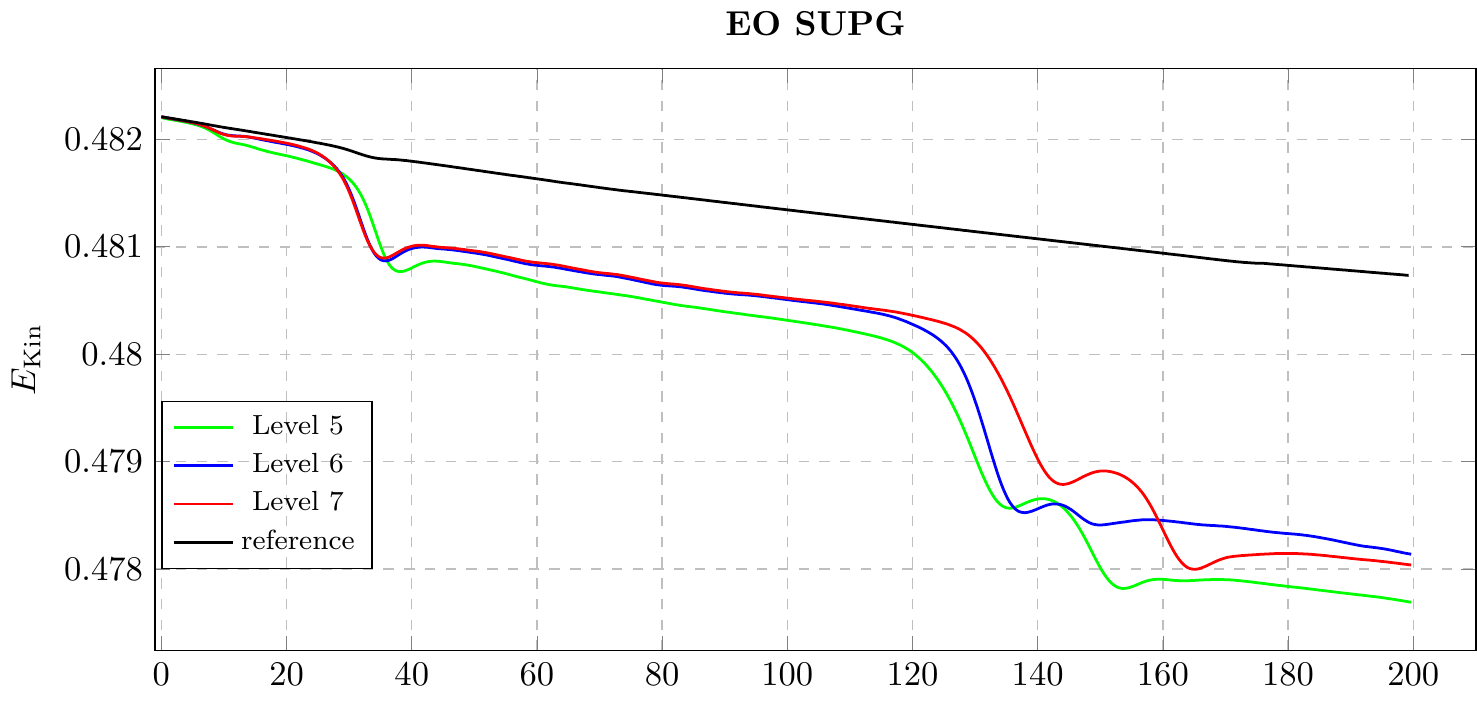}\\
	\includegraphics[scale=0.5]{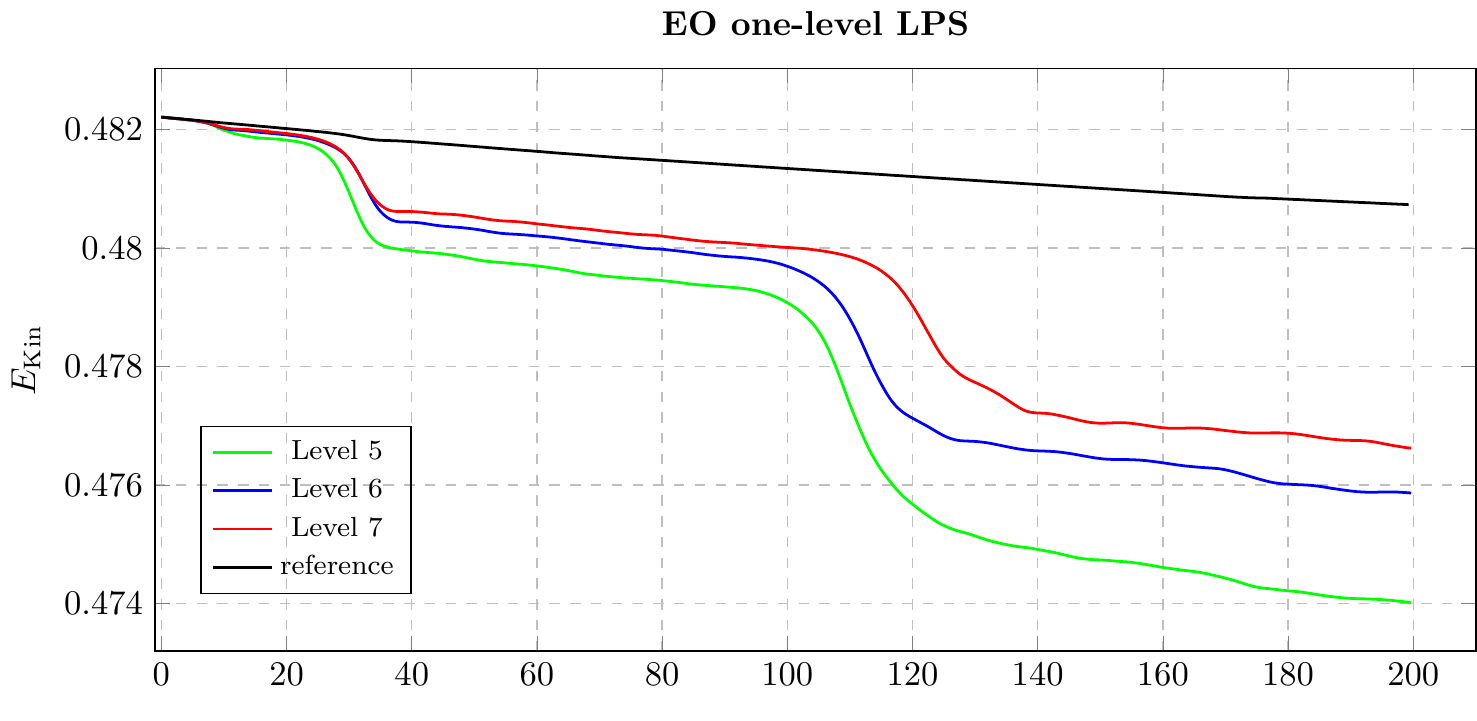}
	\includegraphics[scale=0.5]{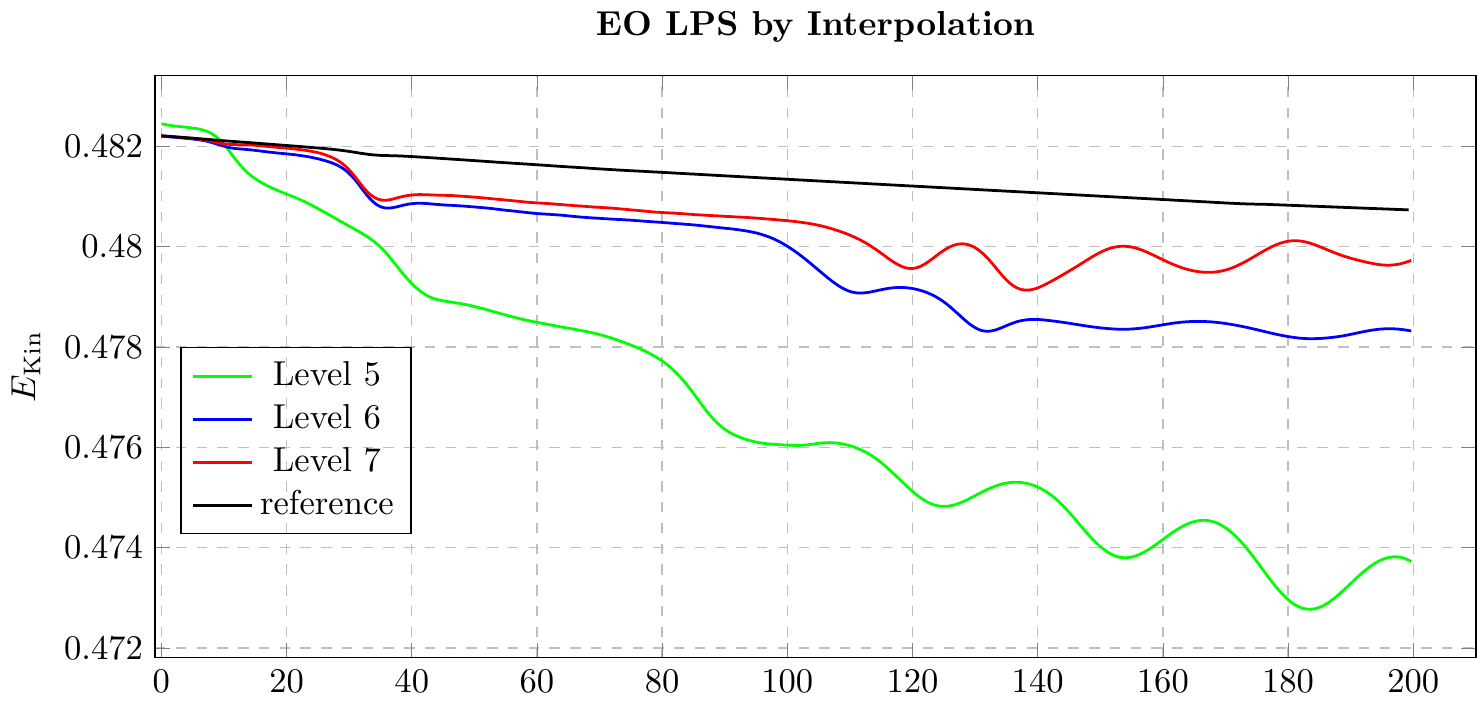}
	\caption{Temporal evolution of kinetic energy with EO--FE: RB-VMS (top left), SUPG (top right), 
		one-level LPS (bottom left), and LPS by interpolation (bottom right), on different mesh 
		refinement levels, $\Delta t=0.0125$. \label{fig:ke_eo_0.0125}}
\end{figure}

\begin{figure}[t!]
	\centering
\includegraphics[scale=0.5]{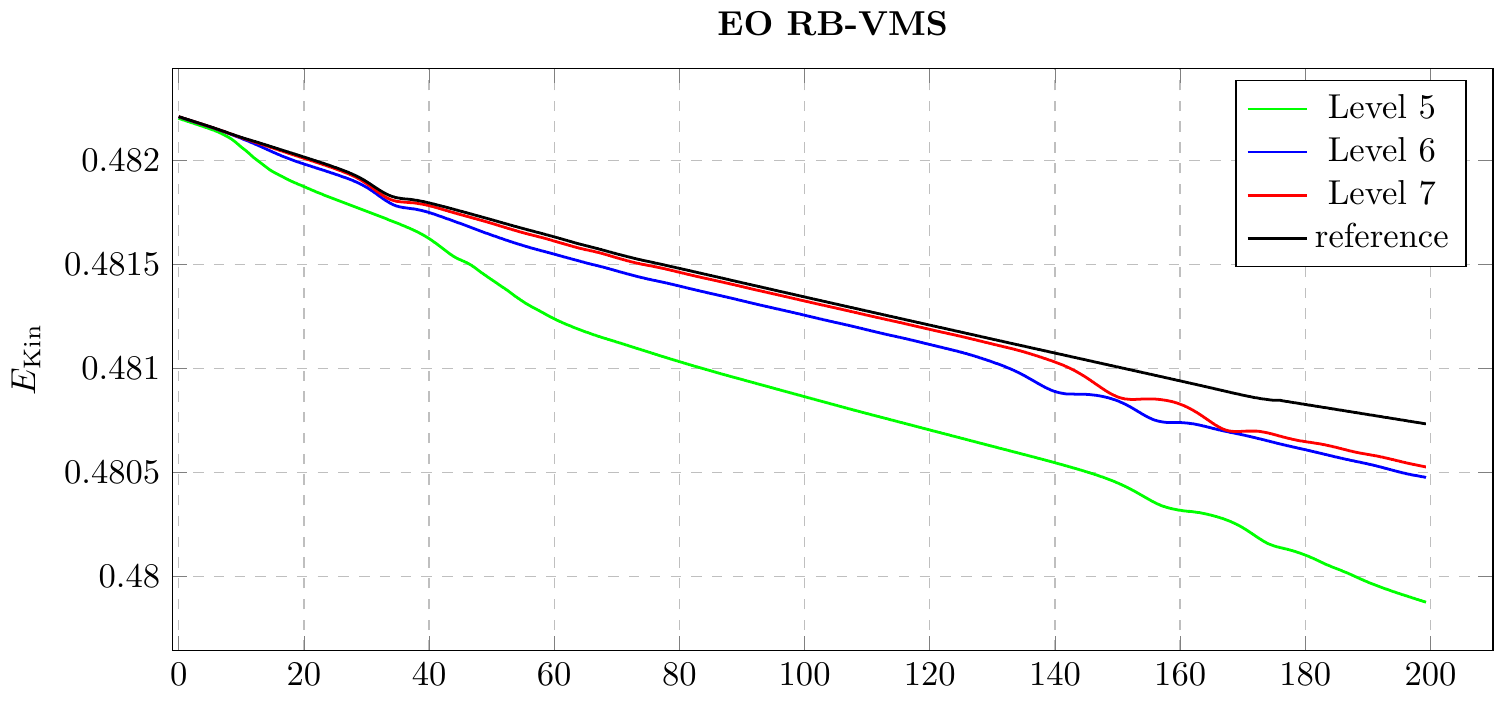}
\includegraphics[scale=0.5]{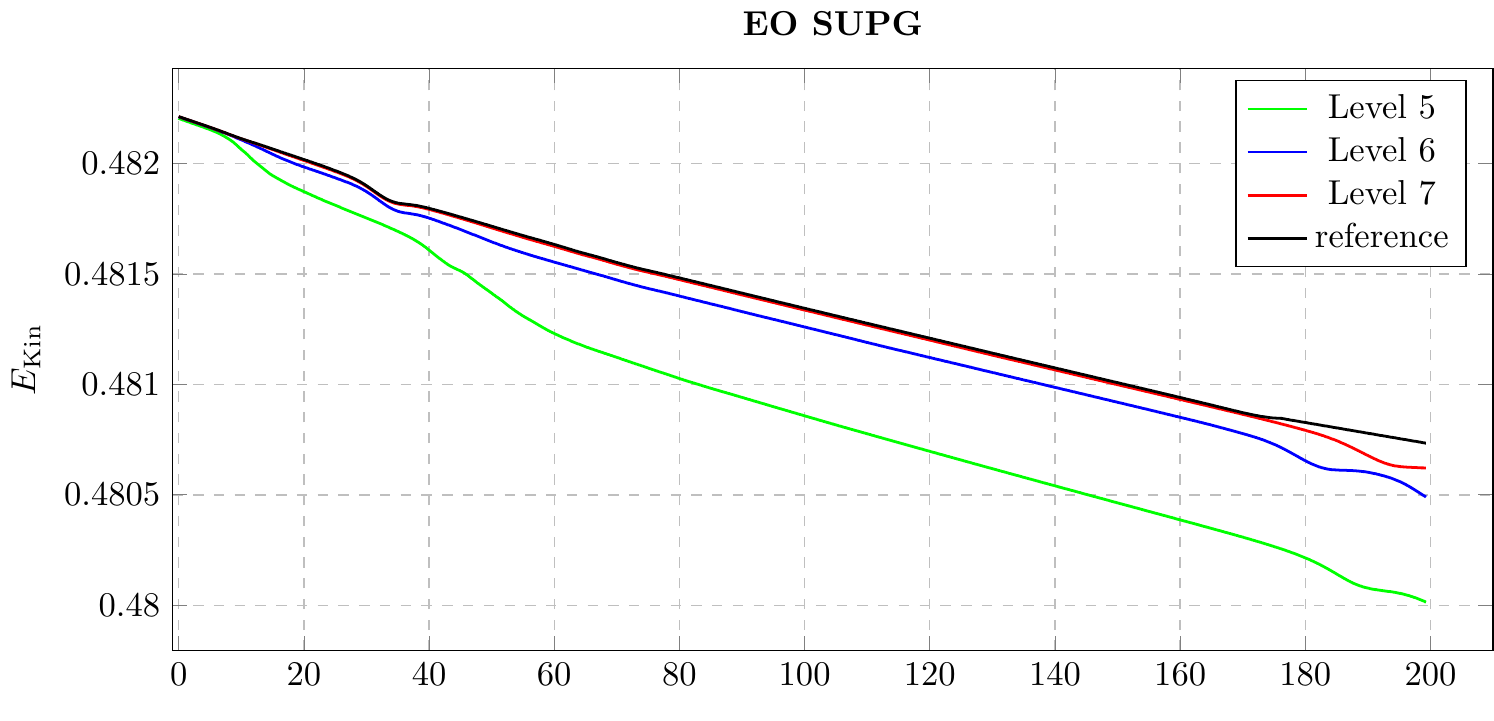}\\
\includegraphics[scale=0.5]{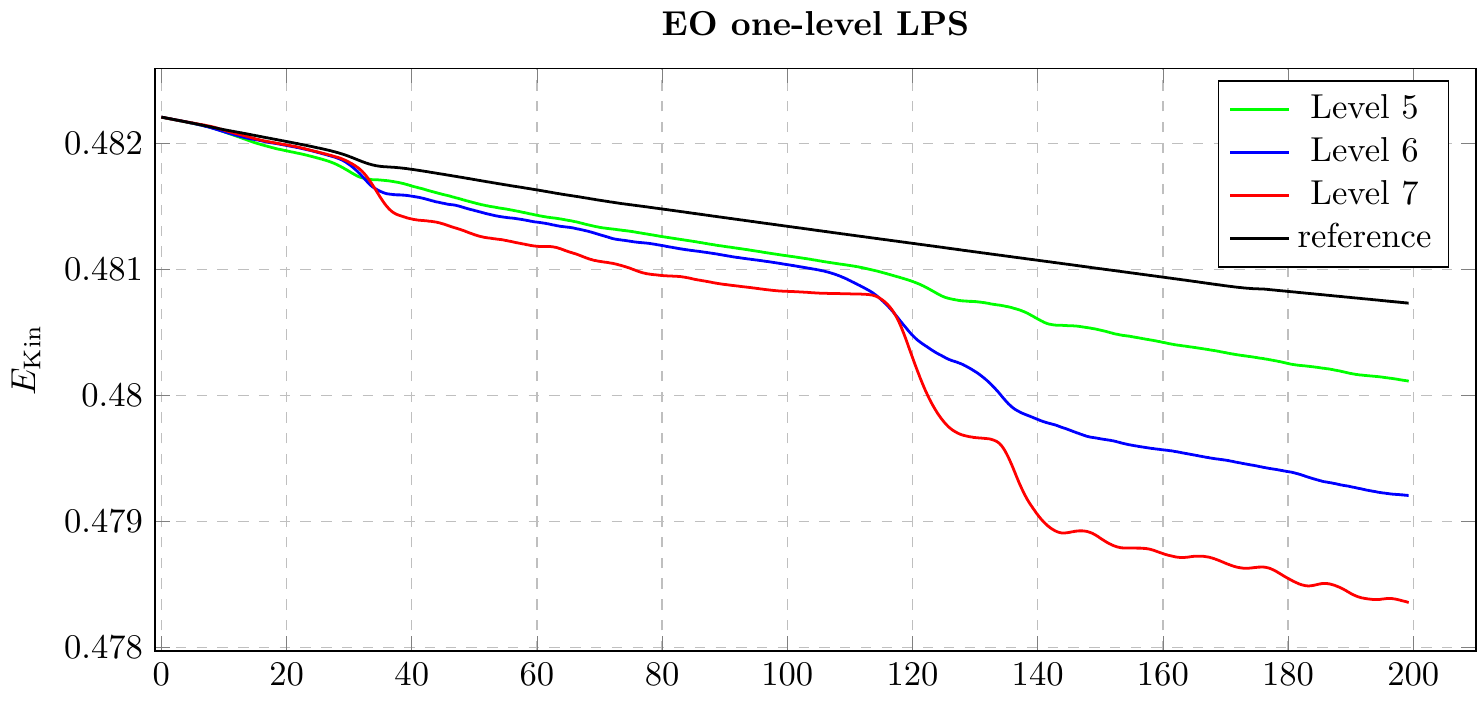}
\includegraphics[scale=0.5]{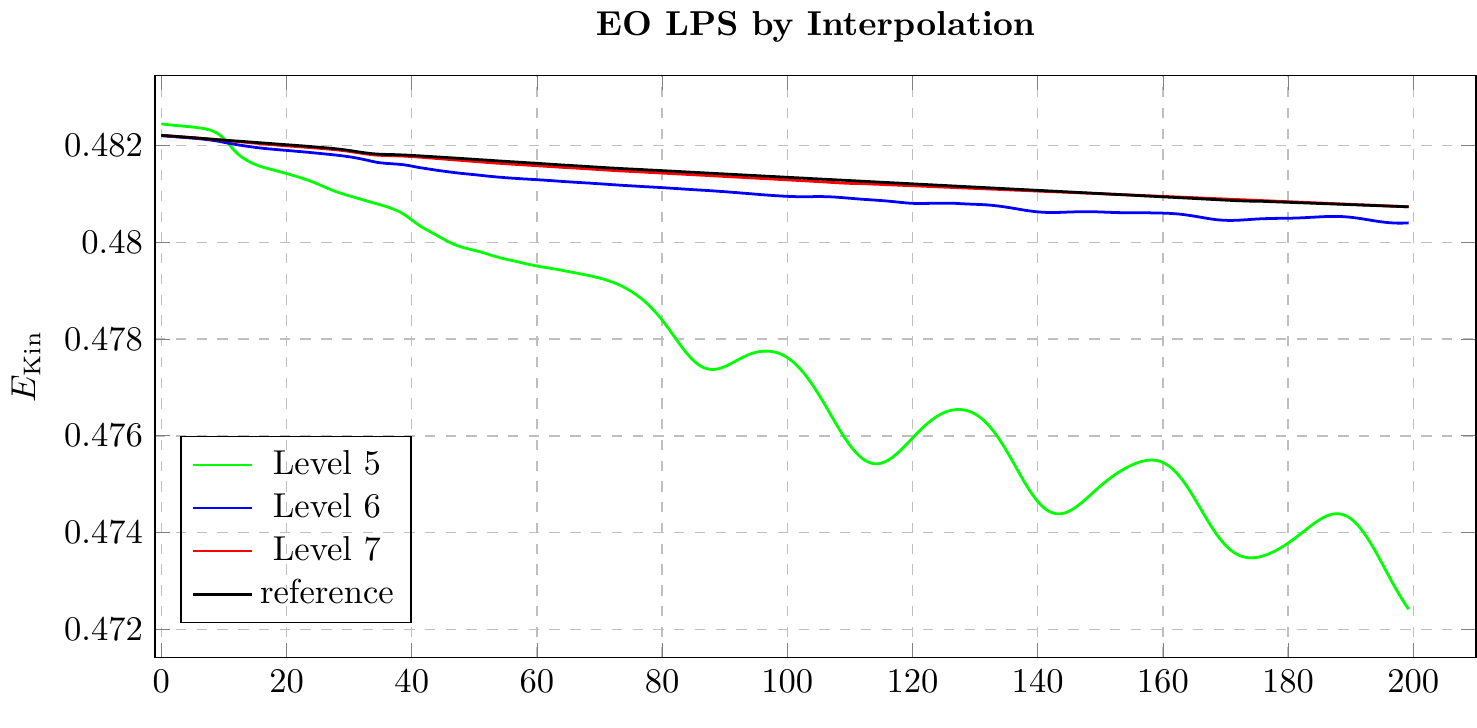}
	\caption{Temporal evolution of kinetic energy with EO--FE: RB-VMS (top left), SUPG (top right), 
	one-level LPS (bottom left), and LPS by interpolation (bottom right), on different mesh refinement 
	levels, $\Delta t=0.003125$. \label{fig:ke_eo_0.003125}}
\end{figure}

\begin{figure}[t!]
	\centering
		\includegraphics[scale=0.5]{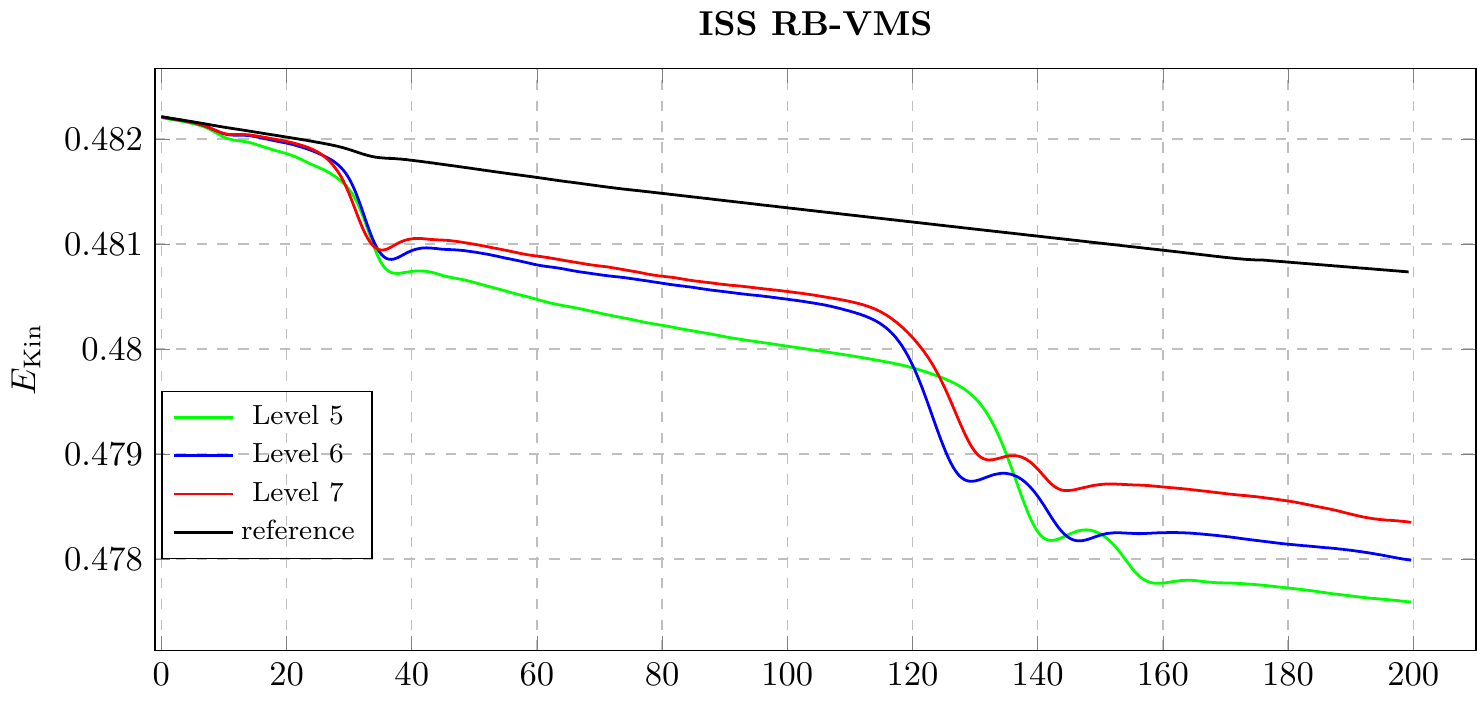}
		\includegraphics[scale=0.5]{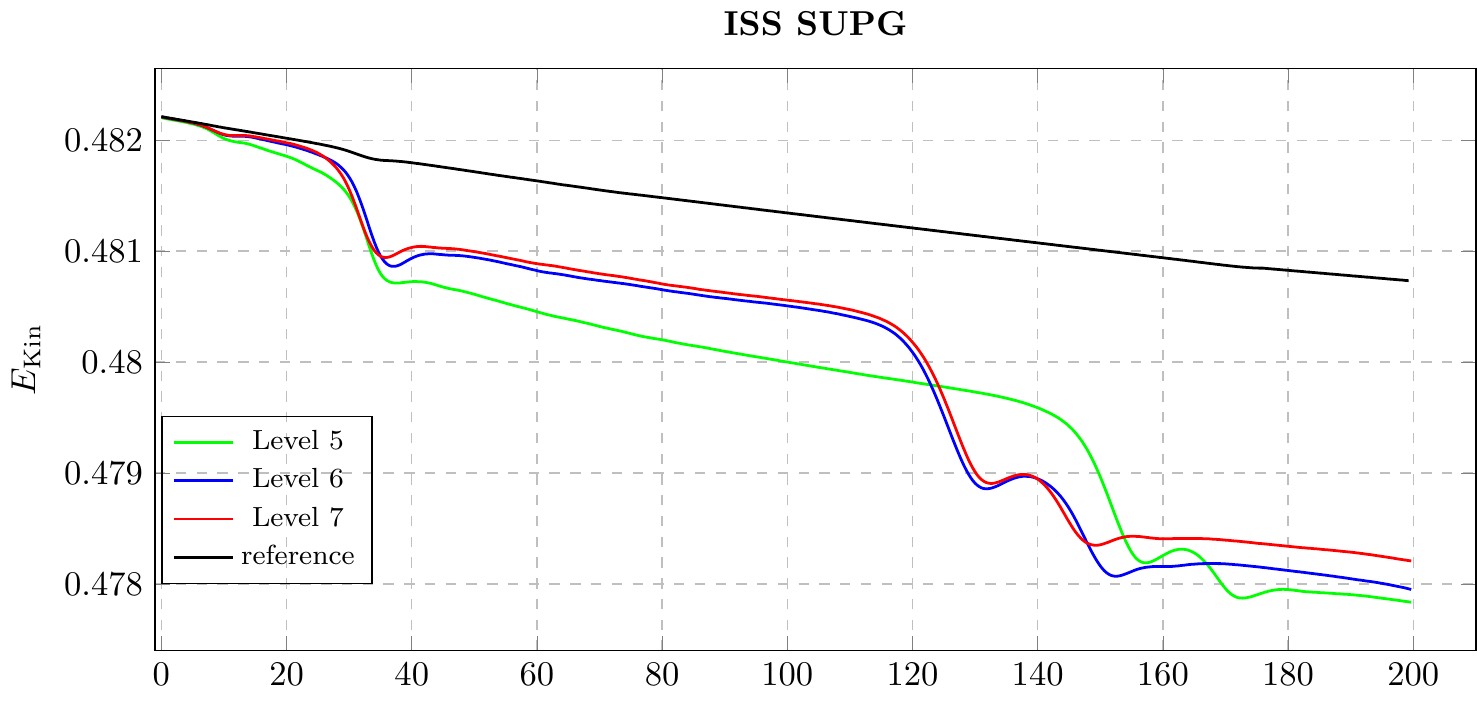}\\
		\includegraphics[scale=0.5]{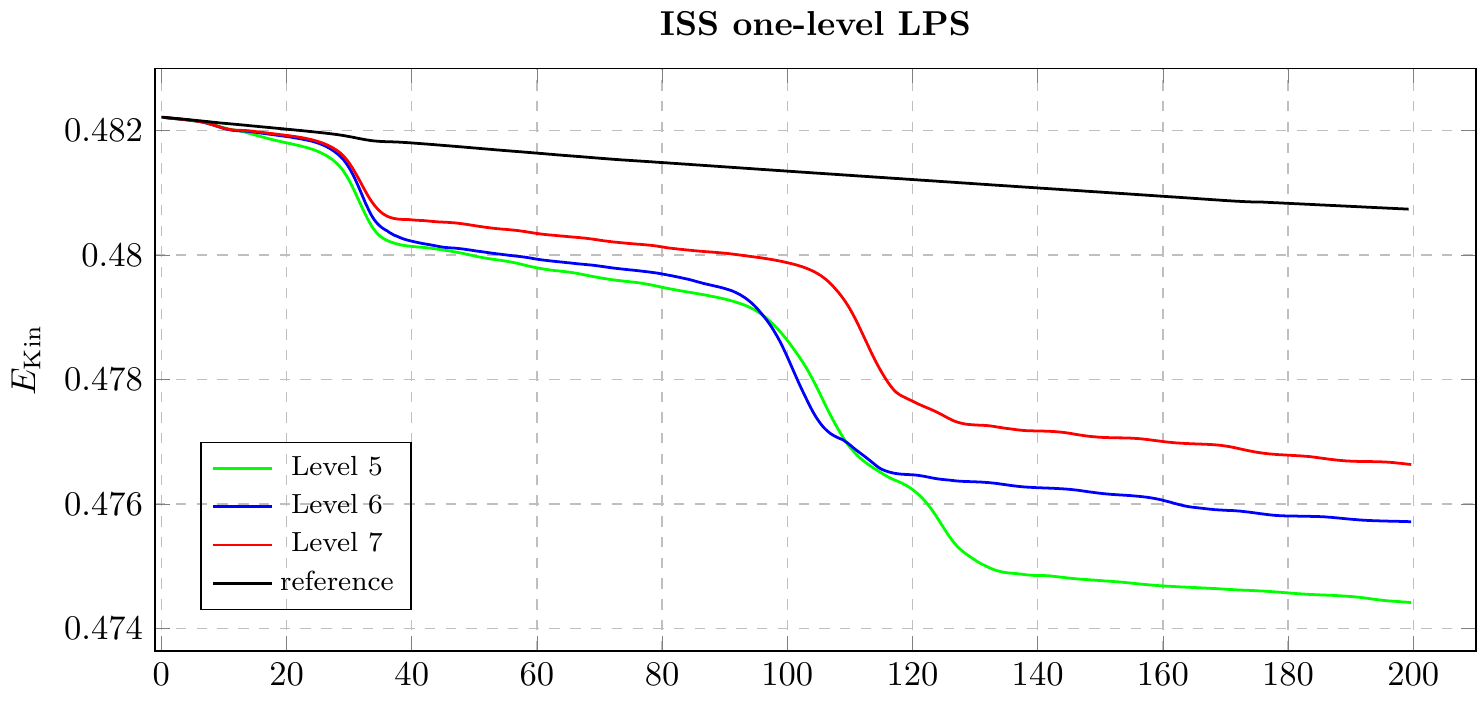}
		\includegraphics[scale=0.5]{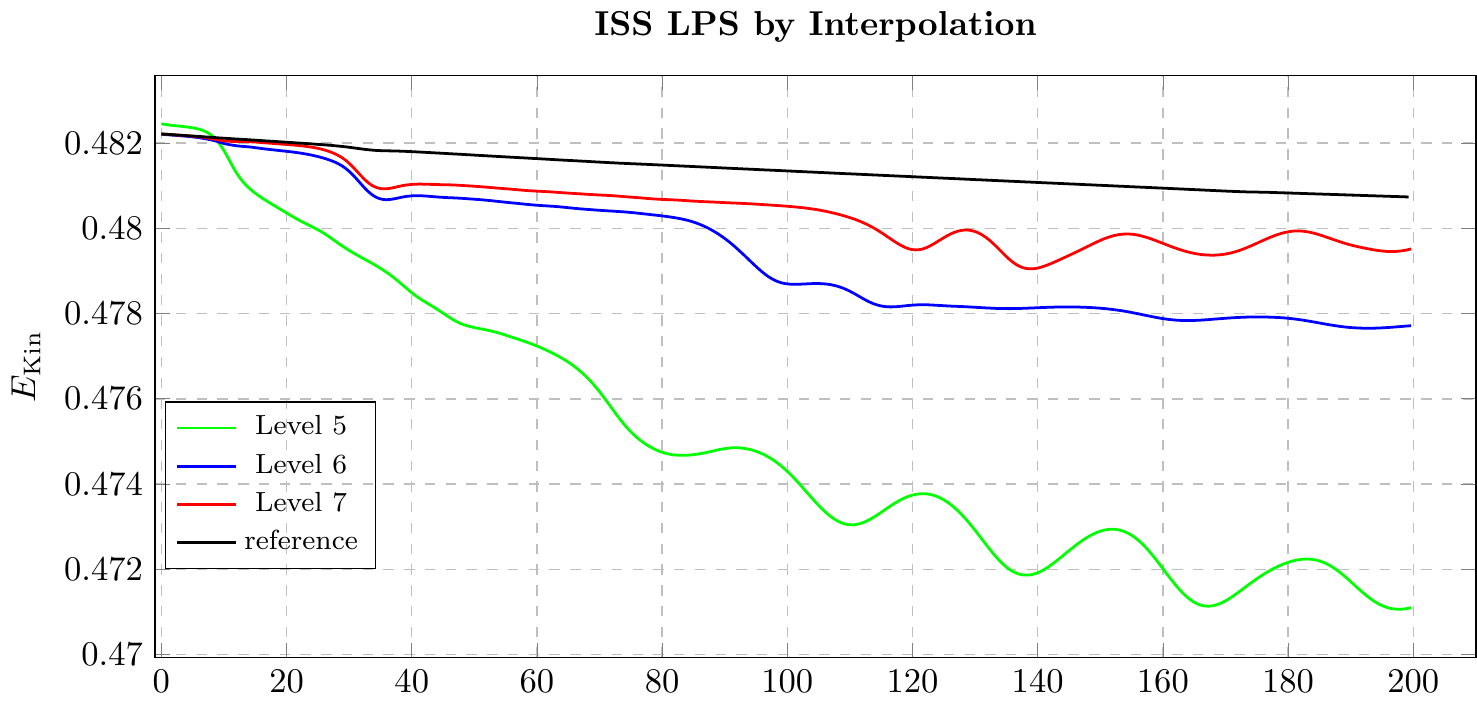}
	\caption{Temporal evolution of kinetic energy with ISS--FE: RB-VMS (top left), SUPG (top right), 
	one-level LPS (bottom left), and LPS by interpolation (bottom right), on different mesh refinement
	levels, $\Delta t=0.0125$. \label{fig:ke_is_0.0125}}
\end{figure}

\begin{figure}[t!]
	\centering
		\includegraphics[scale=0.5]{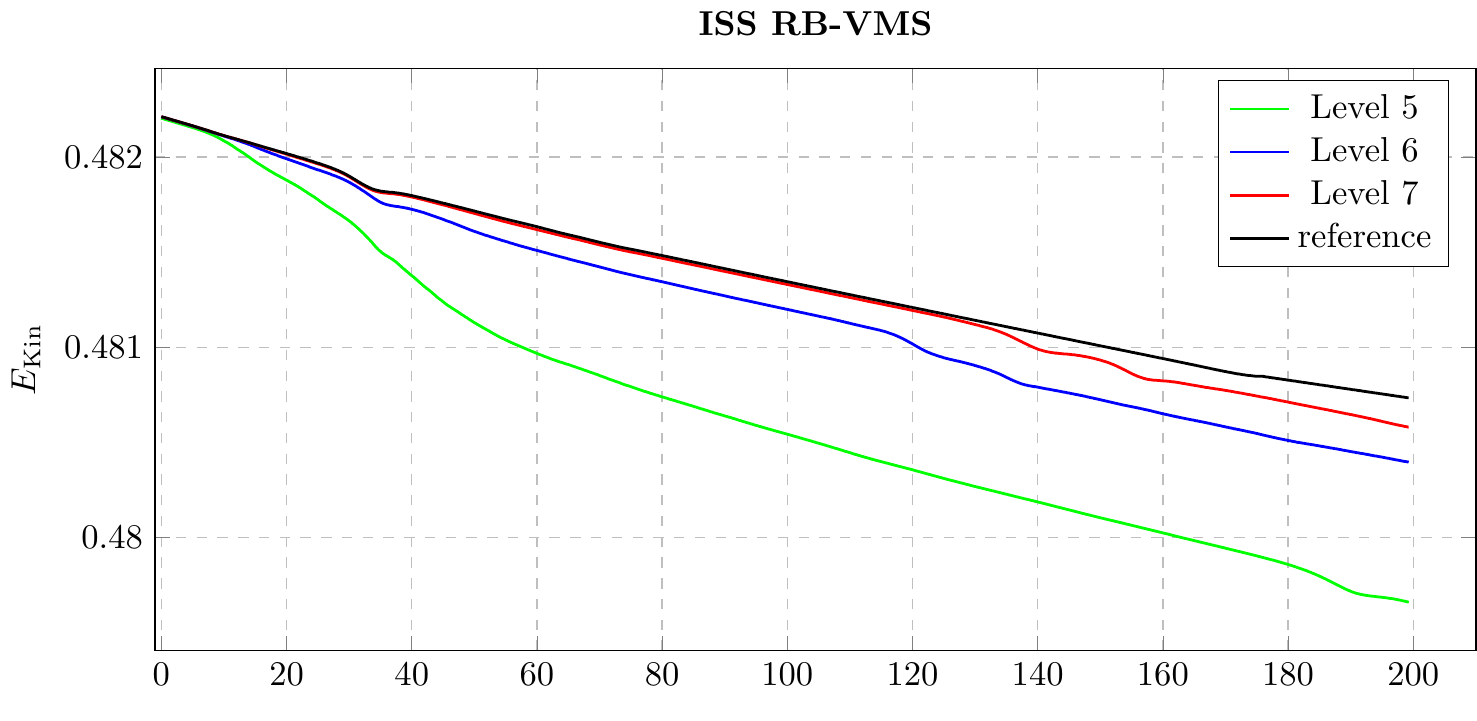}
		\includegraphics[scale=0.5]{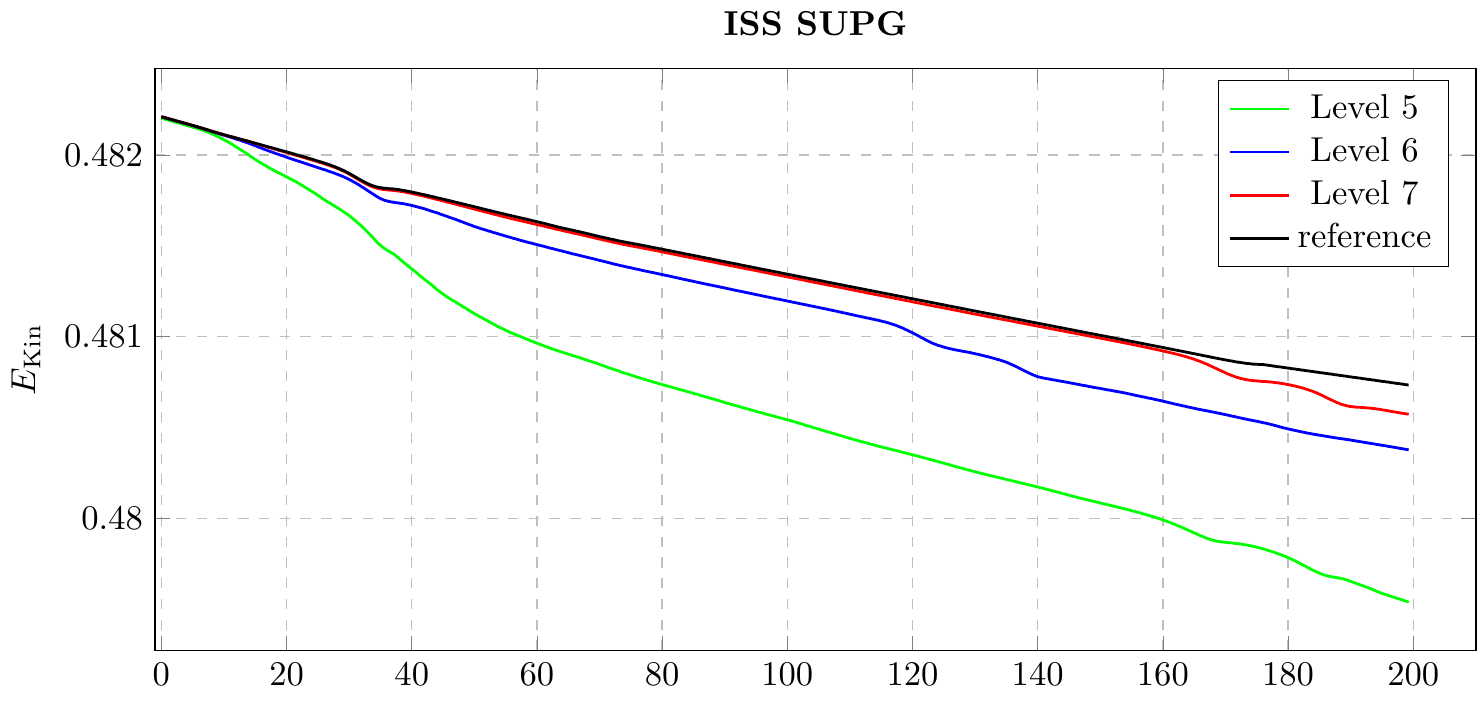}\\
		\includegraphics[scale=0.5]{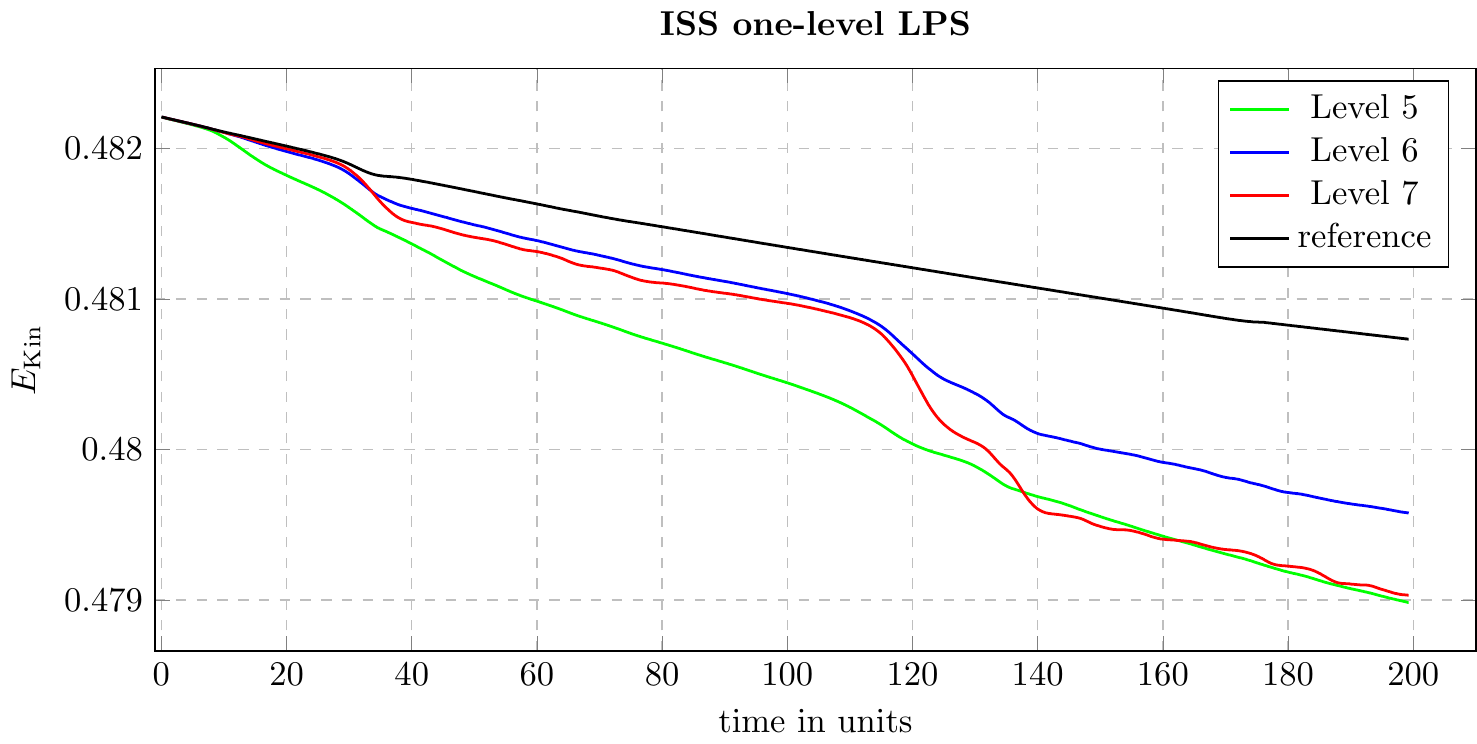}
		\includegraphics[scale=0.5]{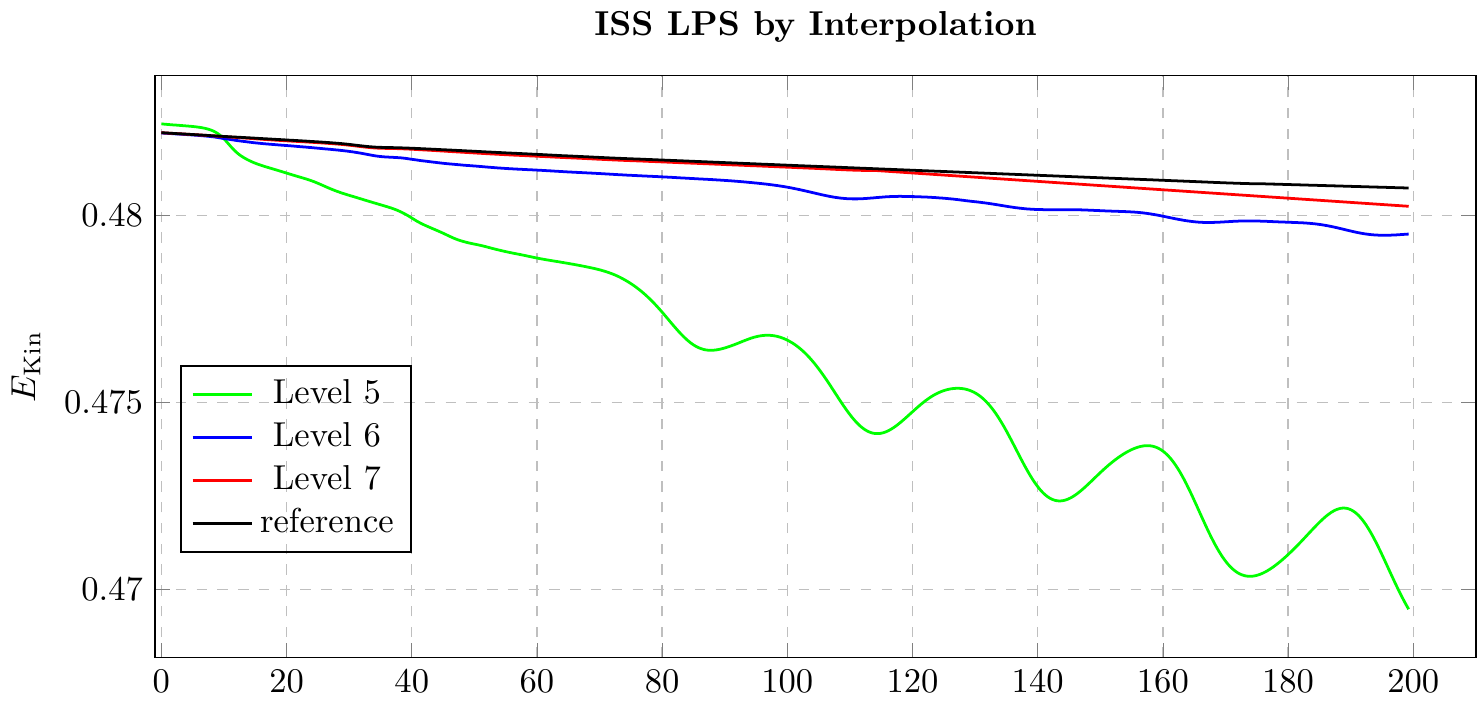}
	\caption{Temporal evolution of kinetic energy with ISS--FE: RB-VMS (top left), SUPG (top right), 
	one-level LPS (bottom left), and LPS by interpolation (bottom right), on different mesh refinement
	levels, $\Delta t=0.003125$. \label{fig:ke_is_0.003125}}
\end{figure}

The temporal evolution of the total kinetic energy for all considered method will be discussed 
in this section. Figures~\ref{fig:ke_eo_0.0125}--\ref{fig:ke_is_0.003125} 
presents the evolution of the total kinetic energy for EO and ISS pair of FE on different refinement 
levels and for different time step lengths. 
In principal,
an evolution exhibiting a monotone decaying total amount of kinetic energy has to be
physically expected, since the initial velocity distribution is subject to a non-zero viscosity,
and no additional energy input is provided.
First, consider the case of the large time step length 
\(\Delta t_1\). The results presented in Figure~\ref{fig:ke_eo_0.0125} shows some increase
in the total kinetic energy in certain time intervals for all methods and grid levels, except for the one-level LPS 
method, being oscillations more evident in the case of the LPS by interpolation method. 
This is clearly not
a physical behavior and there is no mechanism in this problem which could excite it.
In our opinion, this numerical experience suggests that RB-VMS, SUPG, and especially LPS by interpolation methods, 
when coupled with less time consuming semi-implicit discretizations in time, are rather restrictive in terms of 
time step length in order to guarantee stability, that is they would require a finer time step length with respect 
to the one-level LPS method to achieve a correct physical behavior, not influenced by numerical stability issues. 
In any case, note that for all methods, also for the one-level LPS method, the large time step causes an excessive
dissipation (up to a 2\% for the coarsest level), traduced in a higher overall energy loss when compared to the 
corresponding results in~\cite{SL17}, and mesh convergence for the kinetic energy is not achieved for the finer 
meshes. On the contrary, a monotonically decreasing kinetic energy is obtained for all the methods using the small 
time step length \(\Delta t_2\), which can be seen clearly in the Figure~\ref{fig:ke_eo_0.003125}. The only exception 
here is the coarsest level for the LPS by interpolation method, probably due to the fact that for this very simplified
VMS stabilized method this space resolution is to low to guarantee stability. However, for the one-level LPS method, 
mesh convergence is still not reached, which is indeed the case for all other methods on finer grid resolutions, i.e. 
levels 6 and 7. In these cases, the kinetic energy decreases very slowly, around 0.3\%, as in~\cite{SL17}, and results 
are almost comparable on the finest grid to the ones provided by the reference solution, being almost identical for 
SUPG and LPS by interpolation methods (slightly more pronounced differences are noticeable for RB-VMS method). 

Similar conclusions as for EO FE can be drawn for the computational results obtained by using ISS FE, see 
Figures~\ref{fig:ke_is_0.0125} and \ref{fig:ke_is_0.003125}. 

\subsection{Enstrophy}\label{subsec:Enst}
\begin{figure}[t!]
	\centering
		\includegraphics[scale=0.5]{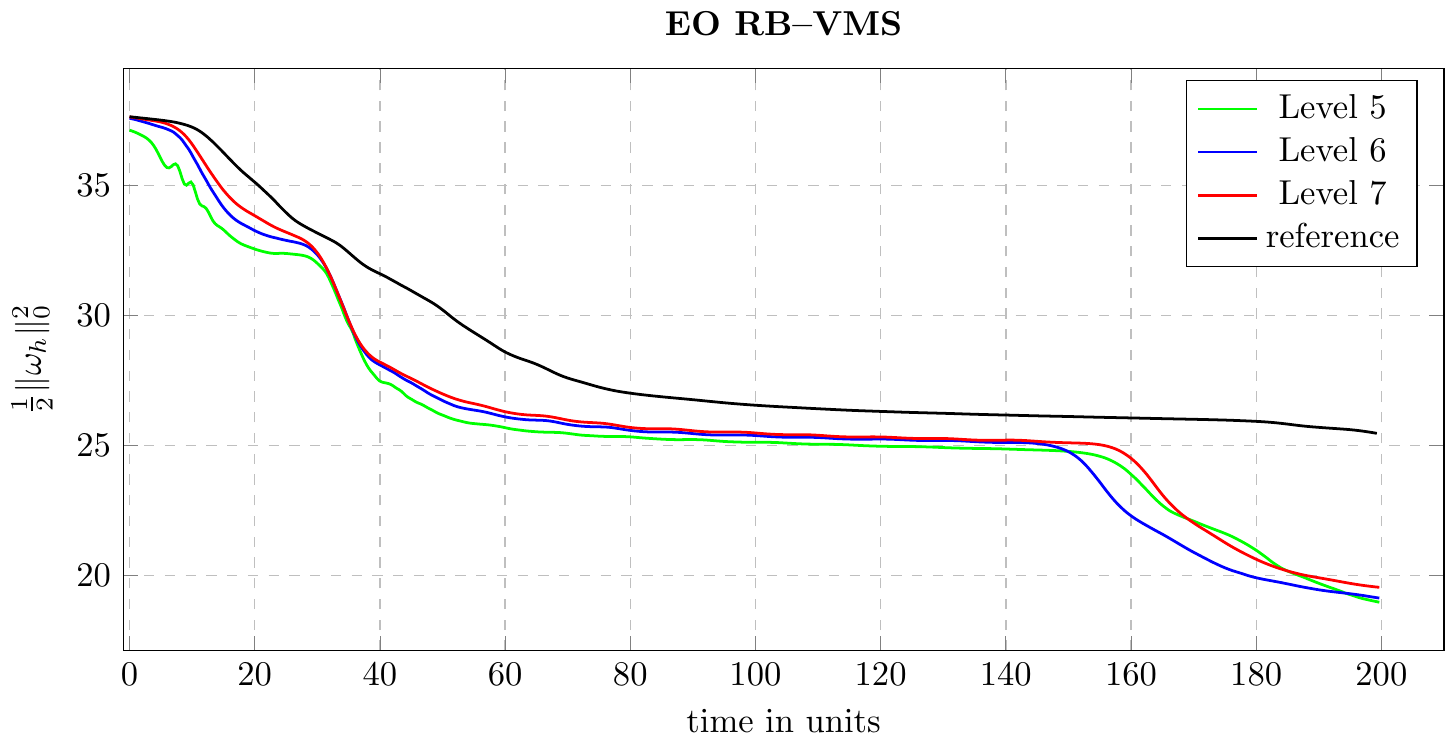}
		\includegraphics[scale=0.5]{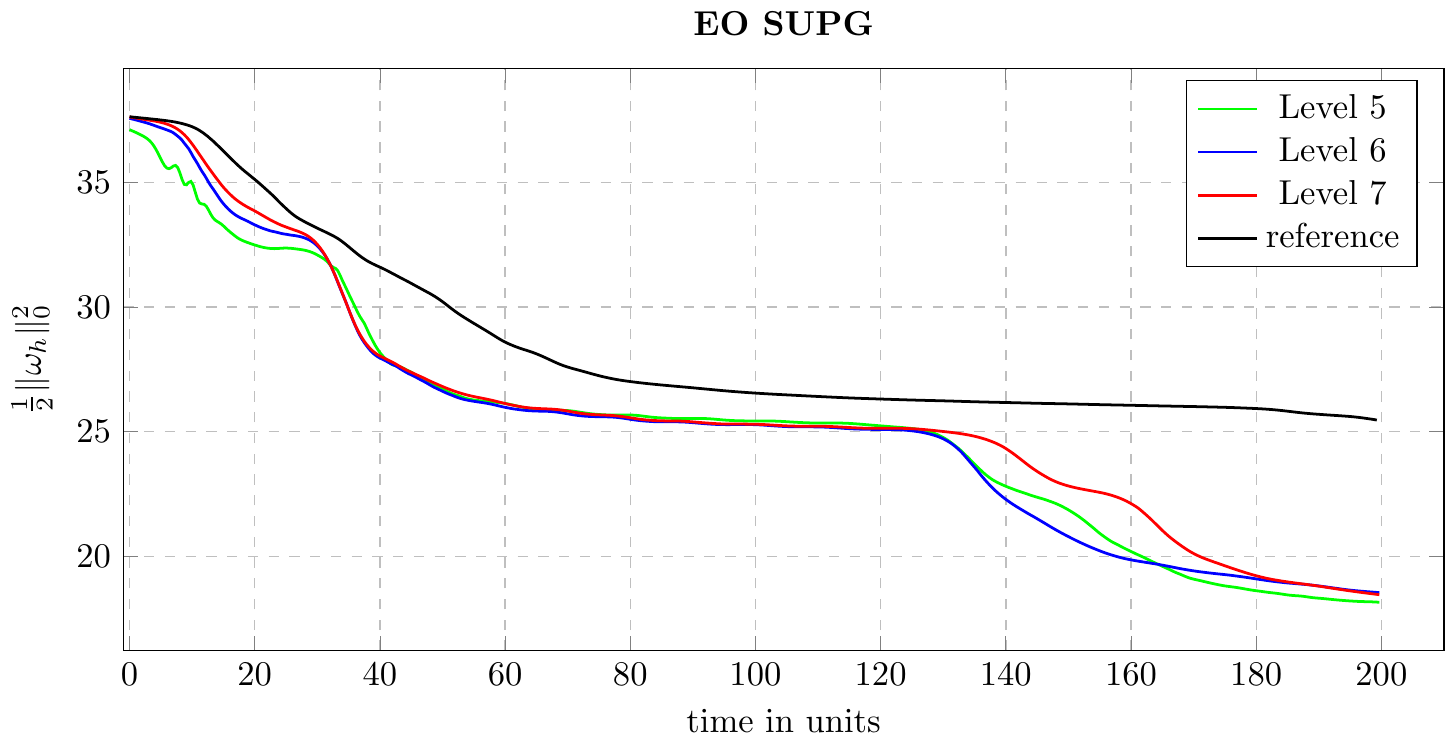}\\
		\includegraphics[scale=0.5]{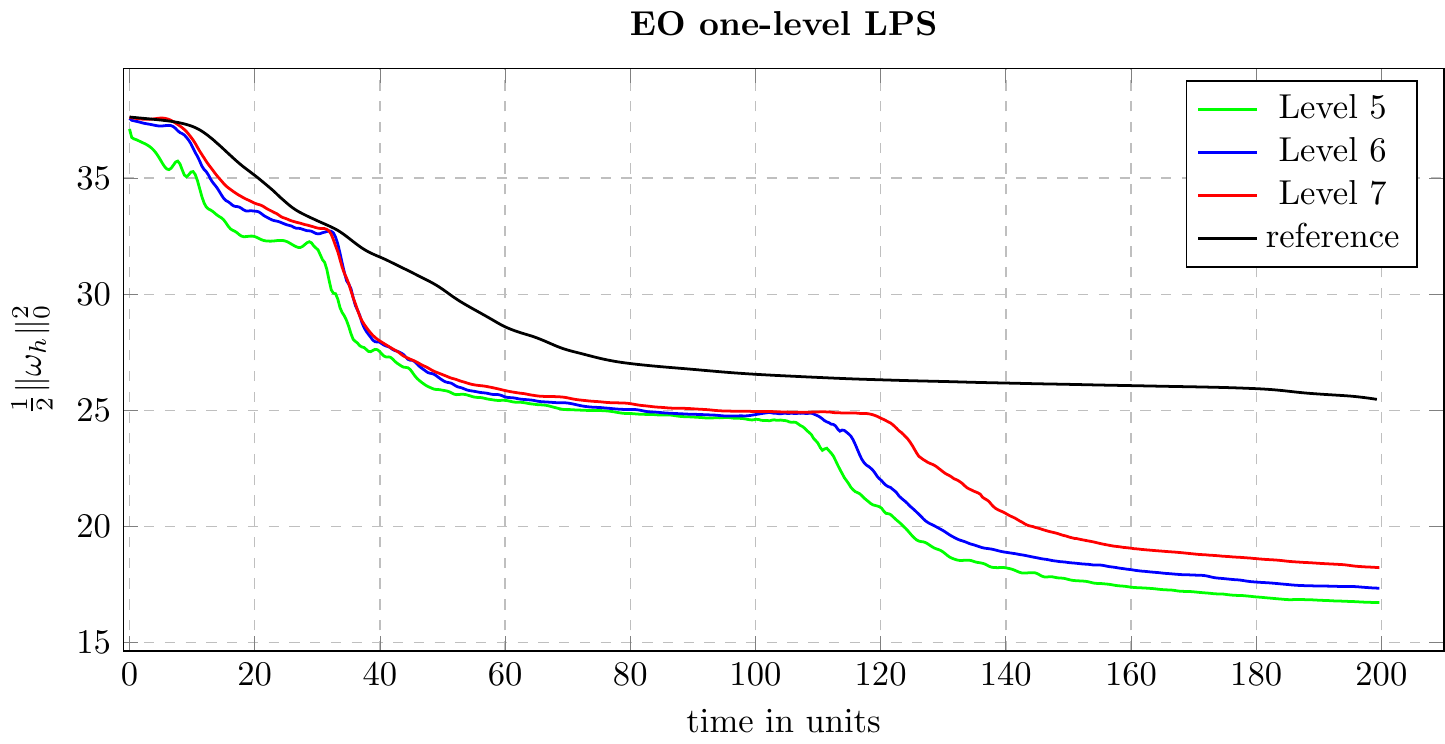}
		\includegraphics[scale=0.5]{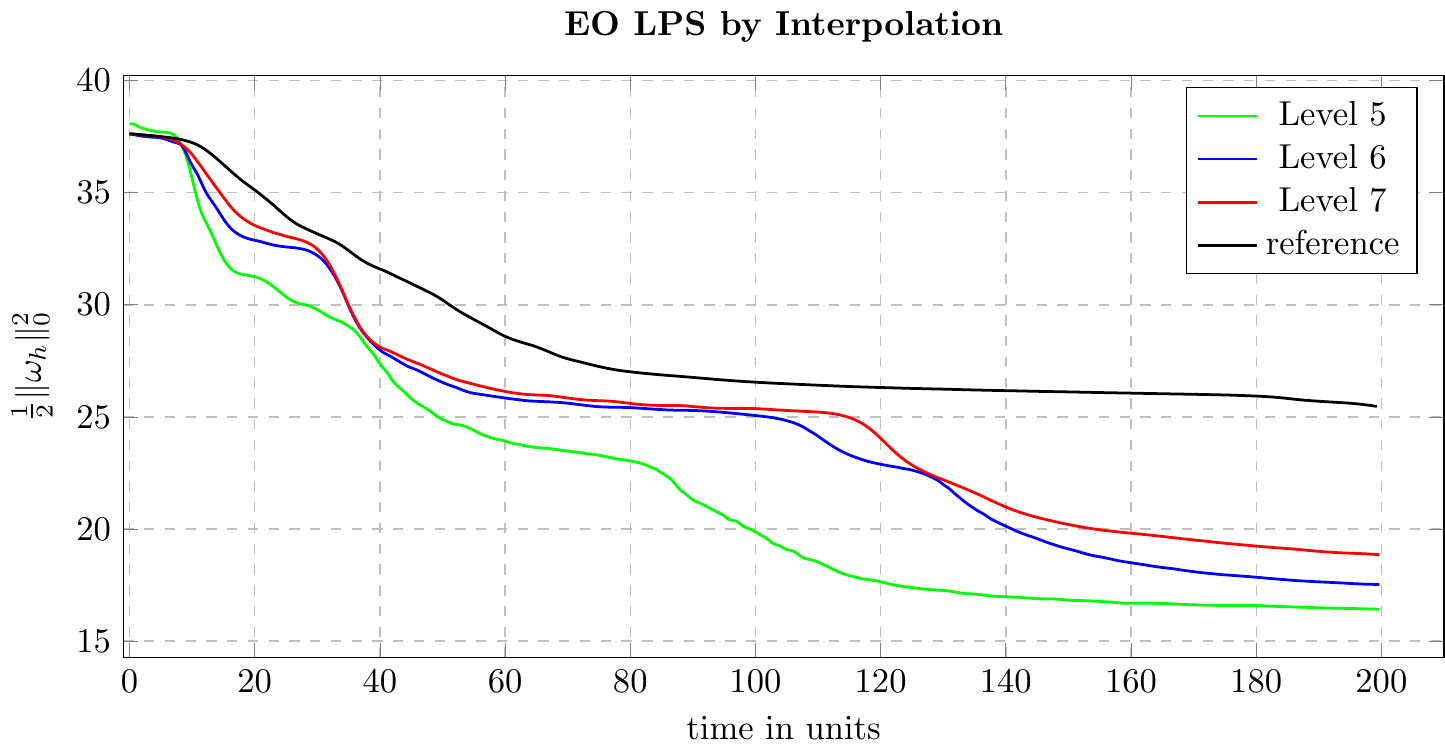}
	\caption{Temporal evolution of enstrophy with EO--FE: RB-VMS (top left), SUPG (top right), 
	one-level LPS (bottom left), and LPS by interpolation (bottom right), on different mesh 
	refinement levels, $\Delta t=0.0125$. \label{fig:enst_eo_0.0125}}
\end{figure}

\begin{figure}[t!]
	\centering
	\includegraphics[scale=0.5]{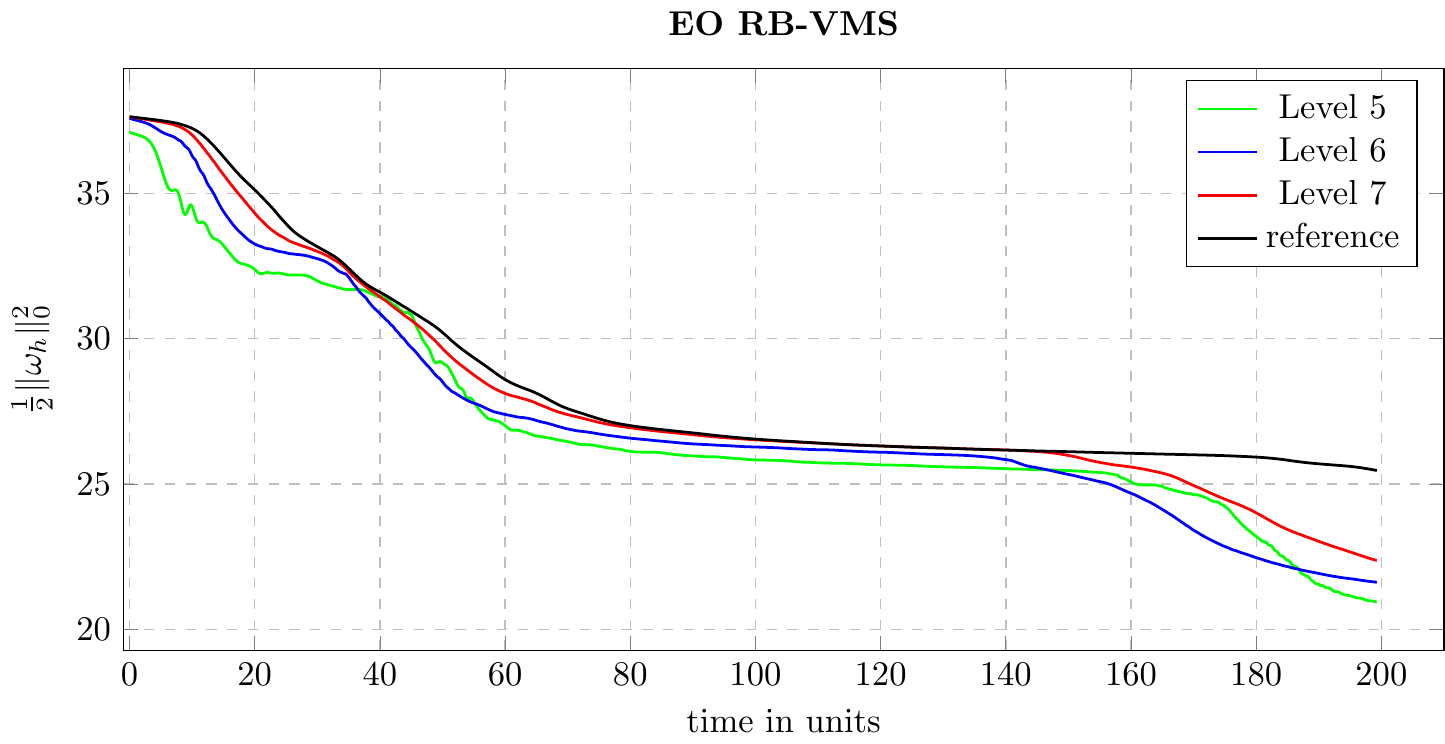}
	\includegraphics[scale=0.5]{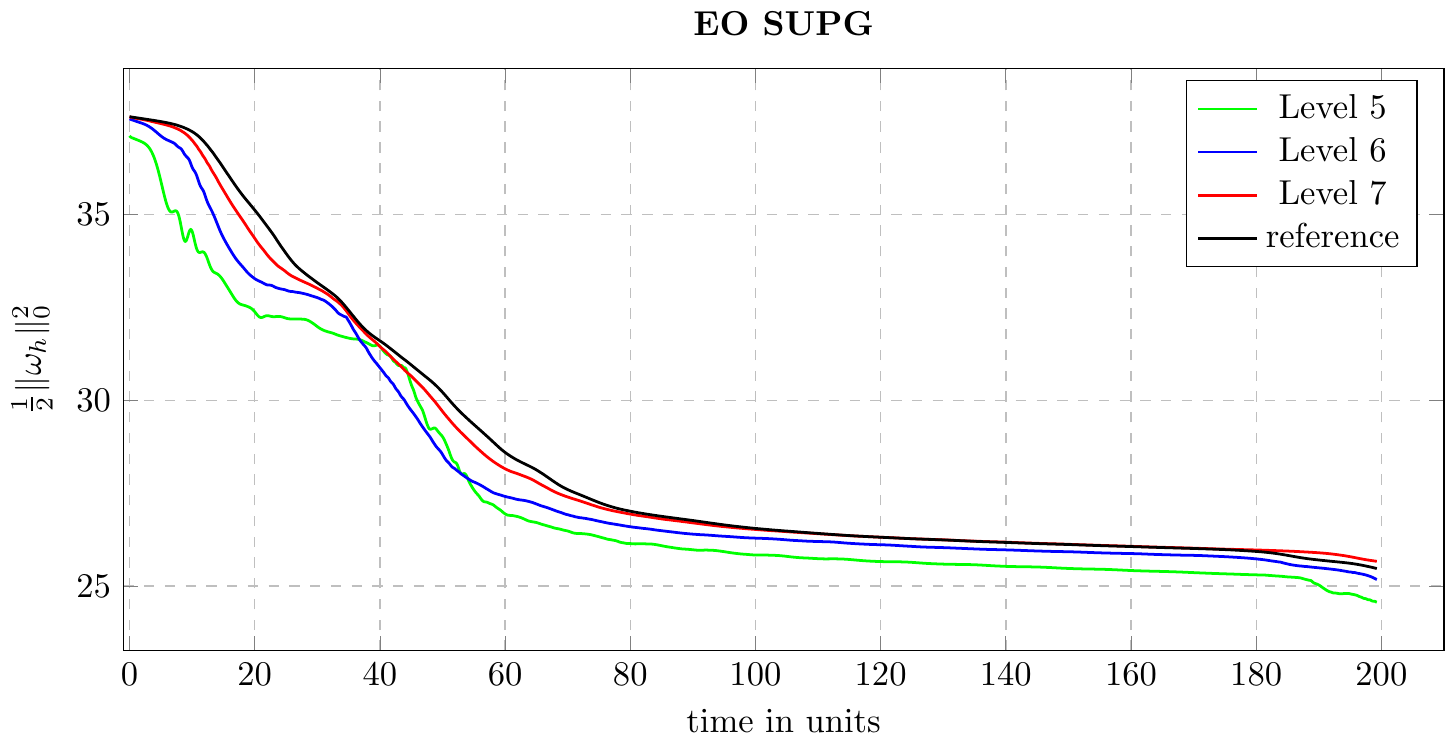}\\
	\includegraphics[scale=0.5]{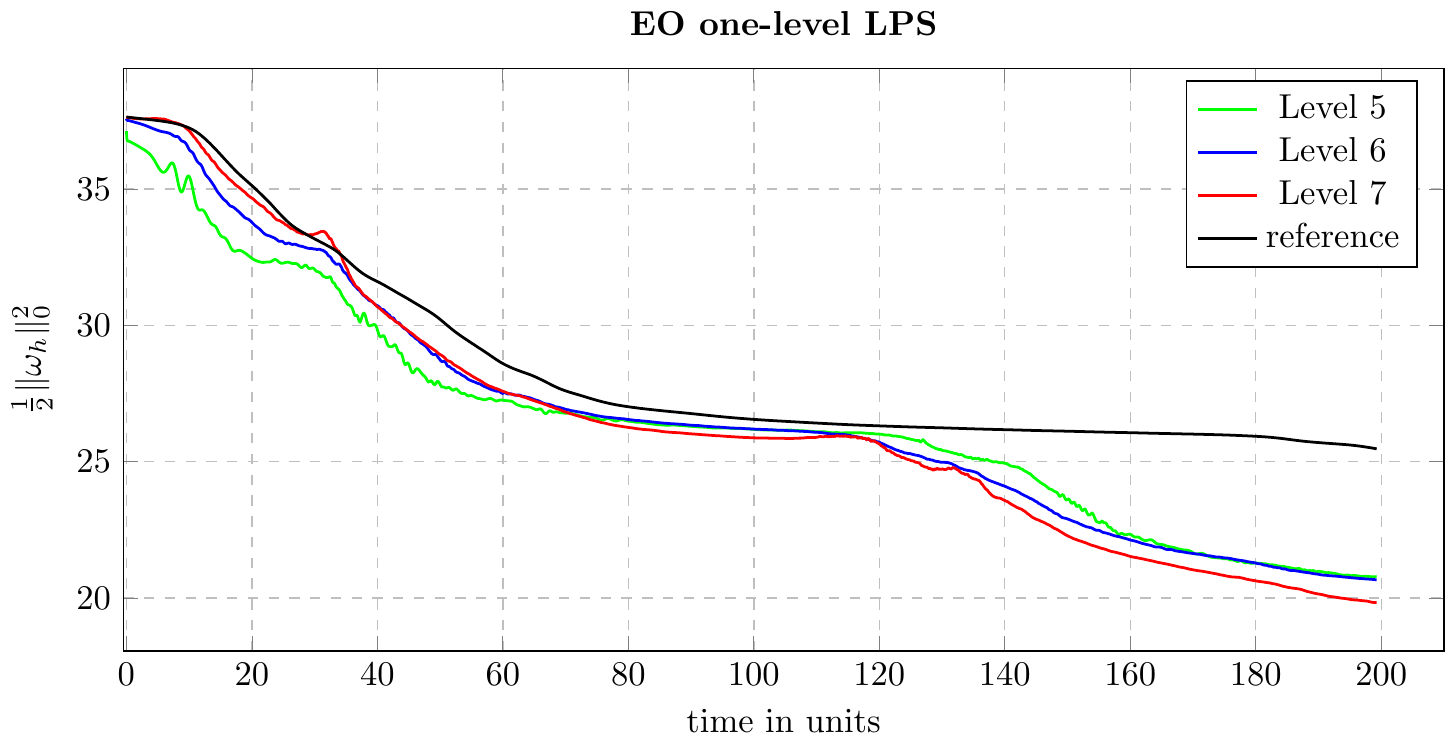}
	\includegraphics[scale=0.5]{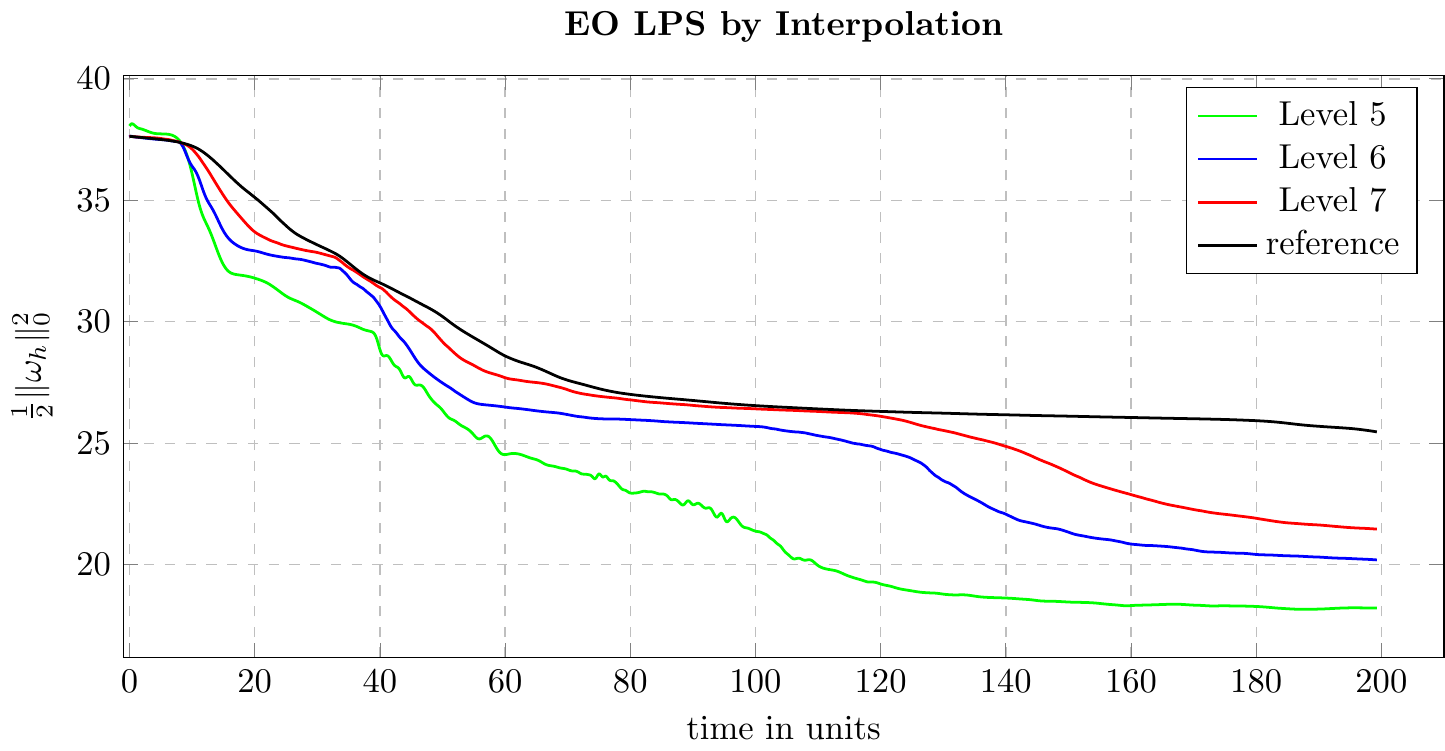}
	\caption{Temporal evolution of enstrophy with EO--FE: RB-VMS (top left), SUPG (top right), 
	one-level LPS (bottom left), and LPS by interpolation (bottom right), on different mesh 
	refinement levels, $\Delta t=0.003125$. \label{fig:enst_eo_0.003125}}
\end{figure}

\begin{figure}[t!]
	\centering
		\includegraphics[scale=0.5]{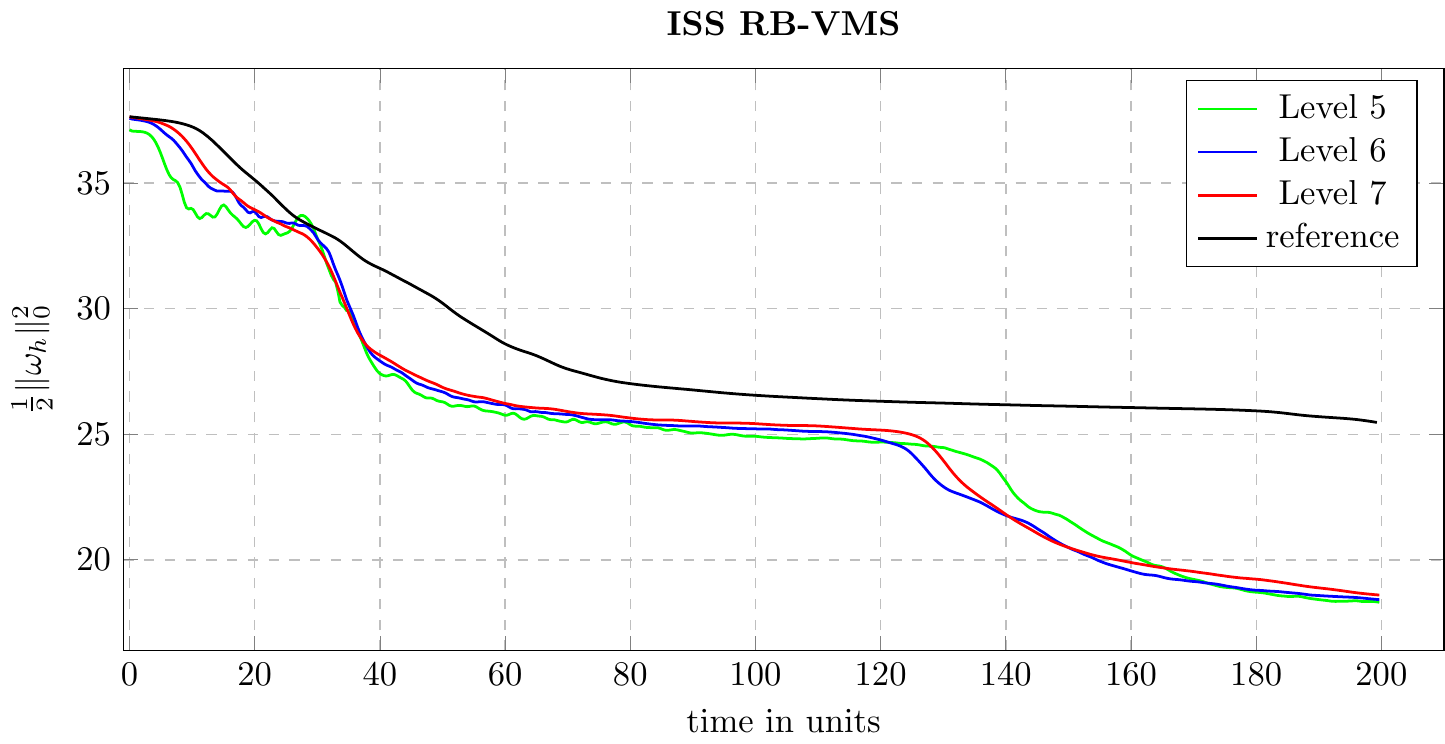}
		\includegraphics[scale=0.5]{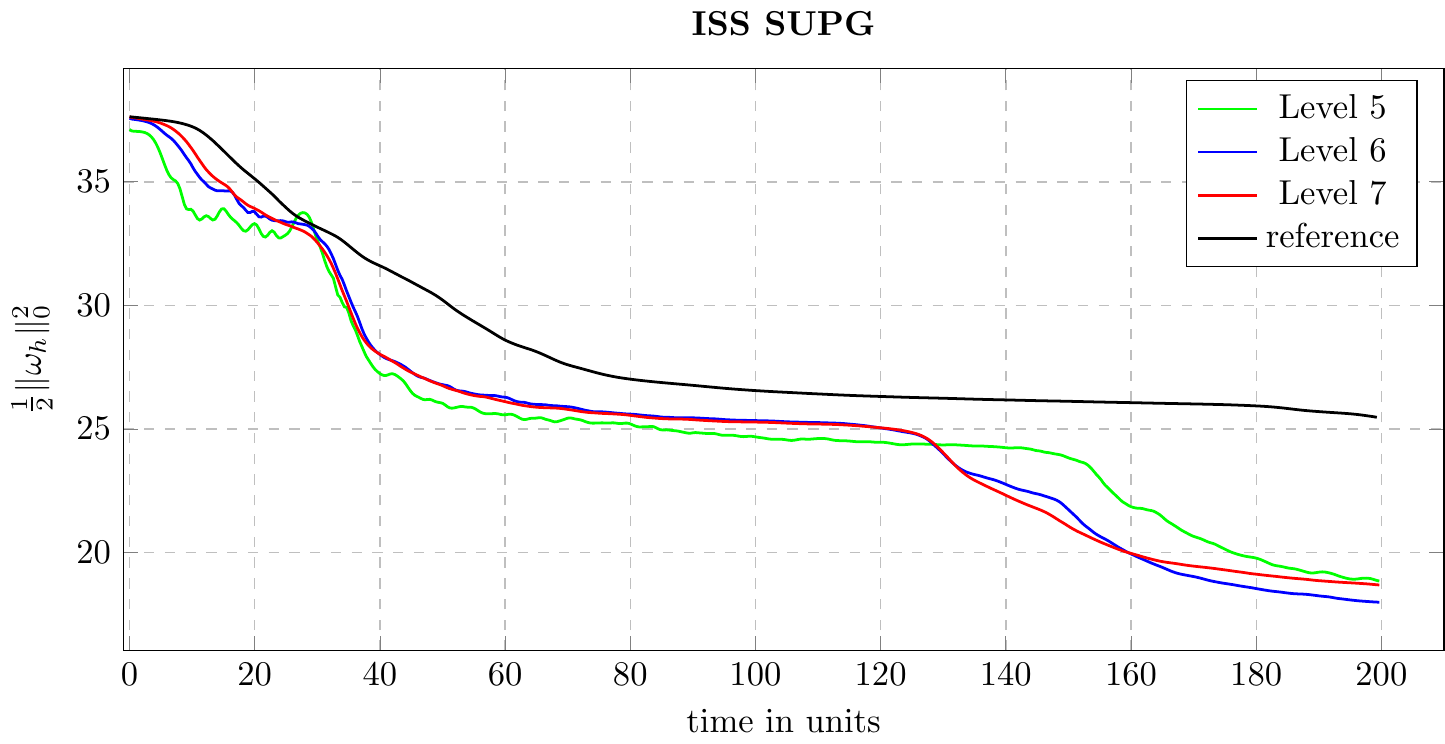}\\
		\includegraphics[scale=0.5]{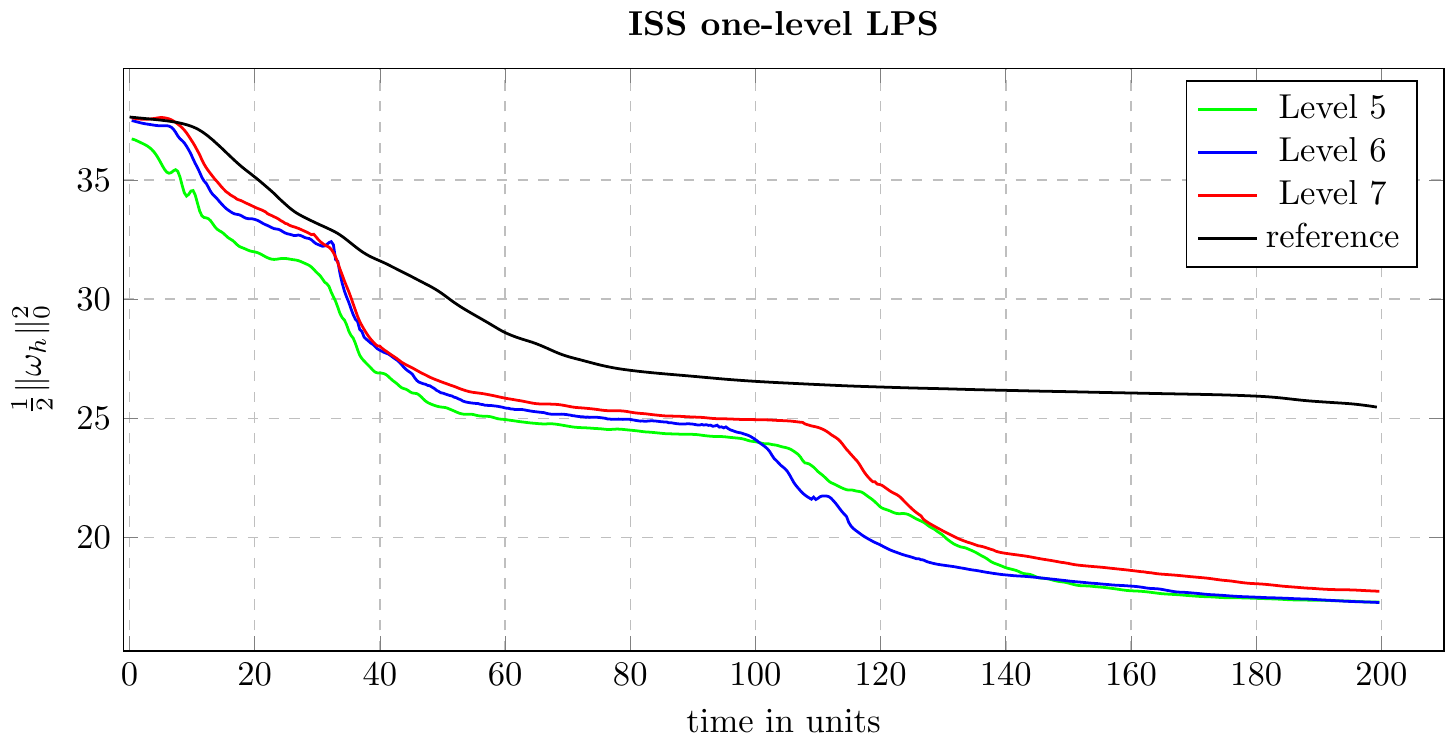}
		\includegraphics[scale=0.5]{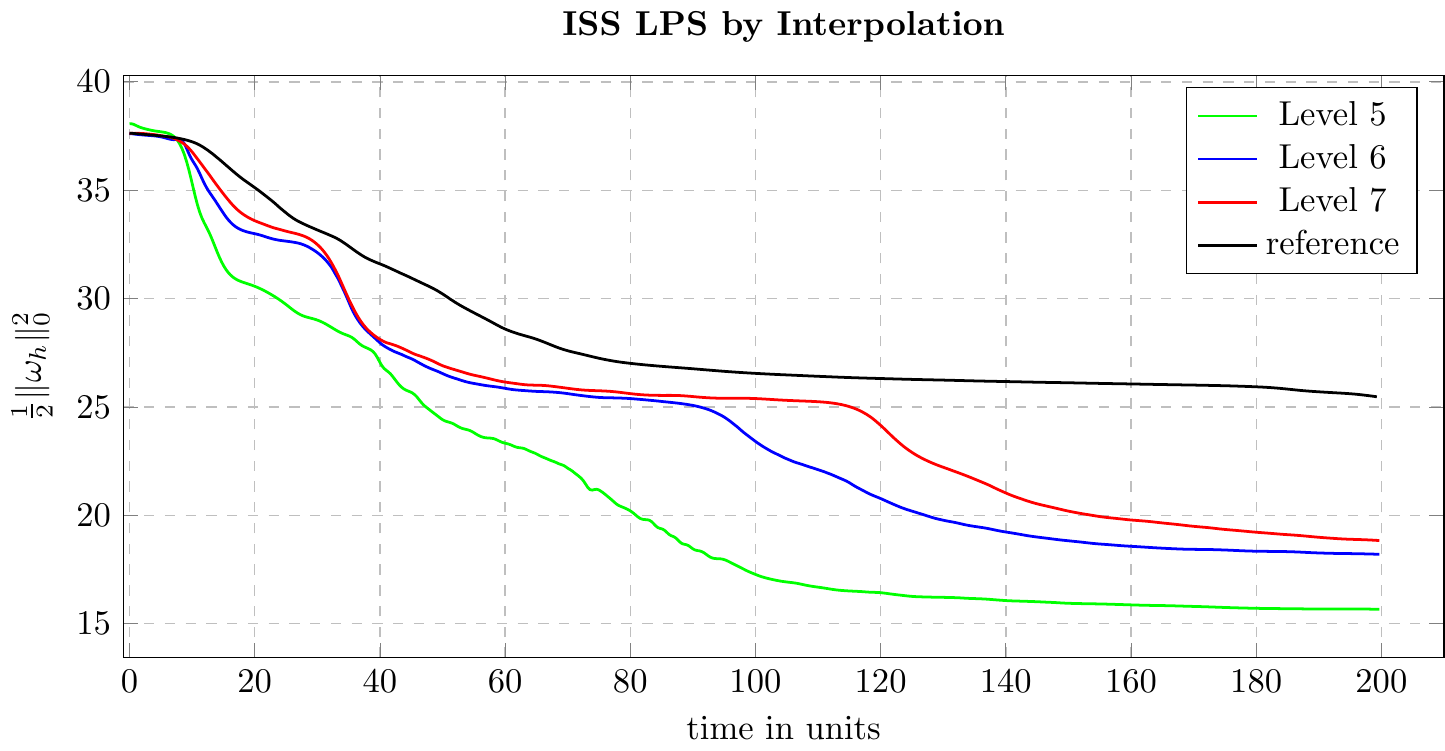}
	\caption{Temporal evolution of enstrophy with ISS--FE: RB-VMS (top left), SUPG (top right), 
	one-level LPS (bottom left), and LPS by interpolation (bottom right), on different mesh 
	refinement levels, $\Delta t=0.0125$. \label{fig:enst_is_0.0125}}
\end{figure}

\begin{figure}[t!]
	\centering
			\includegraphics[scale=0.5]{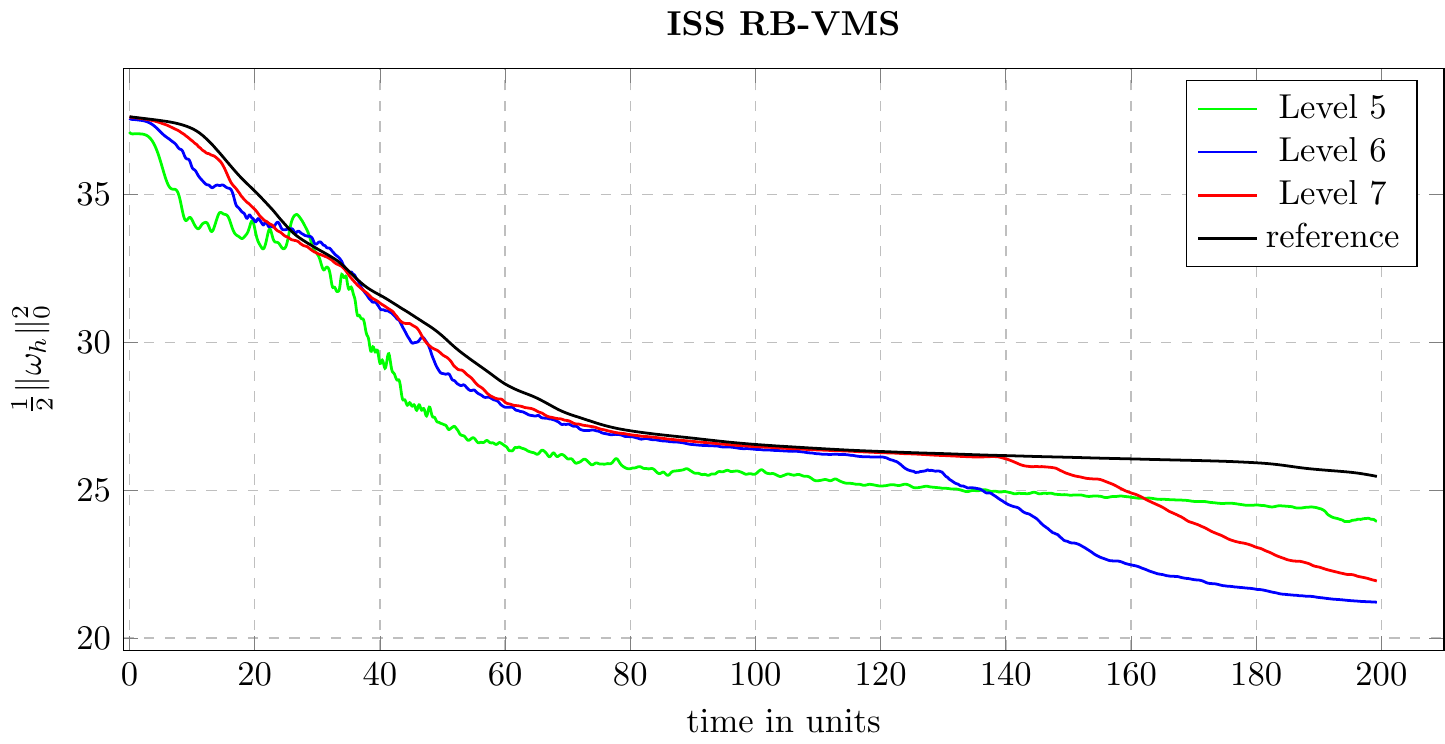}
			\includegraphics[scale=0.5]{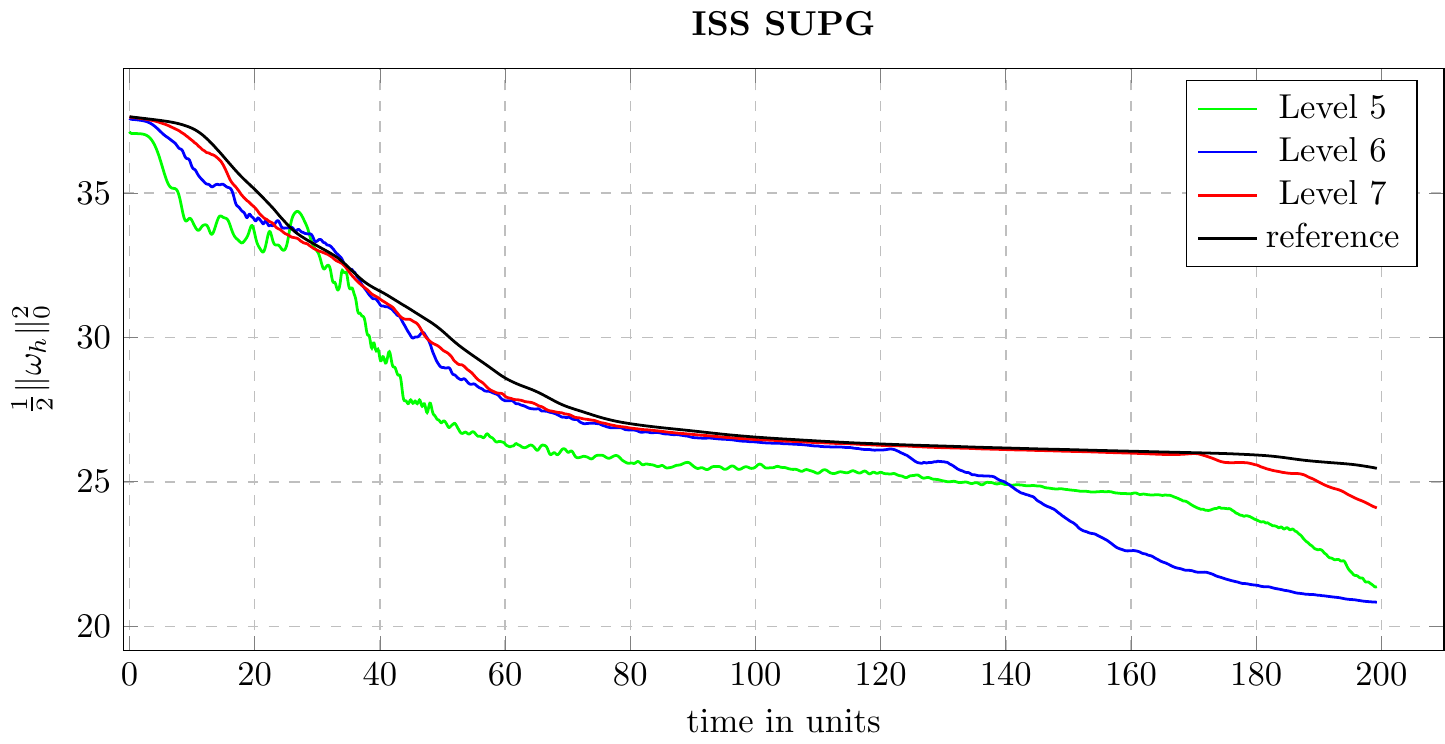}\\
			\includegraphics[scale=0.5]{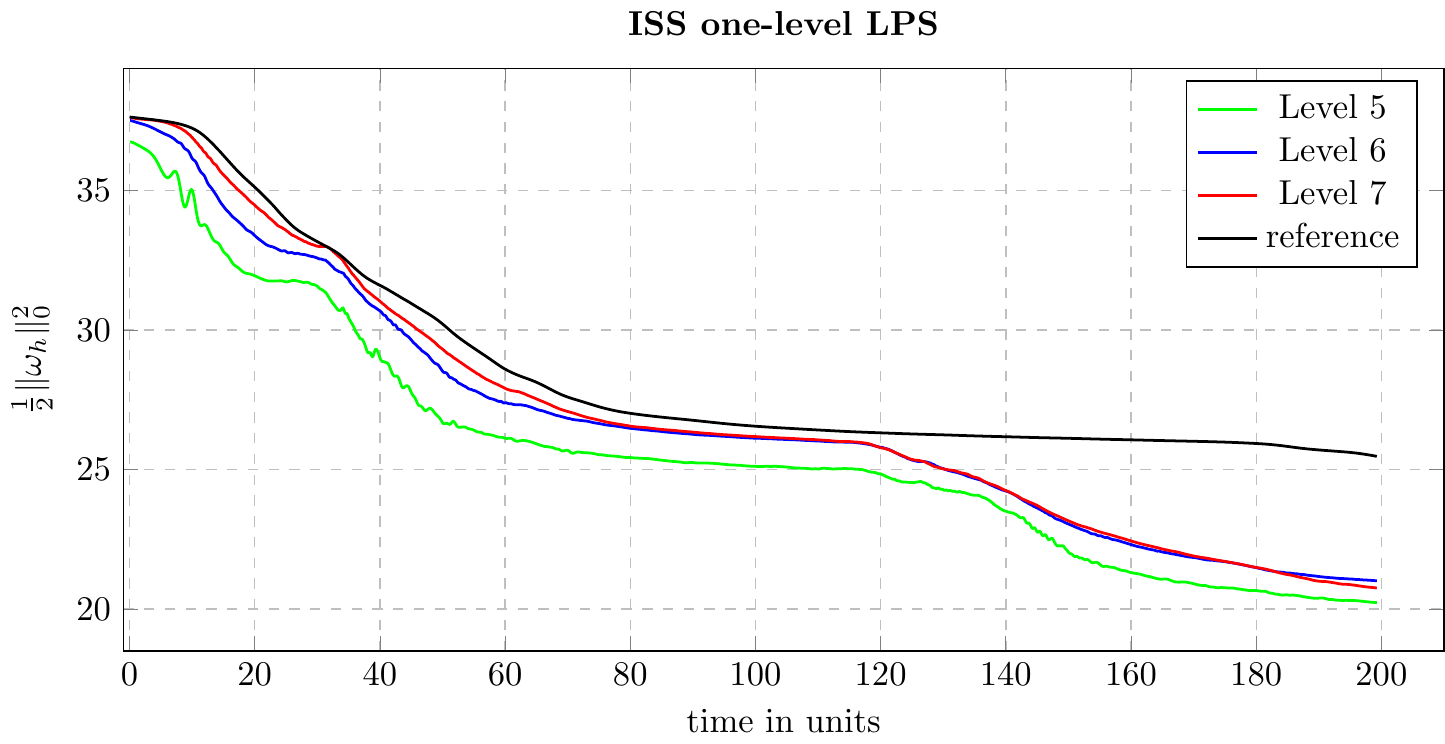}
			\includegraphics[scale=0.5]{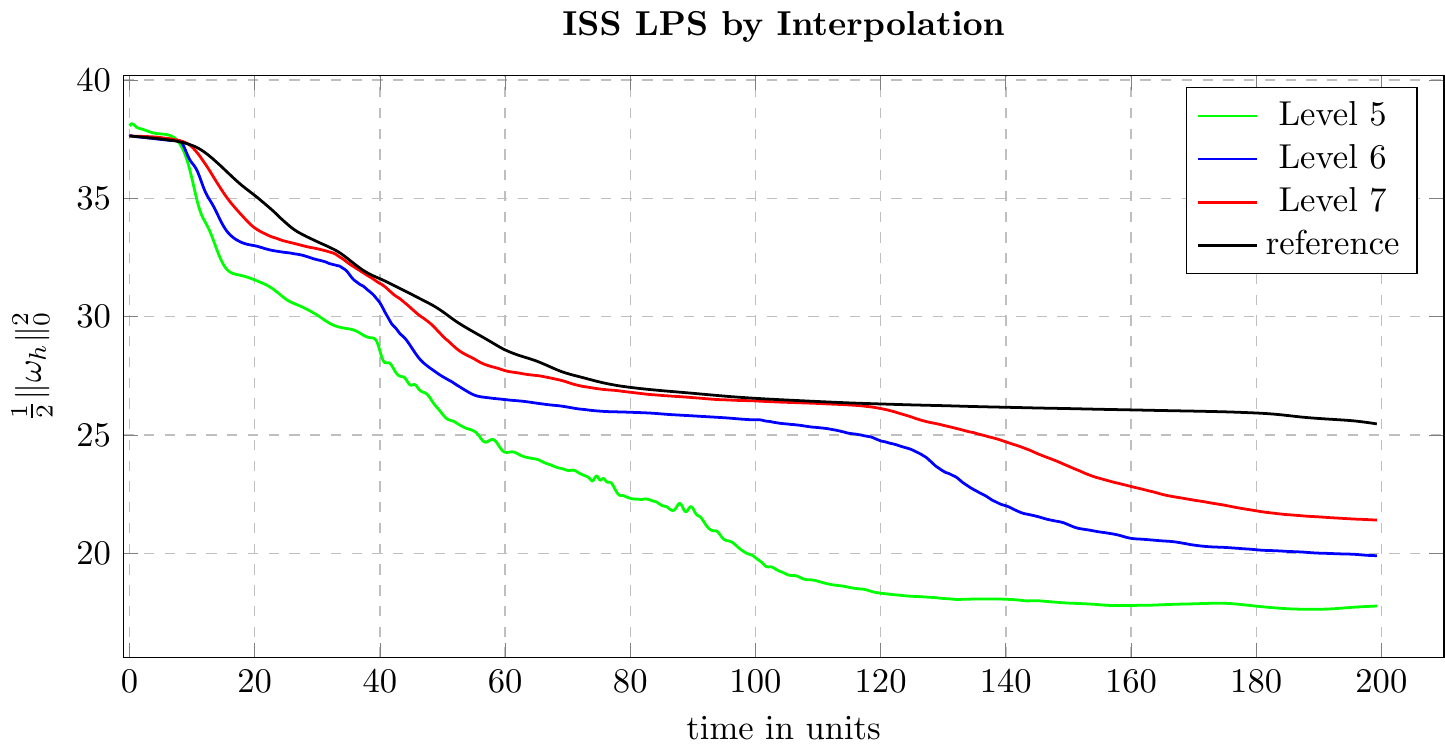}
	\caption{Temporal evolution of enstrophy with ISS--FE: RB-VMS (top left), SUPG (top right), 
	one-level LPS (bottom left), and LPS by interpolation (bottom right), on different mesh 
	refinement levels, $\Delta t=0.003125$. \label{fig:enst_is_0.003125}}
\end{figure}

The temporal evolution of the enstrophy is plotted in 
Figures~\ref{fig:enst_eo_0.0125}--\ref{fig:enst_is_0.003125} for 
all methods. Actually, as for the vorticity thickness, the stages 
of the enstrophy curves are directly connected to the pairing of vortices in the computation. 
Therefore, a more accurate method with higher resolution leads to a later and slower decrease in 
the enstrophy. Similar to the kinetic energy, the amount of the initial enstrophy is 
almost same for all simulations and results into a different final enstrophy for different 
refinements. A sudden decrease in the estrophy can be seen according to the pairing 
of eddies at different times. In agreement to the relative vorticity thickness and 
total kinetic energy, the best results in comparison with the reference solution and finest solution in~\cite{SL17}  
are obtained by EO SUPG method with small time step length \(\Delta t_2\), see Figure~\ref{fig:enst_eo_0.003125} (top right), for which a sort of mesh convergence is again reached. 

As before, note that there are no much noticeable differences between the use of EO and ISS FE also for the enstrophy results.

\subsection{Palinstrophy}\label{subsec:Palinst}
\begin{figure}[t!]
	\centering
\includegraphics[scale=0.5]{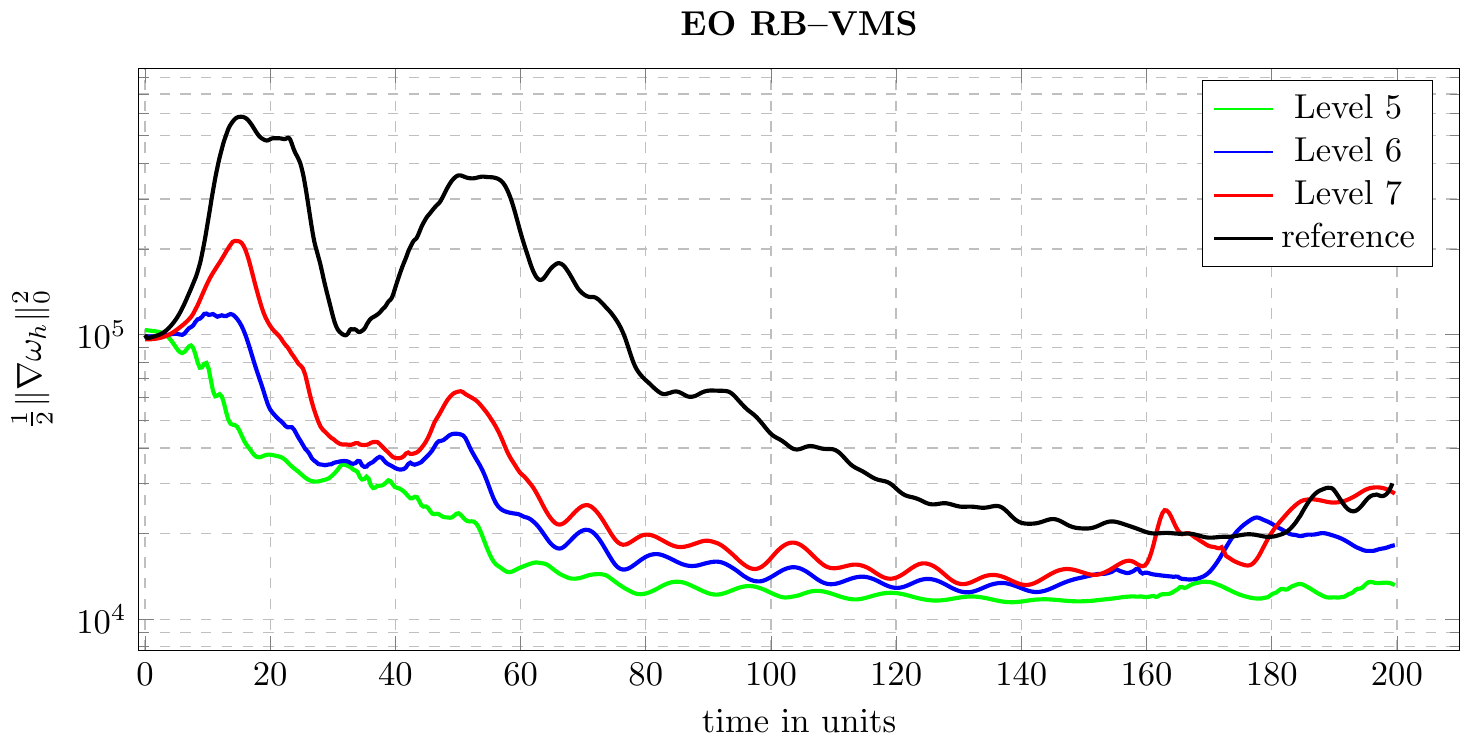}
\includegraphics[scale=0.5]{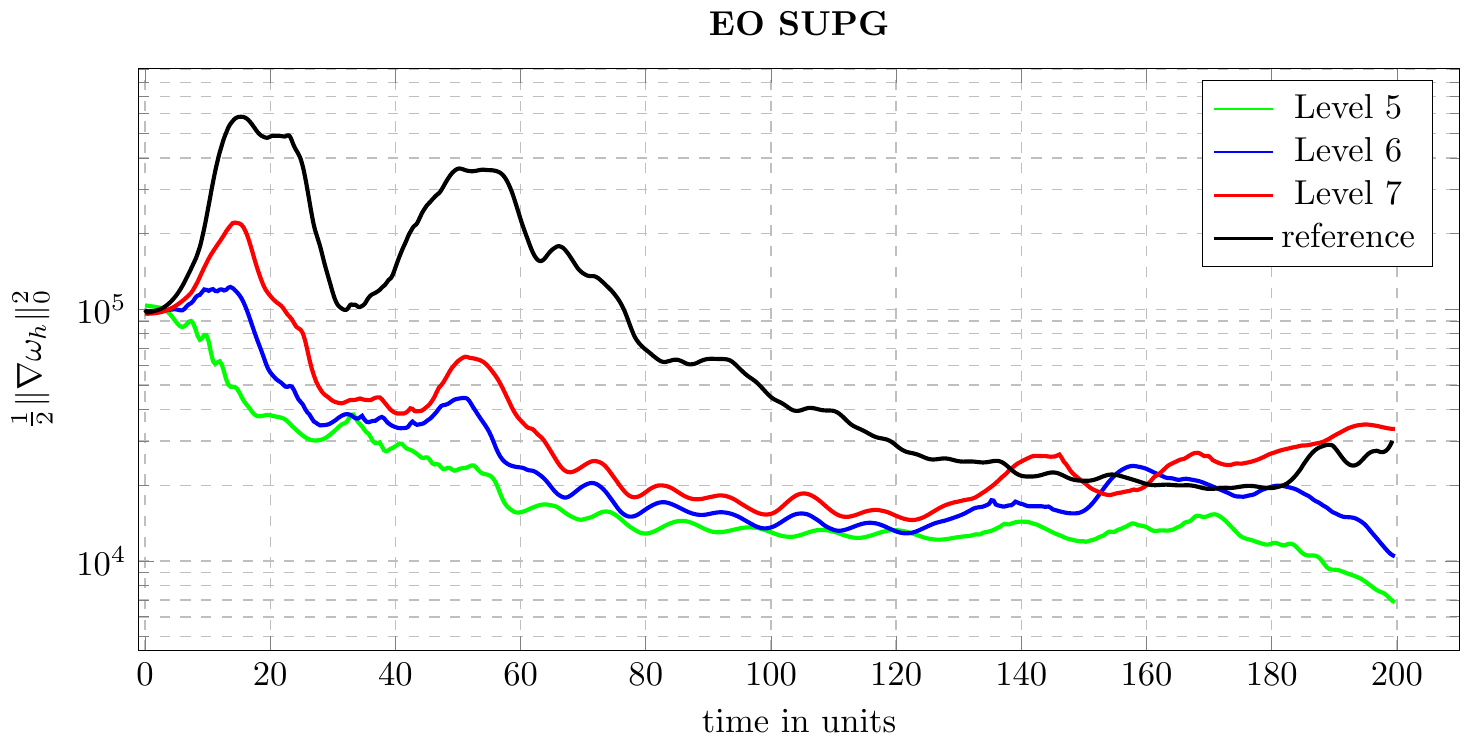}\\
\includegraphics[scale=0.5]{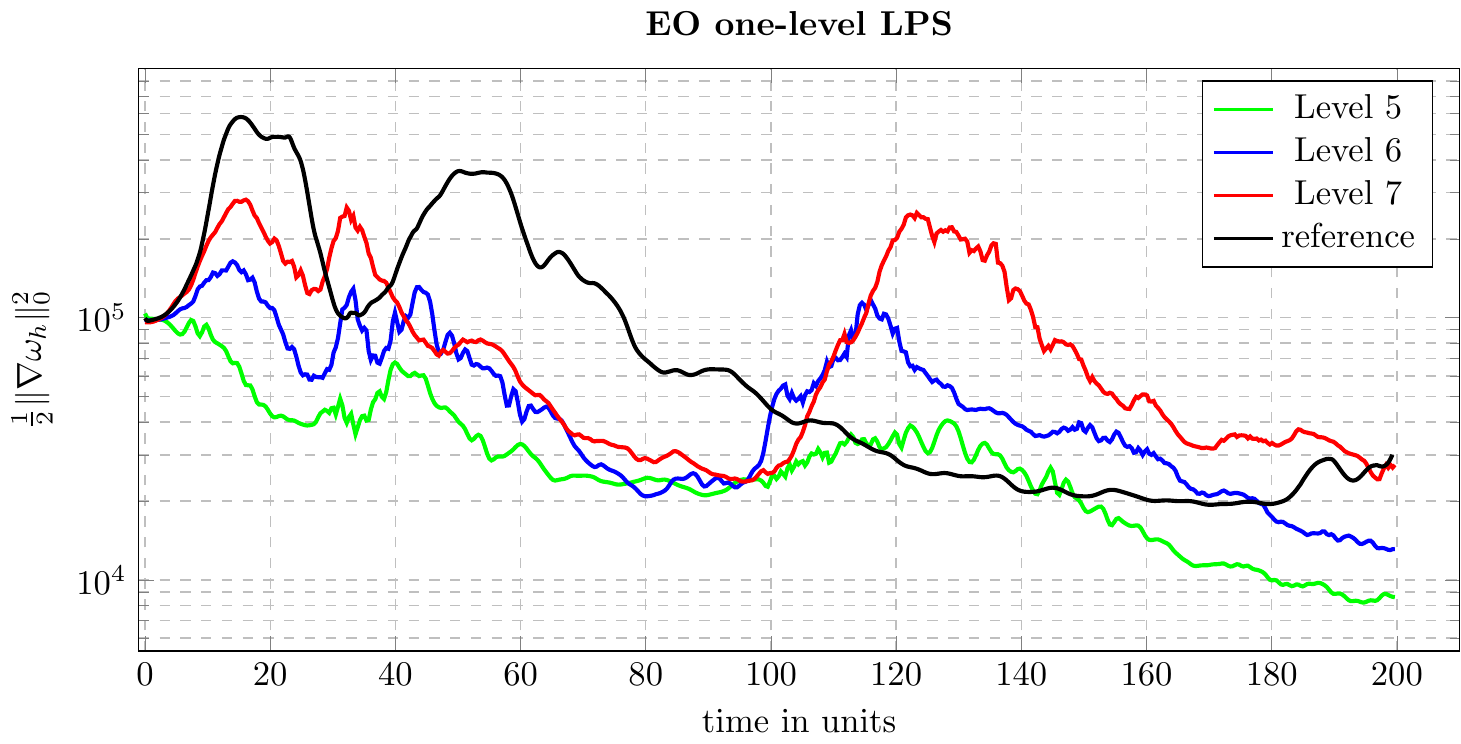}
\includegraphics[scale=0.5]{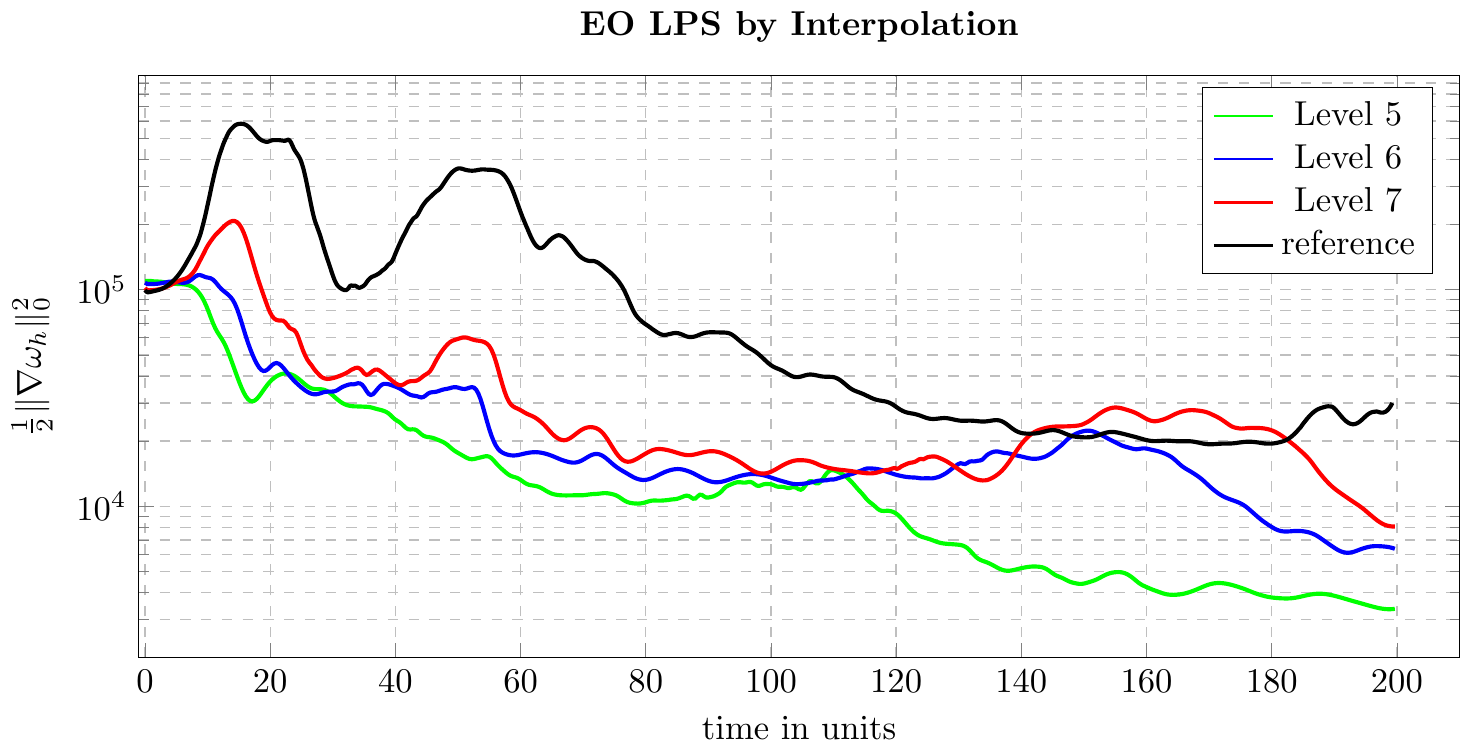}
	\caption{Temporal evolution of palinstrophy with EO--FE: RB-VMS (top left), SUPG (top right), 
	one-level LPS (bottom left), and LPS by interpolation (bottom right), on different mesh refinement 
	levels, $\Delta t=0.0125$. \label{fig:palinst_eo_0.0125}}
\end{figure}

\begin{figure}[t!]
	\centering
	\includegraphics[scale=0.5]{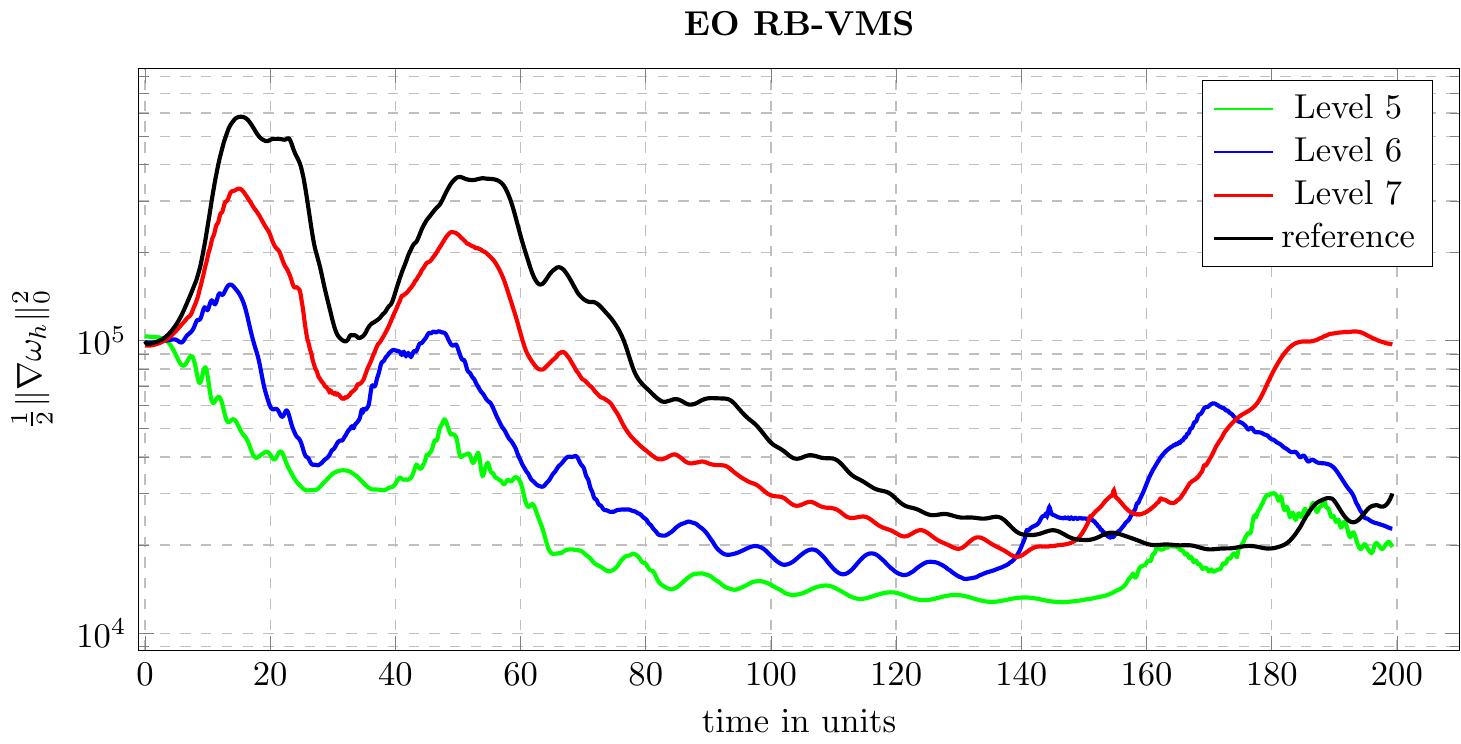}
	\includegraphics[scale=0.5]{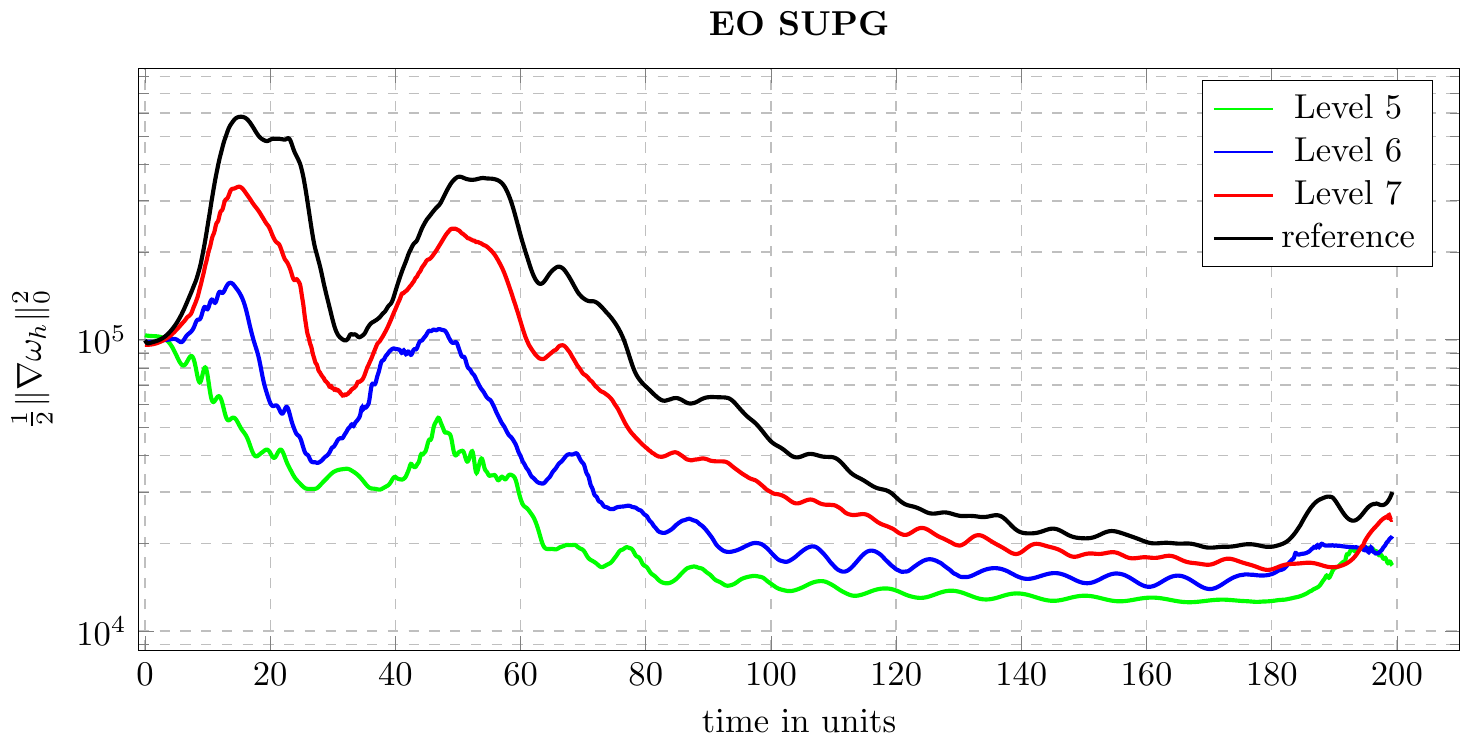}\\
	\includegraphics[scale=0.5]{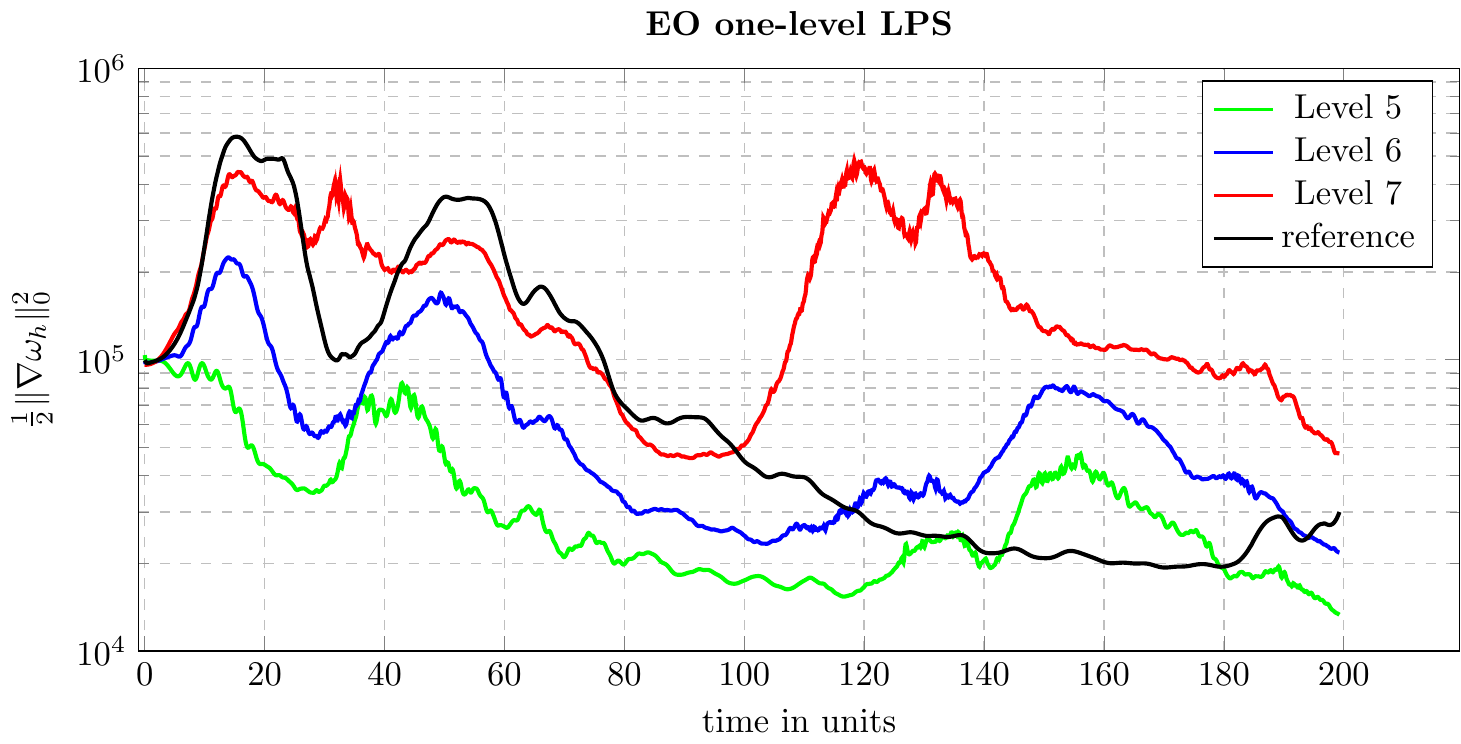}
	\includegraphics[scale=0.5]{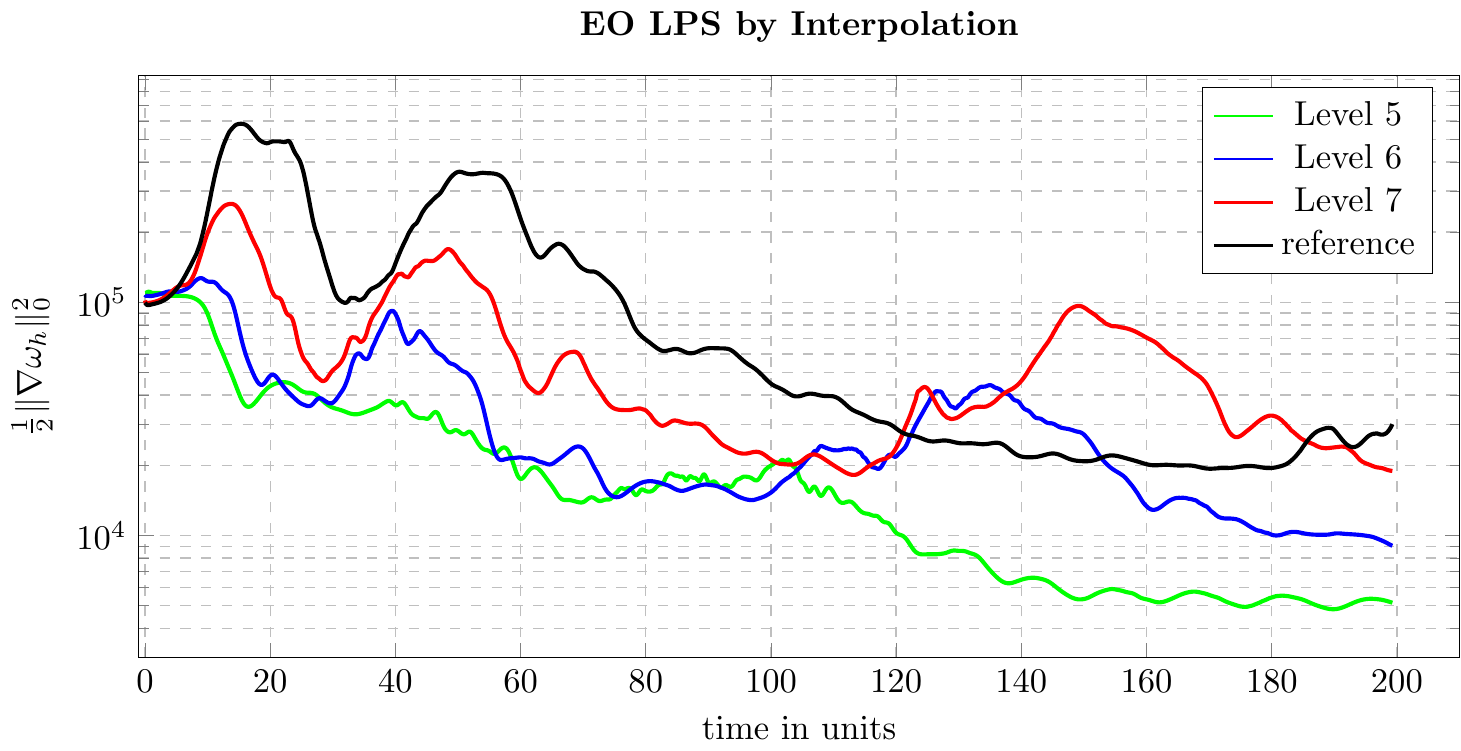}
	\caption{Temporal evolution of palinstrophy with EO--FE: RB-VMS (top left), SUPG (top right), 
	one-level LPS (bottom left), and LPS by interpolation (bottom right), on different mesh refinement 
	levels, $\Delta t=0.003125$. \label{fig:palinst_eo_0.003125}}
\end{figure}

\begin{figure}[t!]
	\centering
	\includegraphics[scale=0.5]{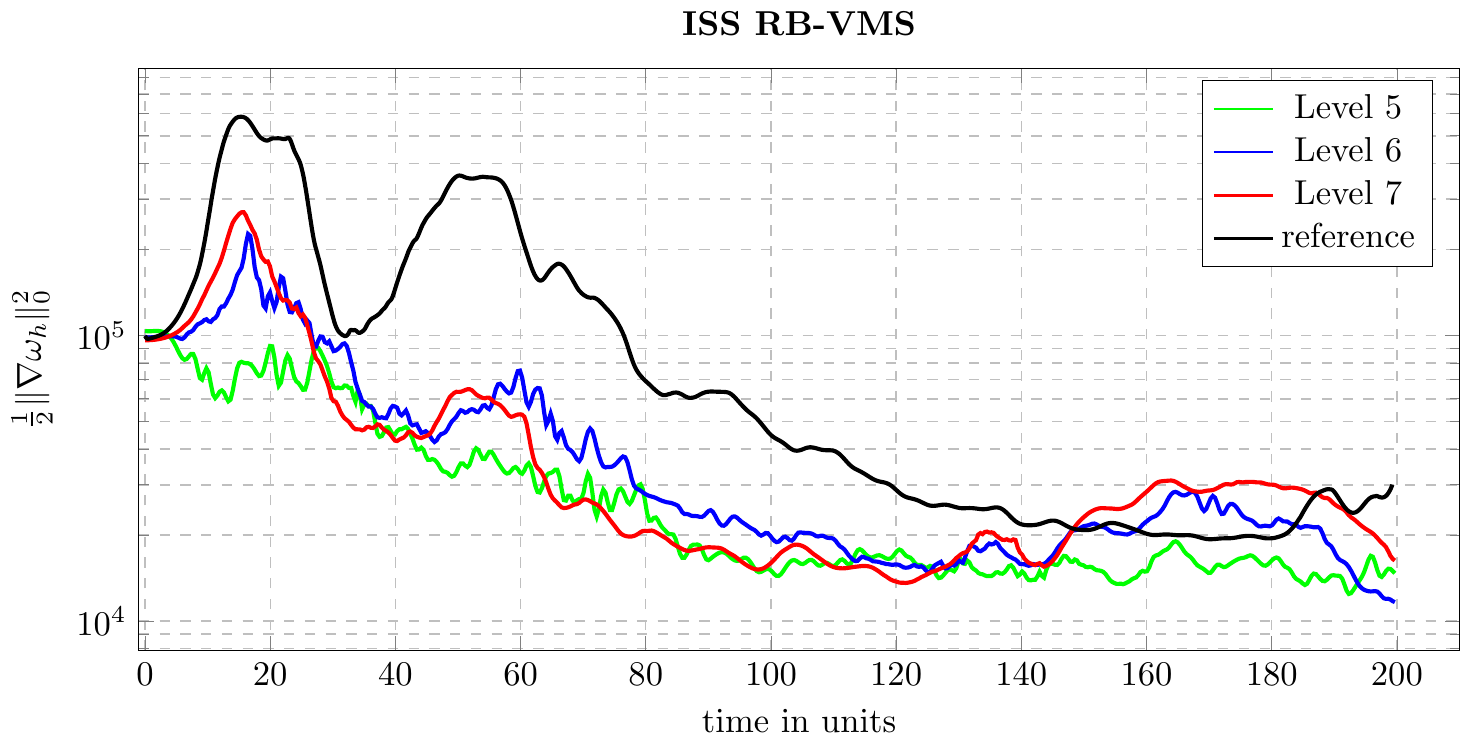}
	\includegraphics[scale=0.5]{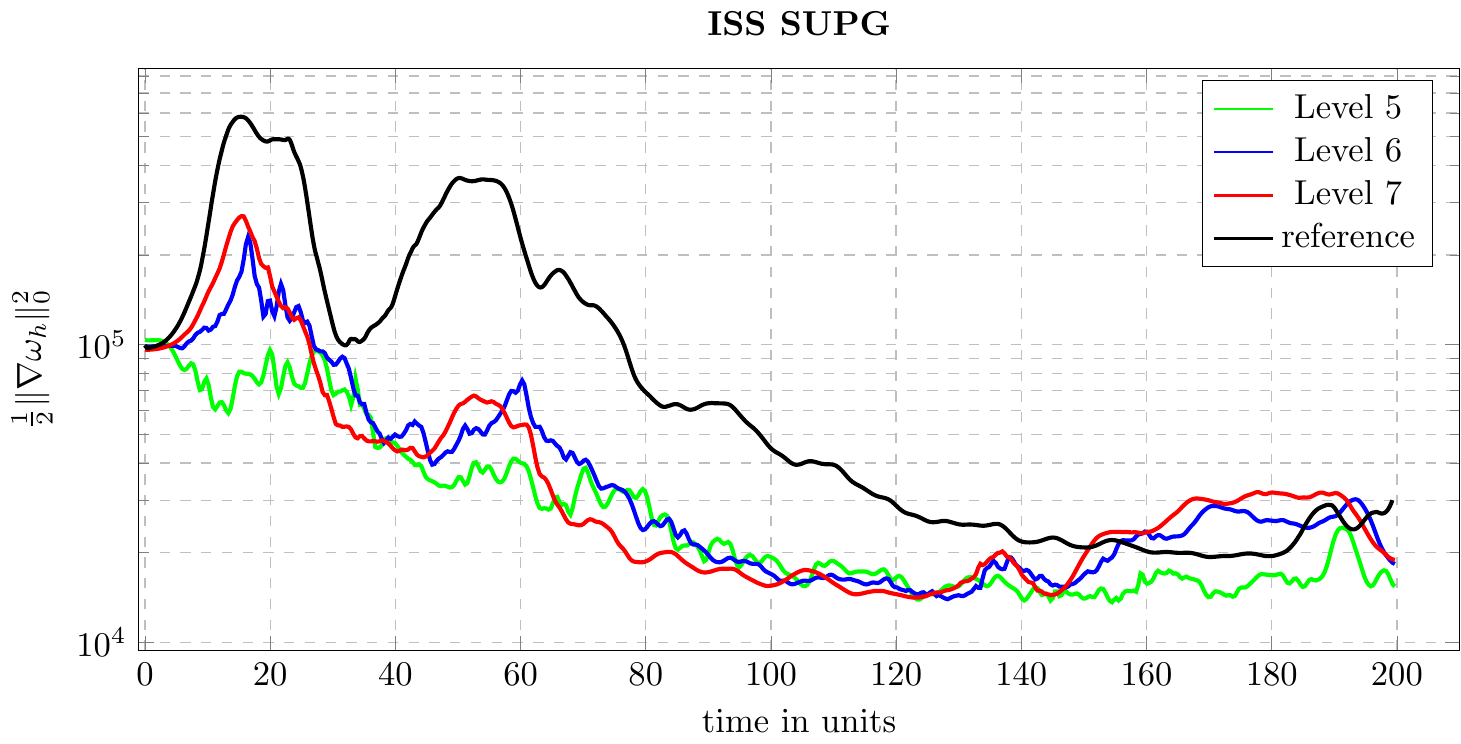}\\
	\includegraphics[scale=0.5]{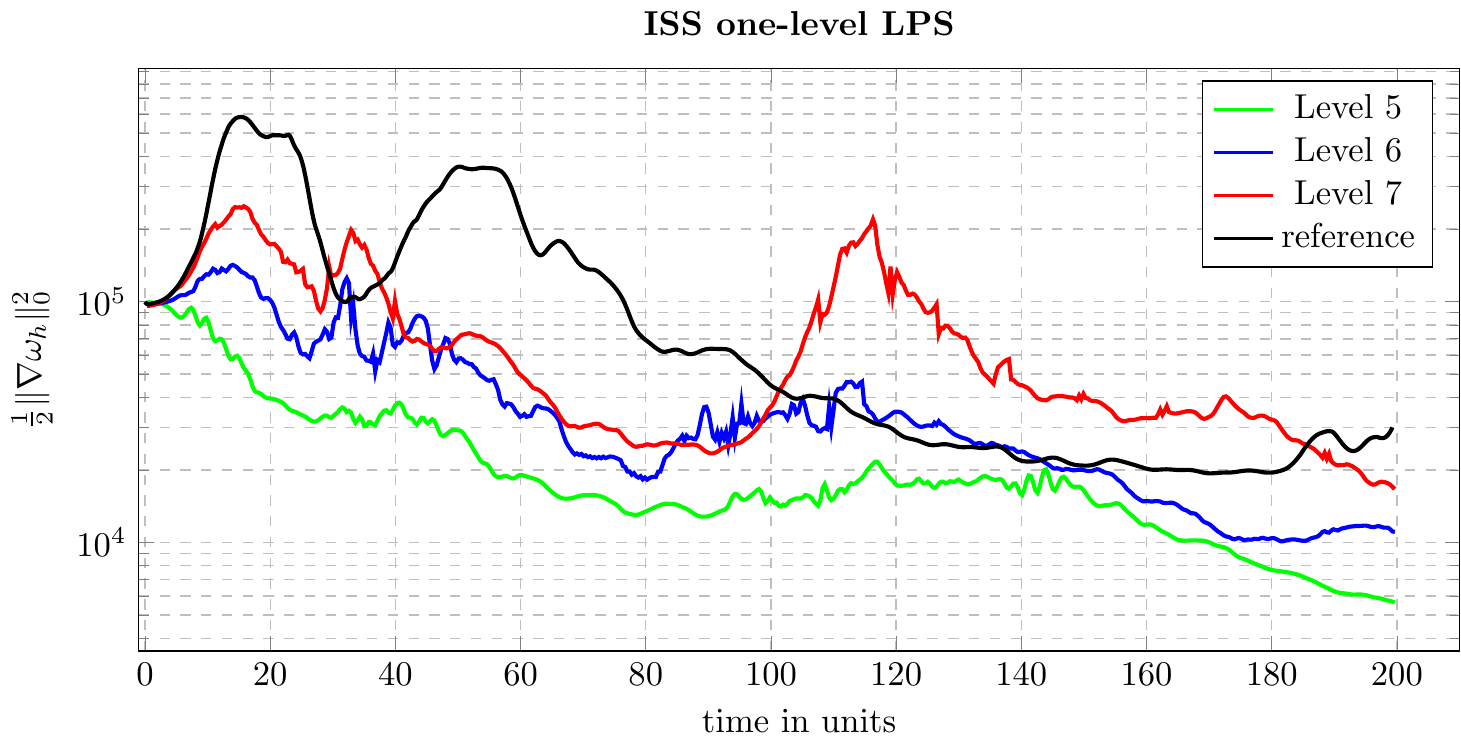}
	\includegraphics[scale=0.5]{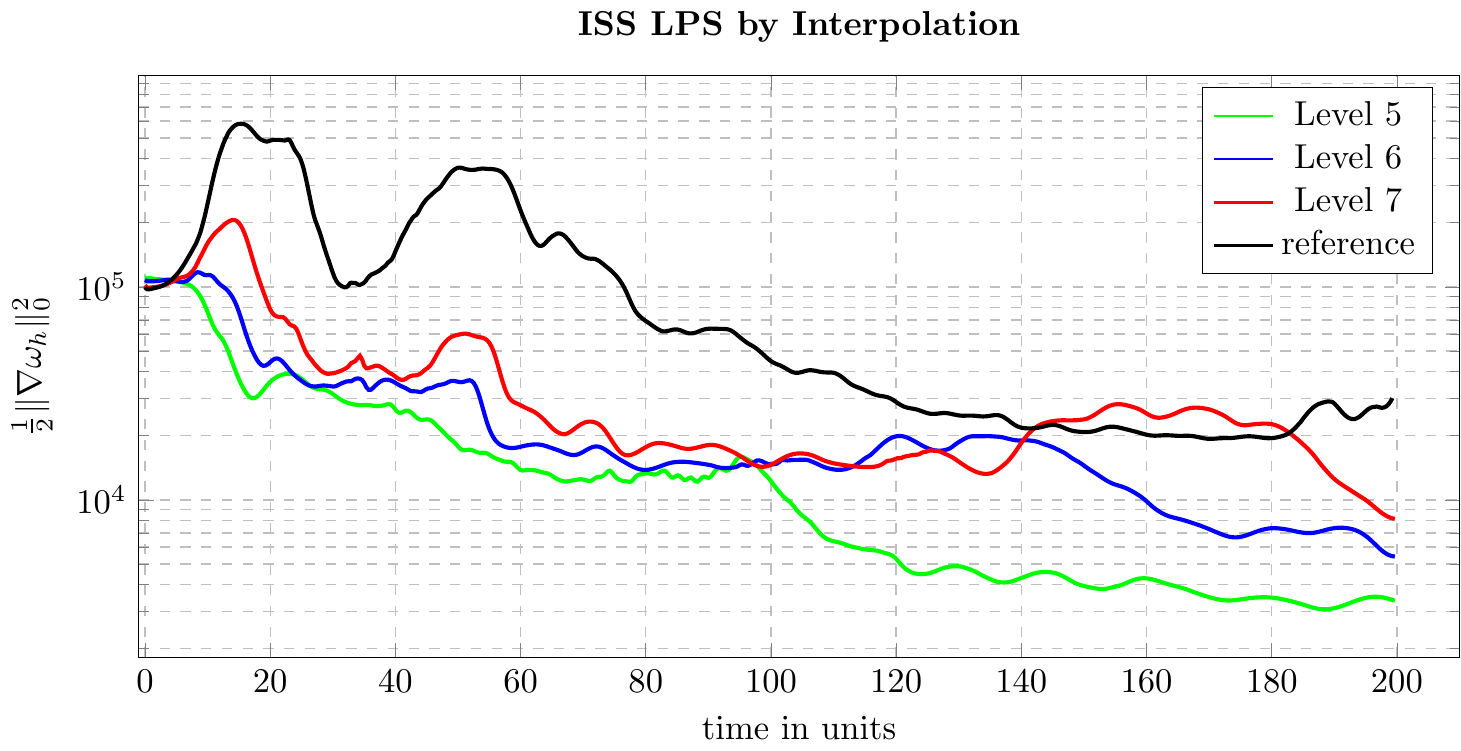}
	\caption{Temporal evolution of palinstrophy with ISS--FE: RB-VMS (top left), SUPG (top right), 
	one-level LPS (bottom left), and LPS by interpolation (bottom right), on different mesh refinement 
	levels, $\Delta t=0.0125$. \label{fig:palinst_is_0.0125}}
\end{figure}

\begin{figure}[t!]
	\centering
	\includegraphics[scale=0.5]{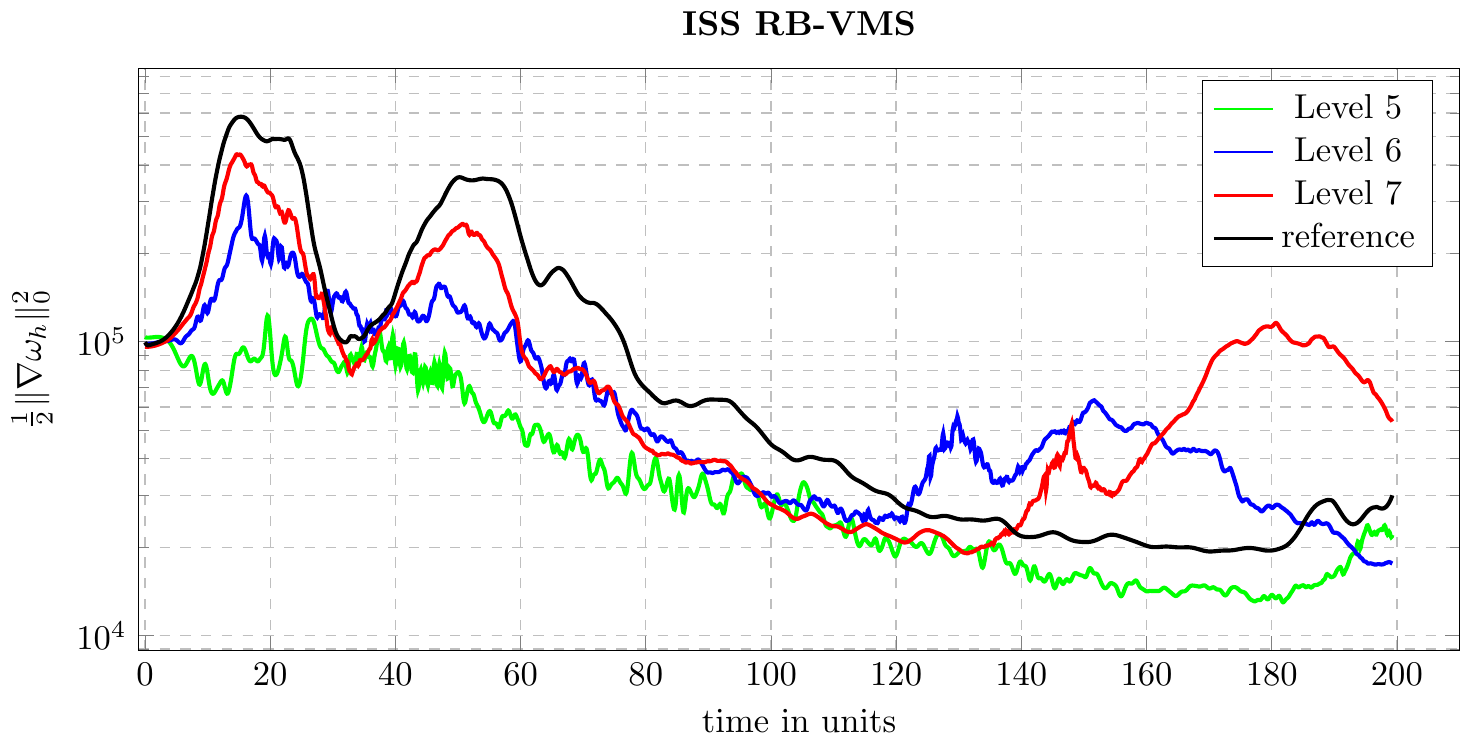}
	\includegraphics[scale=0.5]{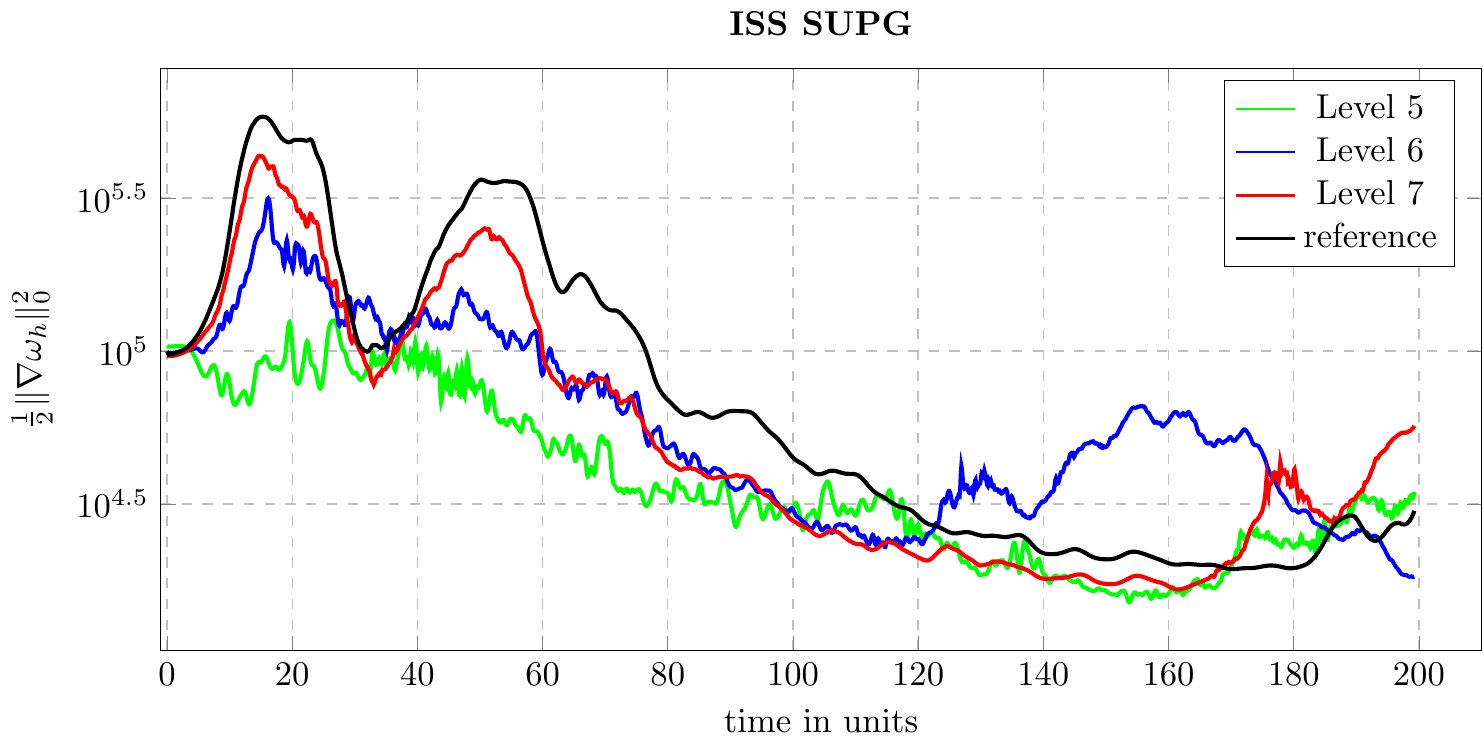}\\
	\includegraphics[scale=0.5]{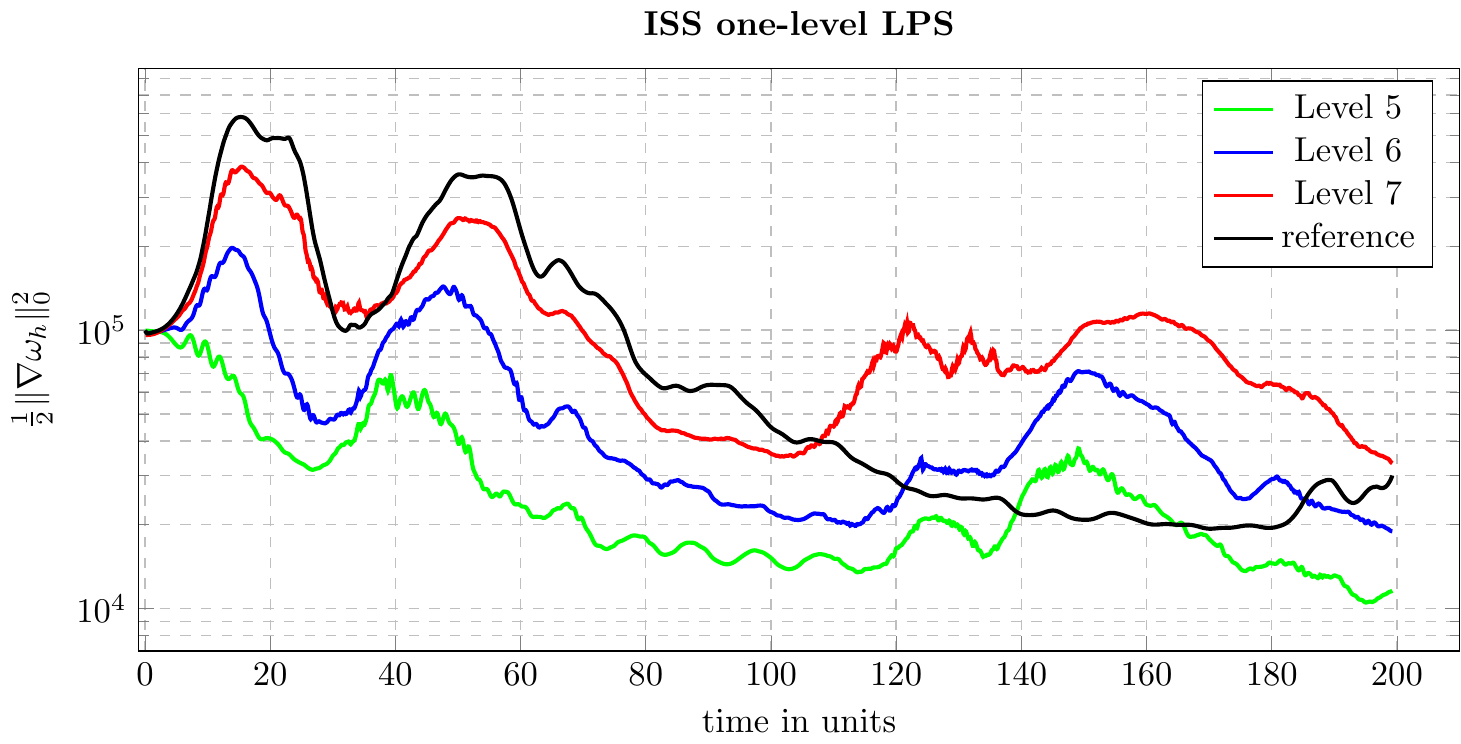}
	\includegraphics[scale=0.5]{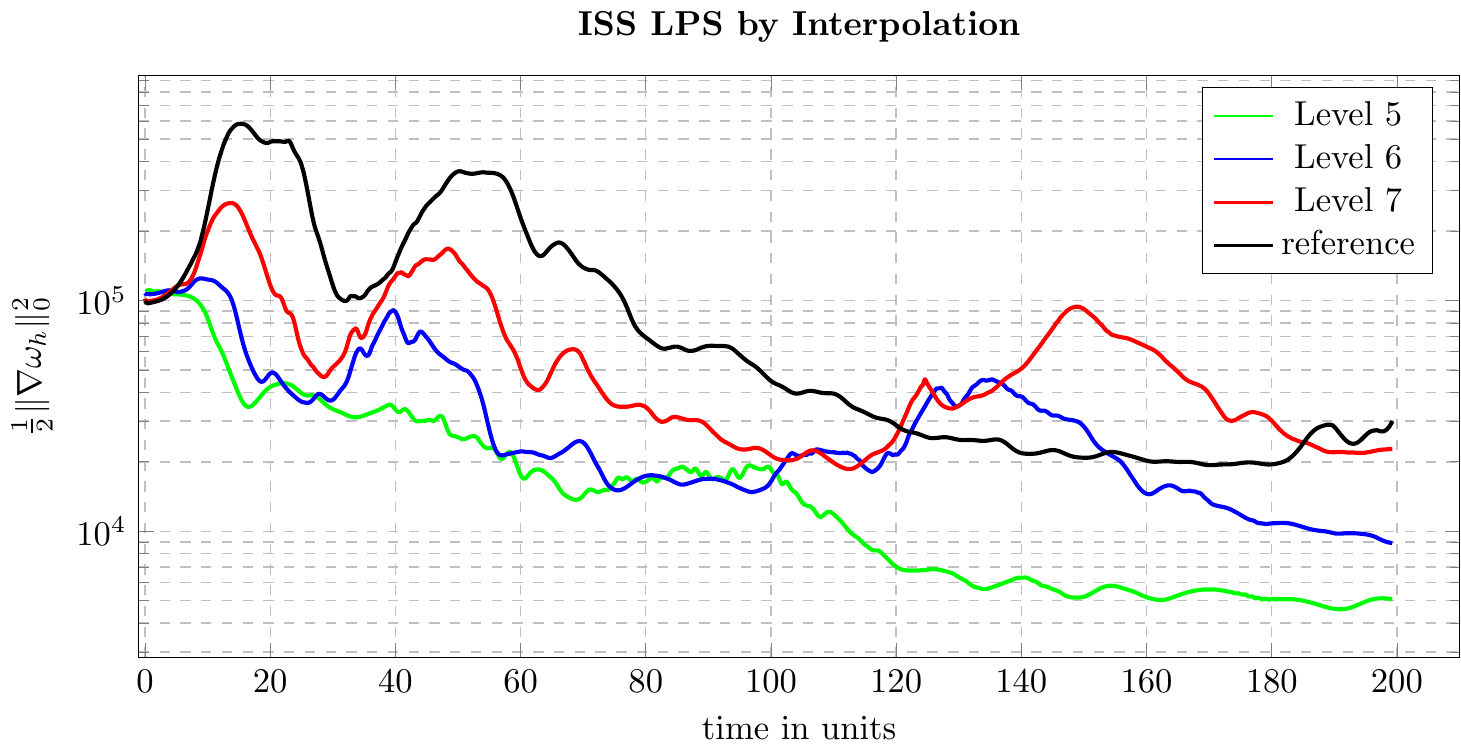}
	\caption{Temporal evolution of palinstrophy with ISS--FE: RB-VMS (top left), SUPG (top right), 
	one-level LPS (bottom left), and LPS by interpolation (bottom right), on different mesh refinement 
	levels, $\Delta t=0.003125$. \label{fig:palinst_is_0.003125}}
\end{figure}

The palinstrophy is one of the most sensitive quantity of interest, which makes it perfect to select the best results among comparisons. In 
Figures~\ref{fig:palinst_eo_0.0125}--\ref{fig:palinst_is_0.003125}, the temporal evaluation of the palinstrophy is presented
for all methods using EO and ISS FE on different refined meshes with time step lengths \(\Delta t_1\) and \(\Delta t_2\). 
In contrast to kinetic energy and enstrophy, palinstrophy can increase, and actually local maxima are almost attained once
merging processes of the vortices terminate. Reference solution indicates that the last merging process does not terminate before \(t=200\bar{t}\).
Over all, also in terms of magnitude, the closest results are obtained again by EO SUPG method with the small time step length \(\Delta t_2\), 
where one can see that the last merging process is still not ended at \(t=\bar{t}=200\). However, we have to stress that this quantity, both in
magnitude and time intervals for local maxima, is highly dependent on the studied method, refinement level, FE pairs, and time step lengths.

\section{Summary and outlook}\label{sec:conclusions}
In this paper, we compared two-scale VMS stabilized FE methods for the simulation of the 
incompressible NSE. These methods are widely used as one of the most 
promising and successful approaches that seek to simulate large-scale structures in turbulent 
flows. The space discretization for the studied methods using both ISS and 
EO FE is combined with a second-order semi-implicit time stepping scheme, based on BDF. 
Relatively coarse grids are chosen for the space discretization, using from large to small 
time step lengths. Several variants of two-scale VMS approaches, from fully residual-based to 
weakly consistent, have been applied to the simulation of 2D Kelvin--Helmholtz instabilities, 
triggered by a plane mixing layer at high Reynolds number \(Re=10^4\). 

Section~\ref{sec:numerics} presents in particular the detailed comparison of RB-VMS, SUPG, one-level 
variant of LPS and LPS by interpolation methods using both EO and ISS pair of FE on rather coarse grid 
levels and with different time steps, with the aim of studying their influence
on the accuracy of the numerical solutions. We discuss the numerical performances of all studied methods, 
by monitoring relevant quantities of interest, such as relative vorticity thickness, kinetic energy, 
enstrophy, and palinstrophy. From the computational point of view, note that this problem is very sensitive 
and results strongly depend on the used methods, grid refinement, and time step lengths. 

Through our numerical experiences, we have shown, for all methods, the need to consider a relatively small time step, both to 
prevent numerical stability issues proper of a less expensive semi-implicit time stepping scheme used here, leading for some 
methods to wrong results from the physical point of view (see, e.g., increase/oscillations in the kinetic energy), and to 
guarantee not excessive numerical dissipation. Altogether, based on the presented numerical studies, it turns out that the 
EO SUPG method with the small time step length outperforms all other studied variants. Closest results to this best performing 
method are attained by RB-VMS method, for which however the extra terms seem to not provide increased accuracy for the studied 
problem on relatively coarse grids, and thus there seems to be no reason to extend the simpler SUPG method by the higher order 
terms of the more complex RB-VMS method in this case. On the other side, LPS methods, which are not fully consistent, but of 
optimal order with respect to the FE interpolation, despite their appealing structure both in terms of practical implementations 
such as to perform the numerical analysis, seems to need higher space resolutions in order to achieve the same accuracy of fully 
residual-based VMS stabilized methods. 

As a future research direction, we plan to compare the selected best performing two-scale VMS stabilized methods towards several 
variants of three-scale VMS methods that use turbulent eddy viscosity (in a more or less sophisticated manner) to model the effect 
of subgrid scales, also on more complex problems presenting genuine 3D turbulent structure, like 3D turbulent channel flow.


\begin{thebibliography}{10}
\expandafter\ifx\csname url\endcsname\relax
  \def\url#1{\texttt{#1}}\fi
\expandafter\ifx\csname urlprefix\endcsname\relax\def\urlprefix{URL }\fi
\expandafter\ifx\csname href\endcsname\relax
  \def\href#1#2{#2} \def\path#1{#1}\fi

\bibitem{ACJR17}
N.~Ahmed, T.~Chac{\'o}n~Rebollo, V.~John, S.~Rubino, A review of variational
  multiscale methods for the simulation of turbulent incompressible flows,
  Archives of Computational Methods in Engineering 24~(1) (2017) 115--164.

\bibitem{JK10_cmame}
V.~John, A.~Kindl, Numerical studies of finite element variational multiscale
  methods for turbulent flow simulations, Comput. Methods Appl. Mech. Engrg.
  199~(13-16) (2010) 841--852.

\bibitem{BCCHRS07}
Y.~Bazilevs, V.~M. Calo, J.~A. Cottrell, T.~J.~R. Hughes, A.~Reali,
  G.~Scovazzi, Variational multiscale residual-based turbulence modeling for
  large eddy simulation of incompressible flows, Comput. Methods Appl. Mech.
  Engrg. 197~(1-4) (2007) 173--201.

\bibitem{BH82}
A.~N. Brooks, T.~J.~R. Hughes, Streamline upwind/{P}etrov-{G}alerkin
  formulations for convection dominated flows with particular emphasis on the
  incompressible {N}avier-{S}tokes equations, Comput. Methods Appl. Mech.
  Engrg. 32~(1-3) (1982) 199--259, fENOMECH '81, Part I (Stuttgart, 1981).

\bibitem{HB79}
T.~J.~R. Hughes, A.~Brooks, A multidimensional upwind scheme with no crosswind
  diffusion, in: Finite element methods for convection dominated flows
  ({P}apers, {W}inter {A}nn. {M}eeting {A}mer. {S}oc. {M}ech. {E}ngrs., {N}ew
  {Y}ork, 1979), Vol.~34 of AMD, Amer. Soc. Mech. Engrs. (ASME), New York,
  1979, pp. 19--35.

\bibitem{beckerbraack}
R.~Becker, M.~Braack, A finite element pressure gradient stabilization for the
  {S}tokes equations based on local projections, Calcolo 38~(4) (2001)
  173--199.

\bibitem{matskytob}
G.~Matthies, P.~Skrzypacz, L.~Tobiska, A unified convergence analysis for local
  projection stabilisations applied to the {O}seen problem, M2AN Math. Model.
  Numer. Anal. 41~(4) (2007) 713--742.

\bibitem{chamavisa}
T.~Chac{\'o}n~Rebollo, M.~G{\'o}mez~M{\'a}rmol, V.~Girault,
  I.~S{\'a}nchez~Mu{\~n}oz, A high order term-by-term stabilization solver for
  incompressible flow problems, IMA J. Numer. Anal. 33~(3) (2013) 974--1007.

\bibitem{Naveed15}
N.~Ahmed, T.~Chac\'on~Rebollo, V.~John, S.~Rubino, Analysis of a full
  space-time discretization of the {N}avier--{S}tokes equations by a local
  projection stabilization method, IMA J. Numer. Anal. 37~(3) (2017)
  1437--1467.

\bibitem{HFB86}
T.~J.~R. Hughes, L.~P. Franca, M.~Balestra, A new finite element formulation
  for computational fluid dynamics. {V}. {C}ircumventing the {B}abu\v
  ska-{B}rezzi condition: a stable {P}etrov-{G}alerkin formulation of the
  {S}tokes problem accommodating equal-order interpolations, Comput. Methods
  Appl. Mech. Engrg. 59~(1) (1986) 85--99.

\bibitem{FortiDede15}
D.~Forti, L.~Ded{\`e}, Semi-implicit {BDF} time discretization of the
  {N}avier-{S}tokes equations with {VMS}-{LES} modeling in a high performance
  computing framework, Comput. \& Fluids 117 (2015) 168--182.

\bibitem{Rubino18}
R.~Haferssas, P.~Jolivet, S.~Rubino, Efficient and scalable discretization of
  the {N}avier--{S}tokes equations with {LPS} modeling, Comput. Methods Appl.
  Mech. Engrg. 333 (2018) 371--394.

\bibitem{Bab71}
I.~Babu{\v{s}}ka, Error-bounds for finite element method, Numer. Math. 16
  (1970/1971) 322--333.

\bibitem{Bre74}
F.~Brezzi, On the existence, uniqueness and approximation of saddle-point
  problems arising from {L}agrangian multipliers, Rev. Fran\c caise Automat.
  Informat. Recherche Op\'erationnelle S\'er. Rouge 8~(R-2) (1974) 129--151.

\bibitem{HT74}
P.~Hood, C.~Taylor, {N}avier--{S}tokes equations using mixed interpolation, in:
  J.~T. Oden, R.~H. Gallagher, O.~C. Zienkiewicz, C.~Taylor (Eds.), Finite
  Element Methods in Flow Problems, University of Alabama in Huntsville Press,
  1974, pp. 121--132.

\bibitem{GWR04}
V.~Gravemeier, W.~A. Wall, E.~Ramm, A three-level finite element method for the
  instationary incompressible {N}avier-{S}tokes equations, Comput. Methods
  Appl. Mech. Engrg. 193~(15-16) (2004) 1323--1366.

\bibitem{GWR05}
V.~Gravemeier, W.~A. Wall, E.~Ramm, Large eddy simulation of turbulent
  incompressible flows by a three-level finite element method, Internat. J.
  Numer. Methods Fluids 48~(10) (2005) 1067--1099.

\bibitem{JKK08}
V.~John, S.~Kaya, A.~Kindl, Finite element error analysis for a
  projection-based variational multiscale method with nonlinear eddy viscosity,
  J. Math. Anal. Appl. 344~(2) (2008) 627--641.

\bibitem{Chacon98}
T.~Chac{{\'o}}n~Rebollo, A term by term stabilization algorithm for finite
  element solution of incompressible flow problems, Numer. Math. 79~(2) (1998)
  283--319.

\bibitem{hetobiska}
L.~He, L.~Tobiska, The two-level local projection stabilization as an enriched
  one-level approach, Adv. Comput. Math. 36~(4) (2012) 503--523.

\bibitem{beckerbraack2}
R.~Becker, M.~Braack, A two-level stabilization scheme for the
  {N}avier-{S}tokes equations, in: Numerical mathematics and advanced
  applications, Springer, Berlin, 2004, pp. 123--130.

\bibitem{SZ90}
L.~R. Scott, S.~Zhang, Finite element interpolation of nonsmooth functions
  satisfying boundary conditions, Math. Comp. 54~(190) (1990) 483--493.

\bibitem{Hecht12}
F.~Hecht, New development in freefem++, J. Numer. Math. 20~(3-4) (2012)
  251--265.

\bibitem{Badia12}
S.~Badia, On stabilized finite element methods based on the {S}cott-{Z}hang
  projector. {C}ircumventing the inf-sup condition for the {S}tokes problem,
  Comput. Methods Appl. Mech. Engrg. 247/248 (2012) 65--72.

\bibitem{BraackBurman06}
M.~Braack, E.~Burman, Local projection stabilization for the {O}seen problem
  and its interpretation as a variational multiscale method, SIAM J. Numer.
  Anal. 43~(6) (2006) 2544--2566.

\bibitem{KnoblochLube09}
P.~Knobloch, G.~Lube, Local projection stabilization for
  advection-diffusion-reaction problems: one-level vs. two-level approach,
  Appl. Numer. Math. 59~(12) (2009) 2891--2907.

\bibitem{Cellier91}
F.~E. Cellier, Continuous system modeling, Springer-Verlag, New York, 1991.

\bibitem{Codina02}
R.~Codina, Stabilized finite element approximation of transient incompressible
  flows using orthogonal subscales, Comput. Methods Appl. Mech. Engrg.
  191~(39-40) (2002) 4295--4321.

\bibitem{CodinaBlasco02}
R.~Codina, J.~Blasco, Analysis of a stabilized finite element approximation of
  the transient convection-diffusion-reaction equation using orthogonal
  subscales, Comput. Vis. Sci. 4~(3) (2002) 167--174.

\bibitem{Codina01}
R.~Codina, A stabilized finite element method for generalized stationary
  incompressible flows, Comput. Methods Appl. Mech. Engrg. 190~(20-21) (2001)
  2681--2706.

\bibitem{ParMooN}
U.~Wilbrandt, C.~Bartsch, N.~Ahmed, N.~Alia, F.~Anker, L.~Blank, A.~Caiazzo,
  S.~Ganesan, S.~Giere, G.~Matthies, R.~Meesala, A.~Shamim, J.~Venkatesan,
  V.~John, Parmoon---a modernized program package based on mapped finite
  elements, Computers \& Mathematics with Applications 74~(1) (2017) 74 -- 88.

\bibitem{SL17}
P.~W. Schroeder, G.~Lube, Divergence-free h(div)-fem for time-dependent
  incompressible flows with applications to high reynolds number vortex
  dynamics, J. Sci. Comput. 75~(2) (2018) 830--858.

\bibitem{Lesieur88}
M.~Lesieur, C.~Staquet, P.~Le~Roy, P.~Comte, The mixing layer and its coherence
  examined from the point of view of two-dimensional turbulence, J. Fluid Mech.
  192 (1988) 511--534.

\bibitem{Boersma97}
B.~J. Boersma, M.~N. Kooper, F.~T.~M. Nieuwstadt, P.~Wesseling, Local grid
  refinement in large-eddy simulations, J. Engrg. Math. 32~(2-3) (1997)
  161--175.

\bibitem{Burman07}
E.~Burman, Interior penalty variational multiscale method for the
  incompressible {N}avier-{S}tokes equation: monitoring artificial dissipation,
  Comput. Methods Appl. Mech. Engrg. 196~(41-44) (2007) 4045--4058.

\bibitem{John05}
V.~John, An assessment of two models for the subgrid scale tensor in the
  rational {LES} model, J. Comput. Appl. Math. 173~(1) (2005) 57--80.

\bibitem{Balaras01}
E.~Balaras, U.~Piomelli, J.~M. Wallace, Self-similar states in turbulent mixing
  layers, J. Fluid Mech. 446 (2001) 1--24.

\bibitem{LM96}
M.~Lesieur, O.~M{\'e}tais, New trends in large-eddy simulations of turbulence,
  in: Annual review of fluid mechanics, {V}ol.\ 28, Annual Reviews, Palo Alto,
  CA, 1996, pp. 45--82.

\bibitem{SF00}
K.~Schneider, M.~Farge, Numerical simulation of a mixing layer in an adaptive
  wavelet basis, Comptes Rendus de l'Acad{\'e}mie des Sciences - Series IIB -
  Mechanics-Physics-Astronomy 328~(3) (2000) 263 -- 269.

\bibitem{EAJ07}
E.~O{\~n}ate, A.~Valls, J.~Garc{\'\i}a, Computation of turbulent flows using a
  finite calculus--finite element formulation, International Journal for
  Numerical Methods in Fluids 54~(6‐8) (2007) 609--637.

\bibitem{DoeringGibbon95}
C.~R. Doering, J.~D. Gibbon, Applied analysis of the {N}avier-{S}tokes
  equations, Cambridge Texts in Applied Mathematics, Cambridge University
  Press, Cambridge, 1995.

\bibitem{AM16}
N.~Ahmed, G.~Matthies, Numerical study of {SUPG} and {LPS} methods combined
  with higher order variational time discretization schemes applied to
  time-dependent linear convection-diffusion-reaction equations, J. Sci.
  Comput. 67~(3) (2016) 988--1018.

\end{thebibliography}
\end{document}